\documentclass[10pt]{article}

\usepackage{amsfonts,latexsym,amstext}
\usepackage{amsmath}

\usepackage{amssymb}
\usepackage[english]{babel}
\usepackage[latin1]{inputenc}

\usepackage[usenames]{color}
\definecolor{red}{rgb}{1.0,0.0,0.0}

\definecolor{blu}{rgb}{0.0,0.0,1.0}

\definecolor{gre}{rgb}{0.03,0.50,0.03}

\definecolor{darkviolet}{rgb}{0.58, 0.0, 0.83}

\usepackage{enumerate}

\newtheorem{theorem}{Theorem}[section]
\newtheorem{lemma}[theorem]{Lemma}
\newtheorem{proposition}[theorem]{Proposition}
\newtheorem{definition}[theorem]{Definition}

\newtheorem{hypothesis}[theorem]{Hypothesis}

\newtheorem{remark}[theorem]{Remark}
\newtheorem{corollary}[theorem]{Corollary}

\numberwithin{equation}{section}

\newcommand{\myref}[1]{(\ref {#1})}

\def\qed{{\hfill\hbox{\enspace${ \square}$}} \smallskip}
\def\sqr#1#2{{\vcenter{\vbox{\hrule height .#2pt \hbox{\vrule
 width .#2pt height#1pt \kern#1pt \vrule
width .#2pt} \hrule height .#2pt}}}}
\def\square{\mathchoice\sqr54\sqr54\sqr{4.1}3\sqr{3.5}3}

\def\eps{\varepsilon}

\def\ds{\begin{displaystyle}}
\def\eds{\end{displaystyle}}
\def\dis{\displaystyle }
\def\<{\left\langle }
\def\>{\right\rangle }

\def\dim{\noindent \hbox{{\bf Proof.} }}

\def\v{\mathbf v}

\def\I{\mathbb I}
\def\R{\mathbb R}
\def\N{\mathbb N}

\def\E{\mathbb E}
\def\P{\mathbb P}

\def\calc{{\cal C}}
\def\cald{{\cal D}}

\def\calf{{\cal F}}

\def\calk{{\cal K}}
\def\call{{\cal L}}
\def\caln{{\cal N}}

\def\calu{{\cal U}}

\def\call{{\cal L}}
\def\cals{{\cal S}}
\def\calo{{\cal O}}

\def\Imm{{\operatorname{Im}}}

\def\to{\rightarrow}
\def\v{\mathbf v}

\begin{document}
\title{Lifting partial smoothing to solve HJB equations
\\
and stochastic control problems}
\date{}

 \author{Fausto Gozzi\footnote{Dipartimento di Economia e Finanza,
Universit\`a \emph{LUISS - Guido Carli} Roma, Italy;
e-mail: fgozzi@luiss.it}
\; and
Federica Masiero\footnote{Dipartimento di Matematica e Applicazioni, Universit\`a di Milano Bicocca, Milano, Italy;
e-mail: federica.masiero@unimib.it}
}

\maketitle

\vspace{-1truecm}

\begin{abstract}
We study a family of stochastic control problems arising in typical applications (such as boundary control and control of delay equations with delay in the control)
with the ultimate aim of finding solutions of the associated HJB equations,
regular enough to find optimal feedback controls.
These problems are difficult to treat since the underlying transition semigroups do not possess good smoothing properties nor the so-called ``structure condition'' which typically allows to apply the backward equations approach.
In the papers \cite{FGFM-I,FGFM-II} and, more recently, \cite{FGFM-III}
we studied such problems developing new partial smoothing techniques which allowed us to obtain the required regularity in the case when the cost functional is independent of the state variable.
This is a somehow strong restriction which is not verified in most applications.
\\
In this paper (which can be considered a continuation of the research of the above papers) we develop a new approach to overcome this restriction.
We extend the partial smoothing result to a wider class of functions which depend on the whole trajectory of the underlying semigroup and we use this as a key tool to improve our regularity result for the HJB equation.
The fact that such class depends on trajectories requires a nontrivial technical work
as we have to lift the original transition semigroup to a space of trajectories, defining a new ``high-level'' environment where our problems can be solved.
\end{abstract}

\textbf{Key words}:
Stochastic boundary control problems;
Stochastic Control of delay equation with delay in the control;
Second order Hamilton-Jacobi-Bellman equations in infinite dimension;
Smoothing properties of transition semigroups.

\smallskip \noindent

\textbf{AMS classification}:
93E20
47D07
49L20
35R15
93C23

\smallskip \noindent

\textbf{Acknowledgements}:
We thank Gianmario Tessitore for useful discussions.
Fausto Gozzi has been supported by the Italian
Ministry of University and Research (MIUR), in the framework of PRIN
project 2017FKHBA8 001 (The Time-Space Evolution of Economic Activities: Mathematical Models and Empirical Applications).
Federica Masiero has been supported
by the Gruppo Nazionale per l'Analisi Matematica, la Probabilit\`a e le loro Applicazioni (GNAMPA) of the Istituto Nazionale di Alta Matematica (INdAM).


\section{Introduction}
This paper is devoted to study a class of stochastic optimal control problems (with finite horizon $T<+\infty$) where the state space is an infinite dimensional Hilbert space $H$ and where the control operator is unbounded, i.e. its image goes out of the state space $H$, taking values in a larger space $\overline{H}$.
Such unboundedness arises commonly in applications, in particular in boundary control and in control of delay equations with delay in the control which are our driving examples and are treated in Section \ref{sec:verifica_ipotesi}.
{Our main goal is to prove existence and uniqueness of solutions of the associated HJB equations
regular enough to define optimal feedback controls
(Theorem \ref{esistenzaHJB-progr})
and to prove a verification theorem and the existence of optimal feedbacks (Theorems \ref{teorema controllo} and \ref{teo su controllo feedback}).}

We studied these types of problems in the recent paper \cite{FGFM-III}
(see also  \cite{FGFM-I,FGFM-II})
obtaining the required regularity in the case when the cost functional is independent of the state variable. This is a somehow strong restriction which is not verified in most applications (see e.g. \cite{CarmonaEtAl18,ChenWu20,EvansKendallTheodoru21,
GozziLeocata21,GMSJOTA,MouraFathy13,Sritharanbook}) and which comes from a difficult technical issue which we explain later in the paper, see, in particular, Remarks \ref{rm:prima-di sezHJB} and \ref{rm:dopoproofmain}.

In this paper (which can be considered a continuation of \cite{FGFM-III}) we develop a new approach to overcome such unpleasant restriction.
The core idea of \cite{FGFM-III} was the proof of partial smoothing results for a class of Ornstein-Uhlenbeck semigroups which includes our main examples.
Such results were proved for data of the form
$\phi(x)=\hat{\phi}(Px)$ for a suitable operator $P$
which is typically a projection on some finite dimensional subspace.

Our new approach is based, roughly speaking, on the extension of the partial smoothing results of
\cite{FGFM-III} to data $\phi$ depending, for every $\bar x\in \overline{H}$, on the
whole path $\{Pe^{tA}x\}_{t> 0}$, i.e.
$\phi(x)=\overline{\phi}\left(\{Pe^{tA}x\}_{t>0}\right)$, for all $x \in \overline{H}$,
where $\overline{\phi}:L^2(0,T;H)\to \R$ (see Subsection 3.2 for the precise setting).

We call this idea ``lifting" as we substitute the base space $H$ with the ``higher level" space
$L^2(0,T;H)$. We did not find any result of this type in the related literature.

Once such an extension is performed we use it as a key tool to improve our regularity result for the HJB equation.
Establishing such results requires a nontrivial technical work
as we have to lift the original transition semigroup to a space of trajectories, defining a new ``high-level'' environment where our problems can be attacked, and solving various technical issues which will be clarified along the way.

{We conclude this introduction briefly recalling the literature on our problems.
\\
First, we observe that some of the problems we treat here have been studied in the framework of viscosity solutions: we mention the papers \cite{deFeo} on problems with delay in the control and \cite{GozziRouySwiech06} on problems with boundary control. Such approach provides results on existence and uniqueness of solutions of HJB equations but do not provide any differentiability of them and no verification theorems.}
\\
{Second, the BSDE approach
(see \cite[Chapter 6]{FabbriGozziSwiech}) for this kind of problems has not been developed due to the lack of the so-called ``structure condition'' (the noise enter the system with the control) which, at the moment, seems crucial to use this approach.}
\\
{Third, the Pontryagin maximum principle has been applied to problem with delay in the control, see e.g. \cite{ChenWu10} where the case of pointwise delay is treated, and \cite{GuaMas} where a more general dependence on the past is considered. We notice that such results provide necessary conditions for optimality without proving existence of optimal control; we also notice that in the present paper we can also prove sufficient conditions and the existence of optimal feedbacks
(see our Theorem \ref{teo su controllo feedback}).
We also recall, the paper \cite{GuaMas13} which prove existence of optimal controls for boundary control problem in presence of boundary noise.
}

The plan of the paper is the following.
Section 2 is devoted to the partial smoothing results. Section 3 contains the setting of our abstract control problem
and our main result on the HJB equation.
{Section 4 presents our two main examples showing that they satisfy the required assumptions.
Section 5 deals with the verification theorem and the existence of optimal feedbacks. It is given after the examples section as here the technique to treat the two examples is different.}
The Appendix A is devoted to recall, for the reader's convenience, some preliminary material on basic notations and generalized derivatives {which can be found, e.g., in \cite{FGFM-III}.} {We suggest the reader to go to Appendix A to see the notations and the basic definitions before starting the reading of the paper.}

\section{Partial smoothing for OU semigroups: extensions}
\label{sec:partsmooth-abstr-setting}

In this section we extend the ``partial smoothing'' properties of the Ornstein-Uhlenbeck semigroup (which we call $R_t$, for $t\ge 0$) proved in \cite[Section 5]{FGFM-III}
to spaces of functions depending on the trajectories of the underlying semigroup $e^{tA}$.


We divide this section in three parts:
\vspace{-0.2truecm}
\begin{itemize}
  \item Subsection \ref{sub:OUknown} where we recall the framework and state the main assumptions and the results of \cite[Section 4]{FGFM-III} which are the starting point of our new results;
\vspace{-0.2truecm}
  \item Subsection \ref{sub:OUnew}
  where we provide our extended partial smoothing result.
\vspace{-0.2truecm}
  \item Subsection \ref{sec:convpartsmooth-abstr-setting}
where we provide our extended partial smoothing result
for convolution integrals: this requires an additional effort
which is crucial to treat our HJB equation in the next Section.
  \end{itemize}

\subsection{Setting and a smoothing result from \cite{FGFM-III}}
\label{sub:OUknown}

The following basic assumption holds throughout this section.
\begin{hypothesis}\label{ip-sde-common}
\begin{enumerate}[(i)]
\item[]
\item
Let $H$, $K$ and $\Xi$ be three real separable Hilbert
spaces\footnote{These will be usually the state space,
the control space and the noise space, respectively.}.
  \item
Let $(\Omega, \calf,(\calf_t)_{t\geq 0}, \P)$ be a filtered probability space
satisfying the usual conditions and let $W$ be an
$(\Omega, \calf,(\calf_t)_{t\geq 0}, \P)$-cylindrical Wiener process in $\Xi$ where $(\calf_t)_{t\geq 0}$ is the augmented filtration generated by $W$.
  \item
Let $A:D(A)\subseteq H \to H$ be the generator of a strongly continuous semigroup
$e^{tA},\, t\geq 0$ in $H$,
  \item
Let $G\in\call(\Xi,H)$ be such that the selfadjoint operator
$Q_t=\int_0^t e^{sA}GG^*e^{sA^*}\,ds$
is trace class. Set $Q=GG^* \in \call(H)$.
\end{enumerate}
\end{hypothesis}
Let $Z(\cdot;x)$ be the Ornstein-Uhlenbeck process
{which solves the following SDE in $H$:
\vspace{-0.3truecm}
\begin{equation}\label{ornstein-gen}
dZ(t)=AZ(t)dt+GdW(t), \qquad X(0)=x.
\vspace{-0.3truecm}
\end{equation}
The process $Z(\cdot;x)$ can be written in mild form as:}
\vspace{-0.3truecm}
\begin{equation}
Z(t;x)  =e^{tA}x +\int_0^te^{(t-s)A}GdW(s)
,\text{ \ \ \ }t\ge 0. \\
  \label{ornstein-mild-gen}
\vspace{-0.2truecm}
\end{equation}
$Z$ is a Gaussian process, namely for every $t>0$, the law of
$Z(t)$ is $\caln (e^{tA}x,Q_t)$, the Gaussian measure with mean $e^{tA}x$ and
covariance operator $Q_t$ defined in Hypothesis
\ref{ip-sde-common}-(iv) .
The convolution $\int_0^te^{(t-s)A}GdW_s$ has law
$\caln (0,Q_t)$ and will be sometimes denoted by $W_A(t)$.
The associated Ornstein-Uhlenbeck transition semigroup $R_t$ is defined by setting, for every $\psi\in B_b(H)$ and $x\in H$,
\vspace{-0.3truecm}
\begin{equation}
 \label{ornstein-sem-gen}
R_t[\psi](x)=\E \psi(Z(t;x))
=\int_H \psi(z+e^{tA}x)\caln(0,Q_t)(dz).
\vspace{-0.2truecm}
\end{equation}
As in \cite[Section 4]{FGFM-III} we assume the following on the control operator $C$ and on the auxiliary operator $P$.

\begin{hypothesis}\label{ip:PC}
\begin{enumerate}[(i)]
\item[]\vspace{-0.2truecm}
\item
Let $\overline H$ be a real Banach space such that
$H\subseteq \overline{H}$ with continuous and dense inclusion
and that the semigroup $e^{tA}$ admits an extension $\overline{e^{tA}}:\overline{H}\to \overline{H}$ which is still a $C_0$ semigroup.\vspace{-0.2truecm}
\item
Let $C\in \call(K, \overline H)$.
\vspace{-0.2truecm}
\item Let $P:H \to H$ be linear and continuous.
Assume that, for every $t>0$ the operator $Pe^{tA}:H \to H$ can be extended to a continuous linear operator $\overline H\to H$, which will be denoted by $\overline{Pe^{tA}}$. With this notation the operator $\overline{Pe^{tA}}C:K\to H$ is well defined and continuous.
\end{enumerate}
\end{hypothesis}\vspace{-0.2truecm}

In \cite[Remark 5.4]{FGFM-III} it is explained that
such assumptions are verified in the key examples of
Dirichlet boundary control and control of delay SDEs
with delay in the control.
%
%
We now recall the spaces where the initial data should belong, for the smoothing result of \cite[Proposition 4.9]{FGFM-III} to hold.

\begin{definition}\label{df:spaziphi1}
We call $B_b^P(H)$ (respectively $C_b^P(H)$) the set of  Borel measurable functions $\phi:H\to \R$ for which there exists $\bar\phi : \operatorname{Im}(P)\to \R$
bounded and Borel measurable and (respectively continuous)\footnote{Here we endow $\operatorname Im P\subseteq H$ with the topology inherited by $H$.} such that
$\phi(x)=\bar\phi(Px)$
$\forall x\in H$.
{Similarly, for $T>0$, we denote by
$B_b^P([0,T]\times H)$ (respectively $C_b^P([0,T]\times H)$) the set of Borel measurable functions $\psi:[0,T]\times H\to \R$ for which there exists
$\bar\psi : [0,T]\times \operatorname{Im}(P)\to \R$
bounded and Borel measurable {(respectively continuous), and} such that
\vspace{-0.2truecm}
\begin{equation}\label{fi-gen-allargatatime}
\psi(t,x)=\bar\psi(t,Px) \quad
\forall (t,x)\in [0,T]\times H.
\end{equation} }
\end{definition}


\begin{remark}\label{rm:ipPCunifcont}
As noted in \cite[Remark 5.6]{FGFM-III},
$B_b^P(H)$ and $C_b^P(H)$ are linear subspaces of
$B_b(H)$ (respectively $C_b(H)$). Similarly for the time dependent case.
When the image of $P$ is closed (as it is in our examples see Section \ref{SSE:HEATEQUATION}, formula \eqref{P-n-gen} and Section \ref{SSE:DELAYEQUATION}), we can identify the space $B_b^P(H)$ with $B_b(\operatorname Im P)$ (and the same for $C_b^P(H)$ and for the time dependent spaces).
In particular, in our examples,
$\Imm P$ is finite dimensional. Hence, for some $n\in\N$
we get $B^{P}_b(H) \sim B_b(\R^n)$, $C^{P}_b(H) \sim C_b(\R^n)$, $B^{P}_b([0,T]\times H) \sim B_b([0,T]\times \R^n)$, $C^{P}_b([0,T]\times H) \sim C_b([0,T]\times \R^n)$.
\end{remark}
\vspace{-0.2truecm}

We now recall the partial smoothing result proved in \cite[Proposition 5.9]{FGFM-III}. We need first the following controllability-like assumption, see \cite[Hypothesis 5.7]{FGFM-III}.

\vspace{-0.2truecm}

\begin{hypothesis}\label{ip:NC}
\begin{itemize}
\item[]
\item [(i)]
We have
$\operatorname{Im}\overline{Pe^{tA}}C\subseteq
\operatorname{Im}(P Q_t P^*)^{1/2}$, $\forall t>0$;
Consequently, by the Closed Graph Theorem, the operator\vspace{-0.2truecm}
$$
\Lambda^{P,C}(t):K\to H, \qquad
\Lambda^{P,C}(t)k:=(P Q_t P^*)^{-1/2}\overline{Pe^{tA}}Ck
\quad \forall k \in K,
$$
is well defined and bounded for all $t>0$.\vspace{-0.2truecm}
\item [(ii)]
There exists $\kappa_1>0$ and $\gamma \in (0,1)$ such that
$\|\Lambda^{P,C}(t)\|_{\call(K,H)} \le
\kappa_1 (t^{-\gamma}\vee 1)
\qquad
\forall t \in (0,T]$.
\end{itemize}
\end{hypothesis}

Note that in \cite[Remark 5.8]{FGFM-III} it is shown that the above Hypothesis \ref{ip:NC} holds in our examples, see Section \ref{sec:verifica_ipotesi}.

\begin{proposition}\label{prop:partsmooth}
Let Hypotheses \ref{ip-sde-common}, \ref{ip:PC} and \ref{ip:NC}-(i) hold true.
Then the semigroup $R_t,\,t>0$ maps functions $\phi\in B_b^P(H)$
into functions which are $C$-Fr\'echet differentiable in $\overline{H}$, and the $C$-derivative is given, for all
$\bar x \in \overline{H}$, by\vspace{-0.2truecm}
\begin{align}\label{eq:formulader-gen-Pold}
 \nabla^C(R_{t}[\phi])(\bar x)k &=\int_{H}\bar\phi\left(z_1+\overline{Pe^{tA}}\bar x\right)
\<\Lambda^{P,C}(t) k,
(PQ_tP^*)^{-1/2}z_1\>_H\caln(0,PQ_tP^*)(dz_1)
\\ \nonumber
&=
\E\left[\bar\phi\left(PX(t;\bar x)\right)
\<\Lambda^{P,C}(t) k,
(PQ_tP^*)^{-1/2}PW_A(t)\>_{H}
\right]
\end{align}
Moreover, for any $\phi\in B^P_b(H)$ and any $\bar x\in \overline{H}$, $k\in K$,\vspace{-0.2truecm}
\begin{equation}\label{norm-Cderold}
\vert \<\nabla^C R_t[\phi](\bar x), k\>\vert \leq
\Vert \Lambda^{P,C}(t) \Vert_{\call(K, H)} \Vert \phi\Vert_\infty \vert k\vert.
\end{equation}
Furthermore, if $\phi\in C^P_b(H)$, then $\nabla^C R_t[\phi]\in C((0,T]\times \overline{H};K)$.
Finally, if also Hypothesis \ref{ip:NC}-(ii) holds,
then the map $(t,\bar x) \to R_t[\phi](\bar x)$ belongs to
$C^{0,1,C}_\gamma([0,T]\times \overline{H})$.
\end{proposition}


\begin{remark}\label{rm:partsmooth-secondold}
As noted in  \cite[Remark 5.10]{FGFM-III}
using the ideas of Proposition 4.5 in \cite{FGFM-I}
it is possible to prove that, if $\phi$ is more regular
(i.e. $\phi\in C^1_b(H)\cap C_b^P(H)$,
also $\nabla^C R_t[\phi]$ has more regularity, i.e.
$\nabla\nabla^{C}R_{t}\left[\phi\right]$,
$\nabla^{C}\nabla R_{t}\left[\phi\right]$ exist, coincide, and
satisfy suitable formulae and estimates.
\end{remark}
\vspace{-0.6truecm}

\subsection{Partial smoothing for functions of paths}
\label{sub:OUnew}

\subsubsection{Setting and spaces}

We now improve the above result to cover a more general class of data $\phi$.
First we introduce a class of spaces which is suitable for our needs.
To do this we need the following preliminary lemma.
\begin{lemma}
Let Hypotheses \ref{ip-sde-common} and \ref{ip:PC}
hold true.
For every $\bar x\in \overline{H}$ and
$0<s\le t$ we have
\vspace{-0.2truecm}
\begin{equation}\label{eq:Petasemigroup}
\overline{Pe^{tA}}\bar x=
\overline{Pe^{sA}}\cdot \overline{e^{(t-s)A}}\bar x.
\vspace{-0.2truecm}
\end{equation}
Moreover the map
$(0,+\infty)\to H$,
$t \mapsto \overline{Pe^{tA}}\bar x$,
is continuous.
\end{lemma}

\dim
Let first $\bar x\in H$. Then, clearly
\eqref{eq:Petasemigroup} holds.
Now, let $\bar x \in \overline{H}$.
Then by density (Hypothesis \ref{ip:PC}-(i))
there exists a sequence $\{x_n\}_{n\in \N}\subset H$
such that $\lim_n x_n=\bar x$ in $\overline{H}$.
Take $0<s\le t$. For every $n\in \N$ we have
${Pe^{tA}}x_n={Pe^{sA}}{e^{(t-s)A}}x_n$.
Now send $n$ to infinity.
By Hypothesis \ref{ip:PC}-(i)
$\lim_n e^{(t-s)A}x_n=\overline{e^{(t-s)A}}\bar x$ in $\overline{H}$.
Moreover using that $\overline{Pe^{tA}}:\overline{H}\to H$ is a continuous linear operator we get
$\lim_n {Pe^{tA}}(e^{(t-s)A}x_n)=
\overline{Pe^{tA}}(\overline{e^{(t-s)A}}\bar x)$
which proves \eqref{eq:Petasemigroup}.
\\
We prove continuity. Let $\bar x \in \overline{H}$ and let
Take a sequence $\{t_k\}_{k\in \N}\subset (0,+\infty)$ converging to $t_\infty>0$.
Without loss of generality we can assume that
$\hat{t}:=\inf_k t_k$ is strictly positive.
Hence, by \eqref{eq:Petasemigroup} we have, for all $k\in \N$
$\overline{Pe^{t_kA}}\bar x=
\overline{Pe^{\hat tA}}\cdot
\overline{Pe^{(t_k-\hat t)A}}\bar x$.
By Hypothesis \ref{ip:PC}-(i) (namely the fact that
the extended semigroup $\overline{e^{tA}}$ is strongly continuous in $\overline{H}$) we get
$\lim_k \overline{e^{(t_k-\hat t)A}}\bar x
=
\overline{e^{(t_\infty-\hat t)A}}\bar x$
in $\overline{H}$.
This provides continuity.
\hfill\qed

\begin{definition}
\label{df:CTCPT}
Take $T\in (0,+\infty]$. For simplicity we write $(0,T]$ for $(0,T]\cap \R$.
Define the set of paths\vspace{-0.2truecm}
$$
\calc^P_A((0,T]; H):= \Big\{
f\in
C ((0,T]; H) \hbox{ such that } \exists\, \bar x\in \overline{H}:\;
f(t)=\overline{Pe^{tA}}\bar x,\; \forall t \in (0,T]\Big\}.\vspace{-0.2truecm}
$$
For every
$\bar x\in \overline{H}$
we call
$y^P_{\bar x}(\cdot)$
the path given by
$y^P_{\bar x}(t)=\overline{Pe^{tA}}{\bar x}$ for all $t \in (0,+\infty)$).
\end{definition}

%

\vspace{-0.2truecm}

In the next lemma we introduce and study a natural operator from $\overline{H}$ to $\calc^P_A((0,T]; H)$.

\vspace{-0.2truecm}

\begin{lemma}
\label{lm:UpsilonPT}
Let Hypotheses \ref{ip-sde-common} and \ref{ip:PC}
hold true.
Let $T\in (0,+\infty]$ and define the map
$\Upsilon^P_T:\overline{H} \to \calc^P_A((0,T]; {H})$
as $\Upsilon^P_T ({\bar x})=y^P_{\bar x}|_{[0,T]}$,
for all ${\bar x} \in \overline{H}$.
$\Upsilon^P_T$ is surjective but not necessarily injective\footnote{Such lack of injectivity arises, e.g., in our examples.}.
Moreover, we have, endowing $C((0,T]; H)$ with the topology of uniform convergence on compact subsets, that
$\bar x_n\to \bar x$ in $\overline{H}$ $\Longrightarrow$ $y^P_{\bar x_n}|_{[0,T]}\to y^P_{\bar x}|_{[0,T]}$
in $C ((0,T]; H)$.
The {converse} is not true in general.
Hence $\Upsilon^P_T$ is continuous if we endow $\calc^P_A((0,T]; {H})$ with the topology inherited by
$C ((0,T]; H)$.
\end{lemma}

\dim
The fact that $\Upsilon^P_T$ is surjective is a direct consequence of the definition of $\calc^P_A((0,T];{H})$, however it is not injective. As an example (taken in a simple case of the Subsection \ref{SSE:HEATEQUATION}),
let $A$ be the Laplace operator with Dirichlet boundary conditions on $H=L^2([0,\pi])$ and let $P:\overline{H}\rightarrow H,\; Px=\<x,v\>v,\, \forall x \in \overline{H}$, where $v \in H$ is an eigenvector of $A$.
In this case for $x\in H$ such that $\<x,v\>=0$ it is immediate to see that
$\overline{Pe^{tA}}x=Pe^{tA}x=e^{tA}Px=0$.
Note that, for $T<+\infty$, ${\bar x}_n\to {\bar x}$ in $\overline{H}$ implies (see Hypothesis \ref{ip:PC}-(ii) and (iv)) that $y^P_{{\bar x}_n}|_{[0,T]}\to y^P_{\bar x}|_{[0,T]}$ in $C_b ([0,T]; H)$
(hence $\Upsilon^P_T$ is continuous).
This holds also when $T=+\infty$ when $A$ is of negative type.
However, differently from above, the {converse} is not true in general, e.g. in the case of the example of Subsection \ref{SSE:HEATEQUATION}: with $P$ as before, we can take a sequence $(x_n)_n\subseteq H$ not converging to $0$ (for example we can consider as $(x_n)_n$ to be the eigenvectors of $A$ othogonal to $v$) such that
$Pe^{tA}x_n=e^{tA}Px_n\equiv 0$ for all $n\in \N$.
\hfill\qed

\begin{remark}
Without explicit notice we will take,
on $\mathcal{C}^P_A((0,T]; {H})$, the topology inherited by
$C ((0,T]; {H})$.
\end{remark}


We consider also the space of functions that are square integrable on $[0,+\infty)$ when multiplied by a suitable weight $e^{-\rho t}, \,\rho>0$:\vspace{-0.2truecm}
$$
L^2_\rho(0,\infty;H):=\left\lbrace f:[0,+\infty)\to H: \quad \|f\|_{L_\rho^2}:=
\displaystyle \int_0^{+\infty}e^{-2\rho t} \vert f(t)\vert^2_H\,dt<+\infty \right\rbrace.\vspace{-0.2truecm}
$$
This is a Hilbert space with the inner product
$\<f,g\>_{L_\rho^2}:= \int_0^{+\infty}
e^{-2\rho t} \<f(t),g(t)\>_H\,dt$.
For simplicity, if no confusion is possible we will write $L^2_\rho$ instead of $L^2_\rho(0,\infty;H)$ and $L^2_T$ for $L^2(0,T;H)$.

Now we add an assumption which allows to prove the analogous of Lemma \ref{lm:UpsilonPT} in the spaces $L^2$.
\begin{hypothesis}
\label{hp:L2}
For every finite $T>0$ and $\bar x \in \overline{H}$ there exist $C_T>0$ and $\eta \in [0,1/2)$ such that
$|\overline{Pe^{tA}}\bar x|_H\le C_T t^{-\eta}|\bar x|_{\overline{H}}$, $\forall (t,\bar{x})\in(0,T]\times \overline{H}$,
{and, as a consequence, the map
$(0,T]\to H$, $t \to \overline{Pe^{tA}}\bar x$,
belongs to $L^2(0,T;H)$.}
\end{hypothesis}

\begin{lemma}
\label{lm:UpsilonL2}
Let Hypotheses \ref{ip-sde-common}, \ref{ip:PC}
and \ref{hp:L2} hold true.
Let $T<+\infty$; then for all
$x\in H$ and
$\bar x\in \overline{H}$, we have
$e^{\cdot A}x|_{[0,T]},
y^P_{\bar{x}}|_{[0,T]}\in L^2(0,T; H)$.
Hence,
$\calc^P_A((0,T]; \overline{H})$
can be seen as a linear subspace of $L^2(0,T; H)$ with continuous embedding given by
$\Upsilon^P_T$.
\\
If $x_n\to x$ in $H$ this implies that
$e^{\cdot A}{x_n}|_{[0,T]}\to
e^{\cdot A}{x}|_{[0,T]}$ in $L^2 (0,T; H)$
but the {converse}  is not true in general.
\\
If $\bar x_n\to \bar x$ in $\overline{H}$ this implies that $y^P_{\bar x_n}|_{[0,T]}\to y^P_{\bar x}|_{[0,T]}$ in $L^2 (0,T; H)$
but the {converse} is not true in general.
\\
Let now $T=+\infty$ and $A$ of type $\omega${\footnote{{We say that $A$ is of type $\omega$ if there exists some $M\geq 1$ such that $\Vert e^{At}\Vert \leq  Me^{\omega t}$ for all $t\ge 0$.}}.} and consider $\rho>\omega$. then for all $x\in H$ and $\bar x\in \overline{H}$, we have
$e^{\cdot A}x, y^P_{\bar x}\in L_\rho^2(0,+\infty; H)$.
Hence,
$\calc^P_A((0,+\infty); \overline{H})$
can be seen as linear subspace of $L_\rho^2(0,\infty; H)$ with continuous embedding given by
$\Upsilon^P_{\infty}$.
\\
If $x_n\to x$ in $H$ this implies that
$e^{\cdot A}x_n\to e^{\cdot A}x$ in $L_\rho^2 (0,\infty; H)$ but the {converse}  is not true in general.
\\
If $\bar x_n\to \bar x$ in $\overline{H}$ this implies that $y^P_{\bar x_n}\to y^P_{\bar x}$ in $L_\rho^2 (0,\infty; H)$ but the {converse} is not true in general.
\end{lemma}

\dim When $T<+\infty$ the assertions follow immediately due to the boundedness of
$e^{\cdot A}x|_{[0,T]}, y^P_{\bar x}|_{[0,T]}$; the stated
continuities can be proved by dominated convergence theorem and also taking into account Hypothesis \ref{hp:L2}.
\\
We now prove that $e^{\cdot A}x_n|_{[0,T]}\to e^{\cdot A}x|_{[0,T]}$ in $L^2 (0,T; H)$ does not imply $x_n\to x$ in $H$.
As an example consider the semigroup $e^{tA^*}$ defined in \eqref{semigroupadjoint}. If we take $x_n=(0,a_n)$ where
$a_n(\xi)=n^\alpha\I_{[-d,-d+1/n]}$ we see that, $|x_n|^2_H=n^{2\alpha-1}$, while, for $T\ge d$,\vspace{-0.2truecm}
$$
\int_{0}^{T} |e^{tA^*}x_n|^2dt=
\int_{0}^{d} |e^{tA^*}x_n|^2dt=
\int_{0}^{d}\int_{-d}^{-t} a_n(s+t)^2 ds dt
$$
$$
=\int_{0}^{d}\int_{-d+t}^{0} n^{2\alpha}\I_{[-d,-d+1/n]}(r) dr dt=
n^{2\alpha}
\int_{0}^{1/n}\left(\frac{1}{n}-t\right) dt=
\frac{1}{2} n^{2\alpha-2}.
$$
Hence, for $\alpha \in (1/2,1)$ we get
$e^{\cdot A}{x_n}|_{[0,T]}\to e^{\cdot A}0|_{[0,T]}$ in $L^2 (0,T; H)$
but $x_n\not\to 0$ in $H$.

We now prove that $y^P_{\bar x_n}|_{[0,T]}\to y^P_{\bar x}|_{[0,T]}$ in $L^2 (0,T; H)$ does not imply $x_n\to x$ in $H$.
Indeed as in Lemma \ref{lm:UpsilonPT} for the example of Subsection \ref{SSE:HEATEQUATION} we can consider $P:\overline{H}\rightarrow H,\; Px=\<x,v\>v,\, \forall x \in \overline{H}$, where $v \in H$ is an eigenvector of $A$ and
 we can take a sequence of orthonormal eigenvectors of $A$, orthogonal to $v$. We have that $x_n\nrightarrow 0$ in $H$ but $\int_0^T \vert Pe^{tA}x_n\vert^2\, dt\equiv 0$.

When $T=+\infty $ the results are proved as in the case $T<+\infty$ since $\rho>\omega$.
\hfill\qed



{
\begin{definition}
\label{df:SP}
We introduce, for $T>T_0\ge 0$ and $\eta \in [0,1)$, the following set of functions (here $Z$ is a Hilbert space):
\begin{multline}
\label{eq:defSbar}
\cals^P_\infty( \overline{H})
 :=\Big\{\phi: \overline{H}\to \R \hbox{ s.t. }
\exists \,\hat \phi:\calc^P_A((0,+\infty); H) \to\R :
\hbox{ $\hat \phi$ bounded,
Borel meas. and }
\phi(\bar x)=\hat \phi(y^P_{\bar x}), \;\forall \bar x \in \overline{H}
\Big\}.
\end{multline}
\vspace{-0.7truecm}
\begin{multline}
\label{eq:defS}
\cals^P_\eta((T_0,T]\times \overline{H};Z)
 :=\Big\{f:(T_0,T]\times \overline{H}\to Z \hbox{ such that }
\exists \,\hat f:(T_0,T]\times \calc^P_A((0,+\infty); H) \to Z :
\\
\hbox{$\hat f$ is Borel meas., }
\hbox{the map $(t,\bar x)\mapsto (t-T_0)^\eta\hat f(t,\bar x)$ is bounded, and }
f(t,\bar x)=\hat f(t,y^P_{\bar x}),\;
\forall (t,\bar x)\in (T_0,T]\times \overline{H}
\Big\}
\end{multline}
\vspace{-0.7truecm}
\begin{multline}
\label{eq:defSprogr}
\cals^P_{\eta,prog}((T_0,T]\times \overline{H};Z)
 :=\{f\in\cals^P_{\eta}((T_0,T]\times \overline{H};Z)
 \hbox{ such that the associated }
\hat f:(T_0,T]\times \calc^P_A((0,+\infty); H) \to Z :
\\
\hbox{ satisfies also }
f(t,\bar x)=\hat f(t,y^P_{\bar x}(\cdot \wedge t)),\;
\forall (t,\bar x)\in (T_0,T]\times \overline{H}
\}
\end{multline}
When $\eta =0$ we omit it simply writing
$\cals^P((T_0,T]\times \overline{H};Z)$ or
$\cals^P_{prog}((T_0,T]\times \overline{H};Z)$.
Finally, as usual, we omit the arrival space in \eqref{eq:defS} when $Z=\R$.
\end{definition}
}

We will need the following properties of $\mathcal{S}^P$
and the adjoints of $\Upsilon^P_\infty,\Upsilon^P_T$.

\begin{proposition}
  \label{pr:propcalS}
Assume Hypotheses \ref{ip-sde-common}, \ref{ip:PC} and \ref{hp:L2}.
We have the following.
\begin{enumerate}
  \item [(i)]
We have $\cals^P_{\infty}(\overline{H})
\subseteq B_b(\overline{H})$. Moreover the restriction to $H$ of a function in $\cals^P_{\infty}(\overline{H})$ belongs to
$B_b(H)$.
  \item [(ii)]
Let $\phi\in \cals^P_\infty(\overline{H})$ and $\hat\phi$
be the function given in the definition of $\cals^P_\infty(\overline{H})$: we have
$\phi=\hat \phi \circ \Upsilon^P_\infty$.
  \item [(iii)]
Assume that $P$ can be extended to a continuous linear operator $\overline{P}:\overline{H}\to H$ such that
$Im \overline{P}=Im P$. Then
we have
$B_b^P(H)\subseteq\cals^P_\infty(\overline{H})$ and
$B_b^P([0,T]\times H)\subseteq \cals^P_0([0,T]\times \overline{H})$ .
Moreover, if also $P$ commutes with $A$ (as in the case of Subsection \ref{SSE:HEATEQUATION}),
then $B_b^P(H)=\cals^P_\infty(\overline{H})$
and
$B_b^P([0,T]\times H)= \cals^P_0([0,T]\times \overline{H})$.
  \item [(iv)]
The adjoint operator
$(\Upsilon^P_\infty)^*:L^2_\rho(0,\infty;H)
\to\overline{H}'\subseteq H$
is given by
$(\Upsilon^P_\infty)^*z(\cdot)=\int_{0}^{+\infty}e^{-\rho s}e^{sA^*}P^*z(s)ds$.
\\
Similarly, for $T<+\infty$,
$(\Upsilon^P_T)^*:L^2(0,T;H)
\to\overline{H}'\subseteq H$
is given by
$(\Upsilon^P_T)^*z(\cdot)=\int_{0}^{T}e^{sA^*}P^*z(s)ds$.
\end{enumerate}
\end{proposition}
\dim
The first statement of (i) is a simple consequence of the definitions \eqref{eq:defSbar}-\eqref{eq:defS}.
The second statement of (i) follows since, by Hypothesis \ref{ip:PC}-(i) $H$ is a measurable subset of $\overline{H}$.
\\
The statement (ii) again follows immediately by the definitions
\eqref{eq:defS}-\eqref{eq:defSbar}.
\\
Concerning statement (iii), under such further assumption on $P$, we see that every element of $\calc^P_A((0,+\infty); H)$
admits a limit for $t\to 0^+$ which is $\overline{P}\bar x$.
Hence in this case the map
$\delta_0:\calc^P_A((0,+\infty); H) \to Im P\subseteq H$
bringing an element of $\calc^P_A((0,+\infty); H)$ to its limit at $0^+$ is measurable and continuous.
Moreover every $\phi \in B_b^P(H)$ is immediately extended to $\overline{H}$ calling $\phi(\bar x)=\bar \phi (\overline{P}\bar x)$.
The above implies that, for $\bar x\in \overline{H}$,
$\phi(\bar x)=\bar \phi(\overline{P}\bar x)=\bar\phi (\delta_0 y^P_x)$.
Then, setting $\hat \phi(z):=\bar\phi (\delta_0 z)$, we see that $\phi\in\cals^P_\infty(\overline H)$.
This proves that $B_b^P(H)\subseteq\cals^P_\infty(\overline H)$.\footnote{On the other hand, if $\phi(x)=Pe^{TA}x$ and $P$ does not commute with $A$ (as in the case of Subsection \ref{SSE:DELAYEQUATION}),
such $\phi$ belongs to $\cals^P_\infty(\overline H)$ but it could not belong to
$B_b^P(H)$, and this is the case Subsection \ref{SSE:DELAYEQUATION}, where
$B_b^P(H)\subsetneqq\cals^P_\infty(\overline H)$.}
\\
Finally, let $P$ commute with $A$, hence with $e^{tA}$ for all $t>0$. Then also $\overline{P}$ commutes with $\overline{e^{tA}}$ for all $t>0$.
Let $\phi \in\cals^P_\infty(\overline H)$ and let $\hat \phi$ be the associated map along \eqref{eq:defSbar}.
Then, for $\bar x\in \overline{H}$ we get
\vspace{-0.2truecm}
$$
\phi(\bar x)=\hat \phi(y^P_{\bar x})=
\hat \phi\left(\overline{P}\; \overline{e^{\cdot A}}\bar x\right)
=
\hat \phi\left(\overline{e^{\cdot A}}\;\overline{P} \bar x\right).
=
\hat \phi\left(e^{\cdot A}\; \overline{P} \bar x\right).
\vspace{-0.2truecm}
$$
Hence, calling
$\bar \phi(z)=\hat \phi\left(e^{\cdot A}z\right)$,
we get $\phi \in B_b^P(H)$.
The proof of (iii) in the time dependent case is analogous
\\
The last point is immediate, using Hypothesis \ref{hp:L2} and the definitions of $\Upsilon^P_\infty$ and $\Upsilon^P_T$ respectively.
\hfill\qed

We conclude this subsection introducing the spaces which will be crucial to perform our fixed point argument to solve our HJB equations.

\vspace{-0.2truecm}

\begin{definition}
\label{df:spaziSeta}
Let $T>T_0\ge 0$, $\eta \in [0,1)$, and consider the following subset of $\cals^P((T_0,T]\times\overline{H})$:
\vspace{-0.2truecm}
\begin{align}
\label{eq:Setadef}
\cals^{1,P}_{\eta}((T_0,T]\times\overline{H}) :=\Big\{&f:[T_0,T]\times \overline{H}\to \R :f\in \cals^P((T_0,T]\times\overline{H})\cap C^{0,1,C}_{\eta}\left([T_0,T]\times \overline{H}\right)
\\
&\nonumber
\hbox{ such that the map } (t,\bar x) \mapsto (t-T_0)^\eta\nabla^C f(t,\bar x)\in \cals^P((T_0,T]\times\overline{H};K)\Big\}.
\vspace{-0.2truecm}
\end{align}
Hence, if a function $f\in \cals^{1,P}_{\eta}((T_0,T]\times\overline{H})$, there exist two functions:
\vspace{-0.2truecm}
$$
\hat f:[T_0,T]\times \calc^P_A((0,+\infty);H) \to \R,
\hbox{such that } f(t,\bar x)=
\hat f(t,y^P_{\bar x}),\;\forall
(t,\bar x)\in (0,T]\times \overline{H}
\vspace{-0.2truecm}
$$
and
\vspace{-0.2truecm}
$$
\hat{\hat f}:[T_0,T]\times \calc^P_A((0,+\infty);H) \to K,
\hbox{such that }
(t-T_0)^\eta\nabla^C f(t,\bar x)=
\hat{\hat f}(t,y^P_{\bar x}),\;
\forall (t,\bar x)\in (T_0,T]\times \overline{H}.
\vspace{-0.2truecm}
$$
When $\eta =0$ we omit it simply writing
$\cals^{1,P}((T_0,T]\times \overline{H})$.
\end{definition}

\begin{definition}
\label{df:spaziSeta-progr}
Let $T> 0$ and consider the following subset of $\cals^P_{prog}((T_0,T]\times\overline{H})$:\vspace{-0.2truecm}
\vspace{-0.2truecm}
\begin{align}
\label{eq:Setadef-prog}
\cals^{1,P}_{\eta,prog} ((T_0,T]\times\overline{H}):=\Big\{&f:[T_0,T]\times \overline{H}\to \R :f\in \cals^P((T_0,T]\times\overline{H})\cap C^{0,1,C}_{\eta}\left([T_0,T]\times \overline{H}\right)
\\
&\nonumber
\hbox{ such that the map } (t,\bar x) \mapsto (t-T_0)^\eta\nabla^C f(t,\bar x)\in \cals^P_{prog}((T_0,T]\times\overline{H};K)\Big\}.
\vspace{-0.2truecm}
\end{align}
Hence, if a function $f\in \cals^{1,P}_{\eta,prog}((T_0,+\infty)\times\overline{H})$, there exist two functions:\vspace{-0.2truecm}
$$
\hat f:[T_0,T]\times \calc^P_A((0,+\infty);H) \to \R,
\hbox{such that } f(t,\bar x)=
\hat f(t,y^P_{\bar x}(\cdot \wedge t)),\;\forall
(t,\bar x)\in (T_0,T]\times \overline{H}
\vspace{-0.2truecm}
$$
and\vspace{-0.2truecm}
$$
\hat{\hat f}:[T_0,T]\times \calc^P_A((0,+\infty);H) \to K',
\hbox{such that }
(t-T_0)^\eta\nabla^C f(t,\bar x)=
\hat{\hat f}(t,y^P_{\bar x}(\cdot \wedge t)),\;
\forall (t,\bar x)\in (T_0,T]\times \overline{H}
\vspace{-0.2truecm}
$$
When $\eta =0$ we omit it simply writing
$\cals^{1,P}_{progr}((T_0,T]\times \overline{H})$.
\end{definition}

It turns out that $\cals^{1,P}_{\eta}((0,T]\times\overline{H})$
and $\cals^{1,P}_{\eta,prog} ((0,T]\times\overline{H})$ are closed subspaces of
$C^{0,1,C}_{{\eta}}([0,T]\times \overline H)$
and so they are Banach spaces if they are endowed with the norm $\Vert \cdot\Vert_{C^{0,1,C}_{{\eta}}}$.

\begin{lemma}\label{lemma:SP}
Let Hypotheses \ref{ip-sde-common}, \ref{ip:PC}
and \ref{hp:L2} hold true.
Then the sets $\cals^{1,P}_{\eta}([T_0,T]\times \overline H)$ and $\cals^{1,P}_{\eta,prog}([T_0,T]\times \overline H)$ are
 closed subspaces of $C^{0,1,C}_{\eta}([T_0,T]\times \overline H)$.
\end{lemma}\vspace{-0.2truecm}
\dim
We make the proof for the case $T_0=0$ as the case $T_0>0$ is completely {analogous}.
It is clear that $\cals^{1,P}_{\eta}([0,T]\times \overline H)$
and $\cals^{1,P}_{\eta,prog}([0,T]\times \overline H)$ are
vector subspaces of $C^{0,1,C}_{\eta}([0,T]\times \overline H)$. We prove now that this is true for $\cals^{1,P}_{\eta}([0,T]\times \overline H)$; the assertion follows in the same way for $\cals^{1,P}_{\eta,prog}([0,T]\times \overline H)$.
Take any sequence $f_n \to f$ in $C^{0,1,C}_{\eta}([0,T]\times \overline H)$ such that
\vspace{-0.2truecm}
$$
\hat f_n:[0,T]\times \calc^P_A((0,T];H) \to \R,
\hbox{    } f_n(t,\bar x)=
\hat f_n(t,y^P_{\bar x}),\;\forall
(t,\bar x)\in (0,T]\times \overline{H}
\vspace{-0.2truecm}
$$
and\vspace{-0.2truecm}
$$
\hat{\hat f}_n:[0,T]\times \calc^P_A((0,T];H) \to \R,
\hbox{  }
t^\eta\nabla^C f_n(t,\bar x)=
\hat{\hat f}_n(t,y^P_{\bar x}),\;
\forall (t,\bar x)\in (0,T]\times \overline{H}.
\vspace{-0.2truecm}
$$
The sequence $\{\hat{ f}_n\}$ is a Cauchy sequence with respect to the $L_\infty$-norm for bounded measurable real functions defined on $[0,T]\times \calc^P_A((0,T];H)$. Indeed, for any $\epsilon >0$ and $n,m\in\N$, take
$(t_{\epsilon,n,m},y_{\epsilon,n,m})\in [0,T]\times \calc^P_A((0,T];H)$ such that\vspace{-0.2truecm}
$$
\sup_{(t,y) \in (0,T]\times  \calc^P_A((0,T];H)}|\hat f_n(t,y)-\hat f_m(t,y)|
<
\epsilon + |\hat f_n(t_{\epsilon,n,m},y_{\epsilon,n,m})-\hat f_m(t_{\epsilon,n,m},y_{\epsilon,n,m})|
\vspace{-0.2truecm}
$$
and let $\overline x_{\epsilon,n,m} \in \overline H$ such that
$y_{\epsilon,n,m}=y^P_{\overline x_{\epsilon,n,m}}$
Hence we get\vspace{-0.2truecm}
$$
\sup_{(t,y) \in [0,T]\times \calc^P_A((0,T];H)} |\hat f_n(t,y)-\hat f_m(t,y)|<
\epsilon + |f_n(t_{\epsilon,n,m},\overline x_{\epsilon,n,m})-f_m(t_{\epsilon,n,m},\overline x_{\epsilon,n,m})|
\vspace{-0.2truecm}
$$\vspace{-0.2truecm}
$$\le
\epsilon + \sup_{(t,x) \in [0,T]\times \overline H}
|f_n(t,\overline x)-f_m(t,\overline x)|.
\vspace{-0.2truecm}
$$
Since $\{f_n\}$ is Cauchy, then $\{\hat f_n\}$ is Cauchy, too. So there exists a function
$\hat f :[0,T]\times \calc^P_A((0,T];H)\to \R $ such that\vspace{-0.2truecm}
$$
\sup_{(t,y) \in [0,T]\times \calc^P_A((0,T];H)} |\hat f_n(t,y)-\hat f(t,y)|\to 0 \hbox{ as }n\to\infty
\vspace{-0.2truecm}
$$
This implies that $f(t,\overline x)=\hat f(t,y^P_{\overline x})$ on $[0,T]\times \overline H$.
With the same argument we get that there exists a function
$\hat {\hat{  f}} $ such that\vspace{-0.2truecm}
$$
\sup_{(t,y) \in [0,T]\times \calc^P_A((0,T];H)} |\hat {\hat{ f}}_n(t,y)-\hat {\hat{  f}}(t,y)|\to 0 \hbox{ as }n\to\infty
\vspace{-0.2truecm}
$$
This implies that $(T-t)^\eta \nabla^C f(t,\overline x)=\hat {\hat{  f}}(t,y^P_{\overline x})$ on $(0,T]\times \overline H$.
\qed

\subsubsection{Partial Smoothing on the class $\cals^P_\infty(\overline{H})$}

We need first the following analogous of Hypothesis \ref{ip:NC} in the new ``lifted" setting.

\begin{hypothesis}
\label{hp:smoothingextension}
\begin{itemize}
\item[]
\item [(i)]
For every $t\ge 0$, $k \in K$ we {have}
$\Upsilon^P_\infty\overline{ e^{tA}}Ck \in \Imm
\left(\Upsilon^P_\infty Q_{t} (\Upsilon^P_\infty)^*\right)^{1/2}$.
Consequently, by the closed graph Theorem, the operator
$\widehat\Lambda^{P,C}(t):K\to L^2_\rho(0,\infty;H)$,
$\widehat\Lambda^{P,C}(t)k:=\left(\Upsilon^P_\infty Q_{t} (\Upsilon^P_\infty)^*\right)^{-1/2}
\Upsilon^P_\infty e^{tA}Ck$,
for all $k \in K$,
is well defined for all $t>0$.
\vspace{-0.2truecm}
\item [(ii)]
There exists $\kappa_0>0$ and $\gamma \in (0,1)$ such that
$\|\widehat\Lambda^{P,C}(t)\|_{\call(K,L^2_\rho)} \le \kappa_0 (t^{-\gamma}\vee 1), \qquad
\forall t \ge 0$.
\end{itemize}
\end{hypothesis}

\vspace{-0.3truecm}


\begin{remark}
  \label{rm:NCwhenPcommute}
Note that, when the operator $P^*$ commutes with the semigroup
$e^{sA^*}$, then Hypothesis \ref{hp:smoothingextension}
is equivalent to Hypothesis \ref{ip:NC}.
To see this, observe first that Hypothesis \ref{ip:NC} is equivalent (see e.g.
\cite[Appendix B]{DP1}
and \cite[Part II, Chapter 2]{Z})
to ask that, for given $\gamma \in (0,1)$, $k_1>0$, and for all $t>0$,
\begin{equation}\label{eq:NCdualbasic}
\left|\left(\overline{Pe^{tA}}C\right)^*h\right|^{2}_K
\le \kappa_1 (t^{-\gamma}\vee 1)
\<Q_tP^*hP^*h\>_H,
\qquad \forall h\in H,
\end{equation}
where, using Remark 5.3 of \cite{FGFM-III},
$\left(\overline{Pe^{tA}}C\right)^*h=
C^*e^{tA^*}P^*h$, for all $h\in H$.
Similarly Hypothesis \ref{hp:smoothingextension} is equivalent to ask that, for given $\gamma \in (0,1)$, $k_1>0$, and for all $t>0$,\vspace{-0.2truecm}
\begin{equation}\label{eq:NCdualhigh}
|\left(\Upsilon^P_\infty \overline{e^{tA}}C\right)^*z|^{2}_K
\le c_t
\<Q_t(\Upsilon^P_\infty)^*z,(\Upsilon^P_\infty)^* z\>_H,
\qquad \forall z\in L^2_\rho,
\end{equation}
where, by Proposition \ref{pr:propcalS}-(iv),
$
\left(\Upsilon^P_\infty \overline{ e^{tA}}C\right)^*z=
C^*e^{tA^*}\int_{0}^{+\infty}e^{-\rho s}
e^{sA^*}P^*z(s)ds$, for all $z\in L^2_\rho$.
Now, if $P^*$ commutes with the semigroup
$e^{sA^*}$, we have, for every $z\in L^2_\rho$,
$\left(\Upsilon^P_\infty\overline{e^{tA}}C\right)^*z=
C^*e^{tA^*}P^*
\int_{0}^{+\infty}e^{-\rho s}e^{sA^*}z(s)ds$,
and\vspace{-0.2truecm}
$$
\<Q_t(\Upsilon^P_\infty)^*z,(\Upsilon^P_\infty)^* z\>_H,
=
\<Q_tP^*\int_{0}^{+\infty}e^{-\rho s}e^{sA^*}z(s)ds ,P^*\int_{0}^{+\infty}e^{-\rho s}e^{sA^*}z(s)ds \>_H,
$$
Hence, if \eqref{eq:NCdualbasic} holds for all $h \in H$, a fortiori, it must be true for $h$ of the form
$\int_{0}^{+\infty}e^{-\rho s}e^{sA^*}z(s)ds$ for some $z\in L^2_\rho$. This means that \eqref{eq:NCdualhigh} holds.
\\
On the other hand, assume that \eqref{eq:NCdualhigh} holds.
This means that \eqref{eq:NCdualbasic} holds for all $h\in H$ of the form $\int_{0}^{+\infty}e^{-\rho s}e^{sA^*}z(s)ds$ for some $z\in L^2_\rho$. Since, when $z$ is constant, we have
$\int_{0}^{+\infty}e^{-\rho s}e^{sA^*}zds
=(\rho-A^*)^{-1}z$
then the set of such $h$ contains $D(A^*)$ and, consequently, is dense in $H$. This immediately implies that  \eqref{eq:NCdualbasic} holds for all $h\in H$.
\end{remark}


We now state and prove an extension of Proposition \ref{prop:partsmooth} concerning the smoothing properties of the semigroup $R_t$ on functions of paths.
First of all we notice that (see \eqref{ornstein-sem-gen})
when $\phi \in\cals^P_\infty(\overline H)$,
$R_t [\phi]$ is defined on the whole $\overline H$.
Indeed, if $\phi \in \cals^P_\infty(\overline H)$,
then there exists $\hat \phi$ such that, for all $\overline x \in \overline H$,
$\phi(\bar x)=\hat \phi(y^P_{\bar x})$. Hence\vspace{-0.2truecm}
$$
R_t[ \phi](\overline x):=
\int_{H} \phi\left(z+\overline{e^{tA}}\bar x\right)\caln(0,Q_t)(dz)
=
\int_{H} \hat\phi\left(\Upsilon^P_\infty(z+\overline{e^{tA}}\bar x )\right)\caln(0,Q_t)(dz),\vspace{-0.2truecm}
$$
where the last equality gives sense to $R_t[ \phi](\overline x)$.

\begin{proposition}\label{prop:partsmoothnew}
Let Hypotheses \ref{ip-sde-common}, \ref{ip:PC}, \ref{hp:L2}, and
\ref{hp:smoothingextension}-(i)
hold true.
Then $R_t,\,t>0$ maps functions $\phi\in \cals^P_\infty(\overline{H})$
into functions which are $C$-Fr\'echet differentiable in $\overline{H}$, and the $C$-derivative is given,
for all $t>0$, $\bar x\in \overline{H}$, by\vspace{-0.2truecm}
\begin{align}\label{eq:formulader-gen-Pnew}
&\nabla^C(R_{t}[\phi])(\bar x)k =
\int_{L^2_\rho}\hat\phi\left(z_1+ \Upsilon^P_\infty \bar x \right)
\<\widehat\Lambda^{P,C}(t) k, (\Upsilon^P_\infty Q_t(\Upsilon^P_\infty)^*)^{-1/2} z_1\>_{L^2_\rho}
\caln(0,(\Upsilon^P_\infty Q_t(\Upsilon^P_\infty)^*))(dz_1)
\notag
\\[2mm]
&=
\E\left[\hat\phi\left(\Upsilon^P_\infty
X(t;\bar x)\right)
\<\widehat\Lambda^{P,C}(t) k,
(\Upsilon^P_\infty Q_t(\Upsilon^P_\infty)^*)^{-1/2}\Upsilon^P_\infty W_A(t)\>_{L^2_\rho}
\right]
\end{align}
Moreover, for any $\phi\in \cals_\infty(\overline{H})$, $t>0$, $\bar x\in \overline{H}$, $k\in K$,\vspace{-0.2truecm}
\begin{equation}\label{norm-Cdernew}
\vert \<\nabla^C R_t[\phi](\bar x), k\>\vert \leq \|\widehat\Lambda^{P,C}(t)\|_{\call(K,L^2_\rho(0,\infty;H)}
\Vert\phi\Vert_\infty \cdot\vert k\vert.
\end{equation}
Furthermore, if $\hat\phi$ associated to $\phi$ in \eqref{eq:defSbar} is continuous, then
the map $(t,\bar x) \to\nabla^C R_t[\phi](\bar x)$
belongs to $C((0,T]\times \overline{H};K)$.
\\
Assume now that, beyond the continuity of $\hat\phi$, Hypothesis \ref{hp:smoothingextension}-(ii) holds.
Then the map $(t,\bar x)\to R_t[\phi](\bar x)$ belongs to
the space
$\cals^{1,P}_{\gamma}((0,T]\times \overline{H})$, hence to
$C_\gamma^{0,1,C}([0,T]\times \overline{H})$.
\\
Finally, if Hypothesis \ref{ip:NC} is verified and $\phi \in C_b^P(H)$ then the map $(t,\bar x)\to R_t[\phi](\bar x)$ belongs to
$\cals^{1,P}_{\gamma,prog}([0,T]\times \overline H)$.
\end{proposition}

\dim
Let $\phi\in \cals^P_\infty(\overline{H})$, $t>0$, $\bar x \in \overline{H}$, $k\in K$, $\alpha \in \R$.
By \eqref{ornstein-sem-gen} and \eqref{eq:defSbar} we have\vspace{-0.2truecm}
$$
R_{t}[\phi](\bar x+\alpha Ck)
=
\int_{H} \phi\left(z+\overline{e^{tA}}(\bar x+\alpha Ck)\right)\caln(0,Q_t)(dz)
=
\int_{H} \hat\phi\left(\Upsilon^P_\infty(z+\overline{e^{tA}}(\bar x+\alpha Ck))\right)\caln(0,Q_t)(dz)
$$
Now we perform the change of variable
$z_1= \Upsilon^P_\infty (z+  \overline{e^{tA}}\alpha Ck)$,
for $z \in H$, $z_1 \in \calc_A^P((0,+\infty);{H})\subset L^2_\rho$.
which gives\vspace{-0.2truecm}
$$
R_{t}[\phi](\bar x+\alpha Ck)
=
\int_{L^2_\rho}
\hat\phi\left(z_1+\Upsilon^P_\infty \overline{e^{tA}} \bar x\right)
\caln(\Upsilon^P_\infty \overline{e^{tA}}\alpha Ck,\Upsilon^P_\infty Q_t(\Upsilon^P_\infty)^*)(dz_1).\vspace{-0.2truecm}
$$
We then compute the incremental ratio in the direction $Ck$ getting,
for $\phi\in \cals^P_\infty(\overline{H})$, $t>0$, $\bar x\in \overline{H}$, $k\in K$, $\alpha \in \R-\{0\}$,\vspace{-0.2truecm}
\begin{align}
\label{eq:rappincrnew}
\frac{1}{\alpha}
&\left[R_{t}[\phi](\bar x+\alpha Ck)-R_{t}[\phi](\bar x)\right]=
\\[2mm]
=&\frac{1}{\alpha}
\left[
\int_{L^2_\rho}
\hat\phi(z_1+\Upsilon^P_\infty \overline{e^{tA}}\bar x))
\caln\left(\Upsilon^P_\infty \overline{e^{tA}}\alpha Ck,\Upsilon^P_\infty Q_t(\Upsilon^P_\infty)^*\right)(dz_1)\right.
\left. -
\int_{L^2_\rho}
\hat\phi(z_1+\Upsilon^P_\infty \overline{e^{tA}}\bar x))
\caln(0,\Upsilon^P_\infty Q_t(\Upsilon^P_\infty)^*)(dz_1)
\right].
\notag
\end{align}
We now apply Cameron-Martin Theorem, see e.g. \cite[Theorem 1.3.6]{DP3}. This gives that the Gaussian measures
\newline$\caln\left(\Upsilon^P_\infty \overline{e^{tA}}\alpha Ck,\Upsilon^P_\infty Q_t(\Upsilon^P_\infty)^*\right)$
and
$\mathcal{N}\left(0,\Upsilon^P_\infty Q_t(\Upsilon^P_\infty)^*\right)$
are equivalent if and only if
$\Upsilon^P_\infty \overline{e^{tA}}\alpha Ck\in\operatorname{Im}\left(\Upsilon^P_\infty Q_t(\Upsilon^P_\infty)^*\right)^{1/2}$, which is guaranteed by Hypothesis \ref{hp:smoothingextension}-(i).
In such case, and for any $y \in \operatorname{Im}\left(\Upsilon^P_\infty Q_t(\Upsilon^P_\infty)^*\right)^{1/2}$, the density is given by\vspace{-0.2truecm}
\begin{align}
&d(t,y,z)     =\frac{d\caln\left(y,\Upsilon^P_\infty Q_t(\Upsilon^P_\infty)^*\right)}
{d\mathcal{N}\left(0,\Upsilon^P_\infty Q_t(\Upsilon^P_\infty)^*\right)  }(z)
\label{eq:density1}
 \\
\notag
&=\exp\left\{\left\langle (\Upsilon^P_\infty Q_t(\Upsilon^P_\infty)^*)^{-1/2}
y,(\Upsilon^P_\infty Q_t(\Upsilon^P_\infty)^*)^{-1/2}z\right\rangle_{L^2_\rho}
-\frac{1}{2}\left|
(\Upsilon^P_\infty Q_t(\Upsilon^P_\infty)^*)^{-1/2}y
\right|_{L^2_\rho}^{2}\right\}  .
\end{align}
This map is defined for all $z\in{ \overline{\operatorname{Im} \left(\Upsilon^P_\infty Q_t(\Upsilon^P_\infty)^*\right)^{-1/2}}}
=\left[\ker\left(\Upsilon^P_\infty Q_t(\Upsilon^P_\infty)^*\right)\right]^\perp\subseteq {L^2_\rho}$, see, e.g., \cite[Proposition 1.59]
{FabbriGozziSwiech}),{ and, see also \cite{DP3}[Section1.2.4], the function $z\mapsto \left\langle (\Upsilon^P_\infty Q_t(\Upsilon^P_\infty)^*)^{-1/2}y,(\Upsilon^P_\infty Q_t(\Upsilon^P_\infty)^*)^{-1/2}z\right\rangle_{L^2_\rho}$ is well defined and square integrable with respect to the measure $\mathcal{N}\left(0,\Upsilon^P_\infty Q_t(\Upsilon^P_\infty)^*\right)$.}
Hence, from \eqref{eq:rappincrnew} we get\vspace{-0.2truecm}
\begin{align}
&
\frac{1}{\alpha}
\left[R_{t}[\phi](\bar x+\alpha Ck)-R_{t}[\phi](\bar x)\right]=
\int_{L^2_\rho}\hat\phi(z_1+\Upsilon^P_\infty \overline{e^{tA}}\bar x)
\frac{d(t,\Upsilon^P_\infty \overline{e^{tA}}\alpha Ck,z_1)-1}{\alpha}
\caln(0,\Upsilon^P_\infty Q_t(\Upsilon^P_\infty)^*)(dz_1).
\label{eq:rappincralfa}
\end{align}
Now we observe that, by the definition of $\widehat\Lambda^{P,C}(t)$,\vspace{-0.2truecm}
\begin{equation}\label{eq:d-1}
\frac{d(t,\alpha\Upsilon^P_\infty \overline{e^{tA}}Ck,z_1)-1}{\alpha}
=\frac{1}{\alpha}\left[\exp\left\{\alpha
\left\langle
\widehat\Lambda^{P,C}(t)k,
(\Upsilon^P_\infty Q_t(\Upsilon^P_\infty)^*)^{-1/2}z_1
\right\rangle_{L^2_\rho}
-\frac{\alpha^2}{2}
\left|\widehat\Lambda^{P,C}(t)k\right|_H^{2}\right\}
-1\right].
\end{equation}
When $\alpha \to 0$,
for every $z_1\in \left[\ker\left(\Upsilon^P_\infty Q_t(\Upsilon^P_\infty)^*\right)\right]^\perp$, the above term converges
to $\left\langle\widehat\Lambda^{P,C}(t)k,
(\Upsilon^P_\infty Q_t(\Upsilon^P_\infty)^*)^{-1/2}
z_1\right\rangle_{L^2_\rho}$.
Now consider the probability space $L^2_\rho$ with the Borel $\sigma$-field
and the Gaussian measure
$\caln(0,\Upsilon^P_\infty
Q_t(\Upsilon^P_\infty)^*)(dz_1)$.
In this space the term $\eqref{eq:d-1}$ converges a.s. for $\alpha \to 0$
and (see, e.g., \cite[Proposition 1.59]{FabbriGozziSwiech})
the map\vspace{-0.2truecm}
$$
z_1 \mapsto \widehat{\mathcal{Q}_t} (z_1):= \left\langle\widehat\Lambda^{P,C}(t)k,
(\Upsilon^P_\infty Q_t(\Upsilon^P_\infty)^*)^{-1/2}
z_1\right\rangle_{L^2_\rho}
$$
is a real valued Gaussian random variable with mean $0$ and variance $\left|\widehat\Lambda^{P,C}(t)k\right|_{L^2_\rho}^{2}$
(see, e.g. \cite[Remark 2.2]{DPZ91}).
So in particular, for all $L>0$, $\E[e^{L|\widehat{\mathcal{Q}_t}| }]<+\infty$.
Now, calling $k(\alpha)$ the exponent in \eqref{eq:d-1},
we set, for $\alpha \in (-1,1)$\vspace{-0.2truecm}
$$
 g(\alpha):=\frac{d(t,\alpha\Upsilon^P_\infty \overline{e^{tA}}Ck,z_1)-1}{\alpha}
=\frac{e^{ k(\alpha)}-1}{\alpha}=
\frac{e^{k(\alpha)}-1}{ k(\alpha)}\cdot \frac{k(\alpha)}{\alpha}
\quad \Longrightarrow \quad
|g(\alpha)|\le \left|
\frac{e^{k(\alpha)}-1}{ k(\alpha)}
\right| \cdot
\left|
\frac{k(\alpha)}{\alpha}\right|.
$$
Now, for $\alpha\in (-1,1)$,
$k(\alpha)\vee (k(\alpha)/\alpha)\le |\widehat{\mathcal{Q}_t}|+\frac12 \left|\widehat
\Lambda^{P,C}(t)k\right|_{L^2_\rho}^{2}$.
Since the function $y\to \frac{e^y-1}{y}$ is increasing over $\R$ then
\vskip-0.2truecm
$$
\left|
\frac{e^{k(\alpha)}-1}{ k(\alpha)}\right|
\le
\left|\frac{e^{|\widehat{\mathcal{Q}_t}|+\frac12 \left|\widehat
\Lambda^{P,C}(t)k\right|_{L^2_\rho}^{2}}-1}
{|\widehat{\mathcal{Q}_t}|+\frac12 \left|\widehat
\Lambda^{P,C}(t)k\right|_{L^2_\rho}^{2}}
\right|
\qquad and \qquad
\left|g(\alpha)\right|
\le
e^{|\widehat{\mathcal{Q}_t}|+\frac12 \left|\widehat
\Lambda^{P,C}(t)k\right|_{L^2_\rho}^{2}}
$$
Hence we can apply the dominated convergence theorem to \eqref{eq:rappincralfa} getting the limit\vspace{-0.2truecm}
\begin{align*}
&
\lim_{\alpha\rightarrow 0}\frac{1}{\alpha}
\left[R_{t}[\phi](\bar x+\alpha Ck)-R_{t}[\phi](\bar x)\right]=
\int_{L^2_\rho}\hat\phi\left(z_1+\Upsilon^P_\infty
\overline{e^{tA}}\bar x\right)
\left\langle\widehat\Lambda^{P,C}(t)k,
(\Upsilon^P_\infty Q_t(\Upsilon^P_\infty)^*)^{-1/2}
z_1\right\rangle_{L^2_\rho}
\caln(0,\Upsilon^P_\infty Q_t(\Upsilon^P_\infty)^*)(dz_1)
\end{align*}
Consequently, along Definition \ref{df4:Gderunbounded}-(i),
there exists the $C$-directional derivative
$\nabla^C R_{t}\left[\phi\right](\bar x;k)$
which is equal to the above right hand side.
Using that $\widehat\Lambda^{P,C}(t)$ is continuous we
see that the above limit is uniform for $k$ in the unit ball of $K$, so (Definition \ref{df4:Gderunbounded}-(iii)) there exists the $C$-Fr\'echet derivative $\nabla^C R_{t}\left[\phi\right](\bar x)$.
From the above and from
\cite[Proposition 1.59]{FabbriGozziSwiech} we get\vspace{-0.2truecm}
\vskip-0.6truecm
\begin{align*}
|\nabla^C R_{t}\left[\phi\right](\bar x;k)|
&\leq
\Vert \hat\phi\Vert_\infty
\left(\int_{L^2_\rho}
\<\widehat\Lambda^{P,C}(t) k,
(\Upsilon^P_\infty Q_t(\Upsilon^P_\infty)^*)^{-1/2}
z_1\>_{L^2_\rho}^2
\caln(0,\Upsilon^P_\infty Q_t(\Upsilon^P_\infty)^*)(dz)\right)^{1/2}
\\[3mm]
\nonumber
&
= \Vert \hat\phi\Vert_\infty
\Vert \widehat\Lambda^{P,C}(t)k  \Vert_{{L^2_\rho}}
  \le \Vert\phi\Vert_\infty
\Vert \widehat\Lambda^{P,C}(t)\Vert_{\call(K;L^2_\rho)} |k|_K . \nonumber
\end{align*}
This gives the required estimate \eqref{norm-Cdernew}.
The statement on continuity follows
using exactly the same argument as in
\cite[Theorem 4.41-(ii)]{FabbriGozziSwiech}.
The latter above, together with Hypothesis \ref{hp:smoothingextension}-(ii), implies, recalling the last part of Definition \ref{df4:Gspaces},
that the map $f$ defined by
$(t,\bar x)\mapsto R_t[\phi](\bar x)$ belongs to
$C_\gamma^{0,1,C}([0,T]\times \overline{H})$.
Now, we prove that such map belongs to
$\cals^{1,P}_{\gamma}((0,T]\times \overline{H})$, i.e. we provide the two functions $\hat{f}$ and $\hat{\hat{f}}$ as required by Definition \ref{df:spaziSeta}.
Indeed observe that, for $t>0$ and $\bar x\in \overline{H}$ we have\vspace{-0.2truecm}
$$
R_{t}[\phi](\bar x)
=
\int_{H} \phi
(z+\overline{e^{tA}}(\bar x))\caln(0,Q_t)(dz)
=
\int_{H} \hat\phi\left(\Upsilon^P_\infty
(z+\overline{e^{tA}}\bar x)\right)\caln(0,Q_t)(dz)
$$\vspace{-0.2truecm}
$$
=
\int_{L^2_\rho}
\hat\phi\left(z_1+\Upsilon^P_\infty \overline{e^{tA}}\bar x\right)
\caln(0,\Upsilon^P_\infty Q_t(\Upsilon^P_\infty)^*)(dz_1).
$$
Hence, we set\vspace{-0.2truecm}
$$
\hat{f}(t,\Upsilon^P_\infty \bar x)
:=
\int_{L^2_\rho}
\hat\phi\left(z_1+\Upsilon^P_\infty \overline{e^{tA}}\bar x\right)
\caln(0,\Upsilon^P_\infty Q_t(\Upsilon^P_\infty)^*)(dz_1).
$$
This satisfies \eqref{eq:defS} thanks to \eqref{eq:Petasemigroup}.
A similar argument holds for
$\nabla^C(R_{t}[\phi])(\bar x)$ using \eqref{eq:formulader-gen-Pnew}.
\\
Assume now that $\phi \in C_b^P(H)$. Then, by Definition \ref{df:spaziphi1} and Hypothesis \ref{ip:PC}-(iii),
we get\vspace{-0.2truecm}
$$
R_{t}[\phi](\bar x)
=
\int_{H} \phi\left(P
(z+\overline{e^{tA}}\bar x)\right)\caln(0,Q_t)(dz)
=
\int_{\overline H}
\bar\phi\left(z_1+\overline{Pe^{tA}}\bar x\right)
\caln(0, PQ_tP^*)(dz_1).
$$
Hence, we set\vspace{-0.2truecm}
$$
\hat{f}(t,y^P_{\bar x}(\cdot\wedge t))
:=
\int_{\overline H}
\bar\phi\left(z_1+\overline{Pe^{tA}}\bar x\right)
\caln(0,P Q_tP^*)(dz_1).
$$
A similar argument holds for
$\nabla^C(R_{t}[\phi])(\bar x)$ using that
Hypothesis \ref{ip:NC}, and then
\eqref{eq:formulader-gen-Pold} are verified.
This implies that such map belongs to
$\cals^{1,P}_{\gamma,prog}([0,T]\times \overline H)$.
\hfill\qed

\vspace{-0.2truecm}
\begin{remark}\label{rm:partsmooth-secondnew}
Generalizing to our setting the ideas of Proposition 4.5 in \cite{FGFM-I}
it is possible to prove that, if $\phi$ is more regular
(i.e. $\phi\in C^1_b(H)\cap C_b^P(H)$,
also $\nabla^C R_t[\phi]$ has more regularity, i.e.
$\nabla\nabla^{C}R_{t}\left[\phi\right]$,
$\nabla^{C}\nabla R_{t}\left[\phi\right]$ exist, coincide, and
satisfy suitable formulae and estimates.
Such type of results could be useful to find optimal feedback controls.
\end{remark}
\vspace{-0.2truecm}

\subsection{New results on convolutions}
\label{sec:convpartsmooth-abstr-setting}

Here we prove that, given any $T\in (0,+\infty]$,
for any element $f$ of a suitable family of functions
in $C_b([0,T]\times H)$,
a partial smoothing property (similar to the one of the previous section) for the convolution integral
$\int_{0}^{t}R_{t-s} [f(s,\cdot)](x) ds$ holds.
We first formulate the assumption we need: it is the analogue of Hypothesis \ref{hp:smoothingextension} in a form which holds in our examples and which allows to prove the partial smoothing property.

\begin{hypothesis}\label{hp:smoothingextension-conv}
Let $T_1>0$ be given.
\begin{itemize}
\item [(i)]
For every $T_1\ge t>s> 0$ and $k \in K$ we have\vspace{-0.2truecm}
\begin{equation}\label{hp:smoothingextension-i}
\Upsilon^P_s \overline{e^{(t-s)A}}Ck \in \Imm
\left(\Upsilon^P_s Q_{t-s} (\Upsilon^P_s)^*\right)^{1/2}
\end{equation}
Consequently, by the closed graph Theorem, the operator\vspace{-0.2truecm}
$$
\widehat\Lambda^{P,C}_s(t-s):K\to L^2(0,s;H),
\qquad
\widehat\Lambda^{P,C}_s(t-s)k:=
\left(\Upsilon^P_sQ_{t-s} (\Upsilon^P_s)^*\right)^{-1/2}
\Upsilon^P_s \overline{e^{(t-s)A}}Ck
\quad \forall k \in K,
$$
is well defined for all $0<s<t\le T_1$.
\item [(ii)]
There exists $\kappa_{T_1} > 0$ and
$\gamma \in (0, 1)$ such that\vspace{-0.2truecm}
\begin{equation}\label{hp:smoothingextension-ii}
\|\widehat\Lambda^{P,C}_s(t-s) \|_{L^2(0,s;H)}
\le \kappa_{T_1}  (t-s)^{-\gamma}|k|_K, \qquad \forall T_1\ge t>s>0, \quad \forall k \in K
\end{equation}
\end{itemize}
\end{hypothesis}

\begin{remark}\label{lemma:ip:contr-sollevata-equiv-1}
\begin{itemize}
  \item [(i)]
Note that, the above Hypothesis \ref{hp:smoothingextension-conv}, needed to prove the regularity of the convolution integral in next Proposition \ref{prop:partsmoothnew-conv}, is not equivalent\footnote{It is indeed weaker in the sense that it is implied by Hypothesis  \ref{hp:smoothingextension}.} to the previous Hypothesis \ref{hp:smoothingextension}
needed to prove the regularity of $R_t[\phi]$.
This may seem a bit strange as
in the paper \cite{FGFM-III} (where the ``lifting"
procedure was not used) there is no need to add a different hypothesis to establish the regularity of the convolution integral
(see \cite[Hypothesis 5.7]{FGFM-III}).
The reason comes from the fact that, in the case of delay in the control (which is one of our motivating examples, see Subsection 5.2 for more details)
Hypothesis \ref{hp:smoothingextension} does not hold while
Hypothesis \ref{hp:smoothingextension-conv} does (only for a small $T_1$, but this is enough for our purposes).\vspace{-0.2truecm}
 \item [(ii)]
Arguing as in Remark \ref{rm:NCwhenPcommute}, it is possible to check that Hypothesis \ref{hp:smoothingextension-conv} is equivalent to ask that for any $0<s <t\le T_1$ and for any $z\in L^2(0,s;H)$
\begin{equation}\label{condiz:ip:contr-sollevata-equiv}
\Vert C^* e^{(t-s)A^*}(\Upsilon^P_s)^* z\Vert_K\leq \kappa_{T_1} (t-s)^{-\gamma}.
\Vert\left(\Upsilon^P_s Q_{t-s}(\Upsilon^P_s)^*\right)^{1/2}z \Vert
\end{equation}
with $\kappa_{T_1}$ and $\gamma$ as in \eqref{hp:smoothingextension-ii}.\vspace{-0.2truecm}
\item[(iii)]
As in Remark \ref{rm:NCwhenPcommute} we have that, when $P^*$ and $e^{tA^*}$ commute, Hypothesis \ref{hp:smoothingextension-conv} and Hypothesis \ref{ip:NC} are equivalent.
 \end{itemize}
\end{remark}


\begin{proposition}\label{prop:partsmoothnew-conv}
Let Hypotheses \ref{ip-sde-common}, \ref{ip:PC}, \ref{hp:L2}
and
\ref{hp:smoothingextension-conv} hold.
Let $T\in (0,T_1]$ (where $T_1$ comes from Hypothesis \ref{hp:smoothingextension-conv}) and $\eta \in [0,1)$.
Assume that $f\in \cals^P_{\eta,prog}((0,T]\times \overline H)$.
Then, for every $t\in [0,T]$, $\bar x \in \overline{H}$, the convolution\vspace{-0.2truecm}
\begin{equation}\label{eq:convgdef}
g(t,\bar x):=\int_{0}^{t}R_{t-s} [f(s,\cdot)](\bar x) ds
\end{equation}
is well defined and $C$-Fr\'echet differentiable: the $C$-derivative belongs to $\cals^P_{\gamma,prog}((0,T]\times \overline H;K)$ and it is given by\vspace{-0.2truecm}
\begin{align}\label{eq:formulader-gen-Pnewconv}
&\nabla^C [g(t,\cdot)](\bar x)\\ \nonumber
&\;=
\int_{0}^{t}\int_{L^2_s}
\hat f\left(s,z_1+\Upsilon^P_s e^{(t-s)A}\bar x\right)
\<\widehat\Lambda_s^{P,C}(t-s) k, (\Upsilon^P_s Q_{t-s}(\Upsilon^P_s)^*)^{-1/2} z_1\>_{L^2_s}
\caln(0,\Upsilon^P_s Q_{t-s}(\Upsilon^P_s)^*)(dz_1)ds
\end{align}
Moreover, for all $f\in \cals^P_{\eta,prog}((0,T]\times \overline H)$ and any $t\in (0,T]$, $\bar x \in \overline{H}$, $k\in K$ (here $\beta(\cdot,\cdot)$ is the Euler beta function),\vspace{-0.2truecm}
\begin{equation}\label{norm-Cderconv}
\vert \<\nabla^C [g(t,\cdot)](\bar x), k\>\vert
\leq
\Vert t^\eta f\Vert_\infty \vert k\vert
\int_{0}^{t} {\kappa}_T s^{-\eta} (t-s)^{-\gamma} ds
\le
\Vert t^\eta f\Vert_\infty \vert k\vert\cdot
{\kappa}_T \beta(\eta,\gamma) t^{-\gamma}
\end{equation}
Furthermore, if the $\hat f$ associated to $f$ in \eqref{eq:defSbar} is continuous, then
the map $g$ and $\nabla^C g$ are continuous and $g$ belongs to the space $\cals^{1,P}_{\gamma,prog}((0,T]\times \overline{H})$, hence to
$C_\gamma^{0,1,C}([0,T]\times \overline{H})$.
This holds also when Hypothesis \ref{ip:NC} is verified and $f$ is such that the map
$(t,x)\mapsto t^\eta f(t,x)$ belongs to $C_b^P([0,T]\times H)$ for some $\eta\in [0,1)$.
\end{proposition}

\dim
First, observe that $g$ is well defined for $(t,\bar x)\in [0,T] \times\overline{H}$ and that $t=0\Rightarrow g\equiv 0$. Let now $t\in (0,T]$. Then\vspace{-0.2truecm}
\begin{align*}
g(t,\bar x)&=\int_{0}^{t}
R_{t-s} \left[f(s,\cdot)\right](\bar x)  ds
=\int_{0}^{t}
\int_H f\left(s,z+\overline{e^{(t-s)A}}\bar x
\right)\caln(0,Q_{t-s})(dz) ds
\\
&=
\int_{0}^{t}
\int_H \hat f\left(s,
\Upsilon^P_s (z+\overline{e^{(t-s)A}} \bar x)\right)
\caln(0,Q_{t-s})(dz) ds
\end{align*}
and the last is well defined.
Now we perform, in the interior integral, the change of variable
$z_1= \Upsilon^P_s z$ ($z \in H$, $z_1 \in L^2(0,s;H)=:L^2_s$).
which gives\vspace{-0.2truecm}
$$
\int_H \hat
f\left(s, \Upsilon^P_s (z+\overline{e^{(t-s)A}}\bar x)\right)\caln(0,Q_{t-s})(dz)
=
\int_{L^2_s} \hat f\left(s, z_1+\Upsilon^P_s
\overline{e^{(t-s)A}}\bar x \right)\caln\left(0,\Upsilon^P_s Q_{t-s} (\Upsilon^P_s)^*\right)(dz_1)
$$
Similarly, for every $t\in (0,T]$, $\bar x \in \overline{H}$,
$k\in K$, $\alpha \in \R$, we get\vspace{-0.2truecm}
\begin{align*}
g(t,\bar x+\alpha Ck)
&
=\int_{0}^{t}
\int_H f(s,z+\overline{e^{(t-s)A}}(\bar x+\alpha Ck)\caln(0,Q_{t-s})(dz) ds
\\
&=
\int_{0}^{t}
\int_H \hat f\left(s, \Upsilon^P_s (z+\overline{e^{(t-s)A}}\bar x+ \overline{e^{(t-s)A}}\alpha Ck)\right)\caln(0,Q_{t-s})(dz)
\\
&=
\int_{0}^{t}
\int_{L^2_s} \hat f\left(s,z_1+\Upsilon^P_s \overline{e^{(t-s)A}}\bar x+\Upsilon^P_s\overline{e^{(t-s)A}}\alpha Ck \right)
\caln\left(0,\Upsilon^P_s Q_{t-s} (\Upsilon^P_s)^*\right)(dz_1) ds
\\
&=
\int_{0}^{t}
\int_{L^2_s} \hat f\left(s, z_1+\Upsilon^P_s \overline{e^{(t-s)A}} \bar x \right)
\caln\left(\Upsilon^P_s\overline{ e^{(t-s)A}}\alpha Ck,\Upsilon^P_s Q_{t-s} (\Upsilon^P_s)^*\right)(dz_1) ds
\end{align*}
The Gaussian measures $\caln\left(0,\Upsilon^P_s Q_{t-s} (\Upsilon^P_s)^*\right)(dz_1)$
and
$\caln\left(\Upsilon^P_s\overline{ e^{(t-s)A}}\alpha Ck,\Upsilon^P_s Q_{t-s} (\Upsilon^P_s)^*\right)(dz_1)$
are equivalent if and only if\footnote{This follows from the Cameron-Martin Theorem, see e.g. \cite[Theorem 1.3.6]{DP3}.}
$\Upsilon^P_s \overline{ e^{(t-s)A}}\alpha Ck\in\operatorname{Im}\left(
\Upsilon^P_s Q_{t-s} (\Upsilon^P_s)^*\right)^{1/2}$, which is guaranteed by Hypothesis \ref{hp:smoothingextension-conv}-(i).
Then, similarly to \eqref{eq:density1}, we set, for $y \in \operatorname{Im}\left(\Upsilon^P_s Q_{t-s} (\Upsilon^P_s)^*\right)^{1/2}$,\vspace{-0.2truecm}
\begin{align}
&d(s,t-s,y,z)     =\frac{d\caln\left(y,\Upsilon^P_s Q_{t-s} (\Upsilon^P_s)^*\right)}
{d\mathcal{N}\left(0,\Upsilon^P_s Q_{t-s} (\Upsilon^P_s)^*\right)  }(z)
\notag \\
 &=\exp\left\{\left\langle (\Upsilon^P_s Q_{t-s} (\Upsilon^P_s)^*)^{-1/2}
y,(\Upsilon^P_s Q_{t-s} (\Upsilon^P_s)^*)^{-1/2}z\right\rangle_{L^2_s}
-\frac{1}{2}\left|(\Upsilon^P_s Q_{t-s} (\Upsilon^P_s)^*)^{-1/2}y\right|_{L^2_s}^{2}\right\}  ,\label{eq:density1-conv}
\end{align}
which, as a function of $z$, is defined for all $z\in {L^2_s}$, see, e.g., \cite[Proposition 1.59]{FabbriGozziSwiech}).
Hence, arguing as in the proof of Proposition  \ref{prop:partsmoothnew} we get\vspace{-0.2truecm}
\begin{align}
&\frac{1}{\alpha}
\left[g(t,\bar x+\alpha Ck)-g(t,\bar x)\right]=
\frac{1}{\alpha}
\left[\int_{0}^{t}
\int_{L^2_s} \hat f\left(s, z_1+\Upsilon^P_s \overline{ e^{(t-s)A}}\bar x \right)
\caln\left(\Upsilon^P_s
\overline{ e^{(t-s)A}}\alpha Ck,\Upsilon^P_s Q_{t-s} (\Upsilon^P_s)^*\right)(dz_1) ds\right.
\notag
\\[2mm]
\notag
&
\left.\qquad \qquad -\int_{0}^{t}\int_{L^2_s}
\hat f\left(s,z_1+\Upsilon^P_s \overline{ e^{(t-s)A}}\bar x \right)
\caln\left(0,\Upsilon^P_s Q_{t-s} (\Upsilon^P_s)^*\right)(dz_1) ds
\right]\\[2mm]
&=
\frac{1}{\alpha}\int_{0}^{t}
\int_{L^2_s} \hat f\left(s, z_1+\Upsilon^P_s \overline{ e^{(t-s)A}}\bar x \right)
\frac{d(s,t-s,\Upsilon^P_s \overline{ e^{(t-s)A}}\alpha Ck,z_1)-1}{\alpha}
\caln\left(0,\Upsilon^P_s Q_{t-s} (\Upsilon^P_s)^*\right)(dz_1) ds
\label{eq:convalphanew}
\end{align}
At this point we can argue exactly as in the proof of Proposition \ref{prop:partsmoothnew} after \eqref{eq:rappincralfa} to send $\alpha \to 0$ in the internal integral of \eqref{eq:convalphanew} applying the dominated convergence theorem. We get, for every $0< s \le t$,\vspace{-0.2truecm}
\begin{align}\label{eq:limitalfaconv}
&
\lim_{\alpha \to 0^+}\frac{1}{\alpha}
\int_{L^2_s} \hat f\left(s, z_1+\Upsilon^P_s \overline{ e^{(t-s)A}}\bar x \right)
\frac{d(s,t-s,\Upsilon^P_s \overline{ e^{(t-s)A}}\alpha Ck,z_1)-1}{\alpha}
\caln\left(0,\Upsilon^P_s Q_{t-s} (\Upsilon^P_s)^*\right)(dz_1)
 \\[2mm]
&
=\int_{L^2_s}
\hat f\left(s,z_1+\Upsilon^P_s \overline{ e^{(t-s)A}}\bar x\right)
\<\widehat\Lambda_s^{P,C}(t-s) k, (\Upsilon^P_s Q_{t-s}(\Upsilon^P_s)^*)^{-1/2} z_1\>_{L^2_s}
\caln(0,\Upsilon^P_s Q_{t-s}(\Upsilon^P_s)^*)(dz_1)
\notag
\end{align}
Now we use Hypothesis \ref{hp:smoothingextension-conv}-(ii)
to guarantee that the limit in \eqref{eq:limitalfaconv}
is integrable in time, getting\vspace{-0.2truecm}
\begin{align*}
&\<\nabla^Cg(t,\bar x) , k\>_K
=
\int_{0}^{t}\int_{L^2_s}
\hat f\left(s,z_1+\Upsilon^P_s \overline{ e^{(t-s)A}}\bar x\right)
\<\widehat\Lambda_s^{P,C}(t-s) k, (\Upsilon^P_s Q_{t-s}(\Upsilon^P_s)^*)^{-1/2} z_1\>_{L^2_s}\caln(0,\Upsilon^P_s Q_{t-s}(\Upsilon^P_s)^*)(dz_1)ds
\end{align*}
Estimate \eqref{norm-Cderconv} follows since, arguing as in the proof of \eqref{norm-Cdernew} we have\vspace{-0.2truecm}
\begin{align*}
&
\left\vert
\int_{L^2_s} s^{-\eta}s^\eta
\hat f\left(s,z_1+\Upsilon^P_s
\overline{ e^{(t-s)A}}\bar x\right)
\<\widehat\Lambda_s^{P,C}(t-s) k, (\Upsilon^P_s Q_{t-s}(\Upsilon^P_s)^*)^{-1/2} z_1\>_H\caln(0,\Upsilon^P_s Q_{t-s}(\Upsilon^P_s)^*)(dz_1)
\right\vert\\
&\leq
s^{-\eta}\Vert\Lambda_s^{P,C}(t-s) k \Vert_{L^2_s}\cdot
\Vert s^\eta f(s,\cdot)\Vert_\infty,
\end{align*}
so the claim follows simply integrating in time
and using Hypothesis \ref{hp:smoothingextension-conv}-(ii).

The statement on continuity follows
using exactly the same argument as in
\cite[Proposition 4.50-(ii)]{FabbriGozziSwiech}.
Such continuity, together with the estimate \eqref{norm-Cderconv}, implies, recalling the last part of Definition \ref{df4:Gspaces},
that the map $g$ belongs to
$C_\gamma^{0,1,C}([0,T]\times \overline{H})$.
Now, we prove that such map belongs to
$\cals^{1,P}_{\gamma,prog}((0,T]\times \overline{H})$, i.e. we provide the two functions $\hat{g}$ and $\hat{\hat{g}}$ as required by Definition \ref{df:spaziSeta-progr}.
Indeed observe that, for $t\ge 0$ and $\bar x\in \overline{H}$ we have\vspace{-0.2truecm}
$$
g(t,\bar x)
=
\int_{0}^{t} \int_{H} f(s,z+\overline{e^{(t-s)A}}\bar x)\caln(0,Q_{t-s})(dz)ds
=
\int_{0}^{t} \int_{H} \hat f\left(s,\Upsilon^P_s
(z+\overline{e^{(t-s)A}}\bar x)\right)\caln(0,Q_{t-s})(dz)ds
$$\vspace{-0.2truecm}
$$
=
\int_{0}^{t}\int_{L^2_s}
\hat f\left(z_1+\Upsilon^P_s \overline{e^{(t-s)A}}\bar x\right)
\caln(0,\Upsilon^P_s Q_{t-s}(\Upsilon^P_s)^*)(dz_1)ds.
$$
Hence, we set\vspace{-0.2truecm}
$$
\hat{g}(t,y^P_{\bar x})
:=
\int_{0}^{t}\int_{L^2_s}
\hat f\left(z_1+\Upsilon^P_s \overline{e^{(t-s)A}}\bar x\right)
\caln(0,\Upsilon^P_s Q_{t-s}(\Upsilon^P_s)^*)(dz_1)ds.\vspace{-0.2truecm}
$$
This satisfies \eqref{eq:defSprogr} thanks to \eqref{eq:Petasemigroup}.
A similar argument holds for
$\nabla^C g$ exploiting \eqref{eq:formulader-gen-Pnewconv}.
\\
Assume now that the map
$(t,x)\mapsto t^\eta f(t,x)$ belongs to $C_b^P([0,T]\times H)$ for some $\eta\in [0,1)$. Then, by Definition \ref{df:spaziphi1} and Hypothesis \ref{ip:PC}-(iii),
we get, with the change of variable $z_1=Pz$,
$$
g(t,\bar x)
=
\int_{0}^{t}\int_{H}\bar f (s,
Pz+\overline{Pe^{(t-s)A}}\bar x)\caln(0,Q_{t-s})(dz)ds
=
\int_{0}^{t}\int_{\overline H}
\bar f\left(s,z_1+\overline{Pe^{(t-s)A}}\bar x\right)
\caln(0, PQ_{t-s}P^*)(dz_1)ds.
$$
Hence, we set
$$
\hat{g}(t,y^P_{\bar x})
:=
\int_{0}^{t}\int_{\overline H}
\bar f\left(s,z_1+\overline{Pe^{(t-s)A}}\bar x\right)
\caln(0, PQ_{t-s}P^*)(dz_1)ds.
$$
A similar argument holds for
$\nabla^C g$ using that
Hypothesis \ref{ip:NC}, and then
\eqref{eq:formulader-gen-Pold} are verified.
\hfill\qed

\section{The stochastic optimal control problem and the HJB equation}
\label{sec-HJB}

\subsection{The setting}
\label{sec-setting}

We consider the following class of stochastic optimal control problems. The state space and the control space are, respectively, the real separable Hilbert spaces $H$ and $K$.
The state equation is the following controlled SDE in $H$:
\vspace{-0.2truecm}
\begin{equation}\label{eq-common-contr}
  \begin{array}{l}
  \dis
d X(s)= AX(s)\,ds+Cu(s)\,ds +GdW(s), \qquad s\in [t,T],
\qquad X(t)=x\in H.
\end{array}
\vspace{-0.2truecm}
\end{equation}
Here $A,\,G$ and $W$ satisfy Hypothesis \ref{ip-sde-common}
and $C$ satisfies Hypothesis \ref{ip:PC}, where, in point (iii), $P:H \to H$ is a given linear and continuous map.
The process $u$  represents the control process that belongs to the set
\vspace{-0.2truecm}
\begin{equation}\label{eq:admcontr}
  \calu:=\left\{
  u:[0,T]\times \Omega \to U \subseteq K, \text{ progressively measurable}
  \right\}.
\vspace{-0.2truecm}
\end{equation}
for a suitable subset $U$ of $K$.
Equation (\ref{eq-common-contr}) is only formal and to make sense  it has to be considered in its mild formulation (using the so-called variation of constants, see e.g. \cite[Chapter 7]{DP1}). There are still some issues,
see \cite[Section 7.1]{FGFM-III}: the mild solution of \eqref{eq-common-contr} is
 still formal and it is given by\vspace{-0.2truecm}
\begin{equation}
X(s)=e^{(s-t)A}x+\int_t^s{e^{(s-r)A}}C u(r) dr +\int_t^se^{(s-r)A}Q^{1/2}dW(r),
\text{ \ \ \ }s\in[t,T]. \\
  \label{eq-mild-common}
\end{equation}
Here the first and the third term belong to $H$, thanks to Hypothesis \ref{ip-sde-common}, while the second in general does not.
By Hypothesis \ref{ip:PC}-(i) we see that the second term can be written as
$\int_t^s\overline{e^{(s-r)A}}Cu(r)dr \in \overline{H}$.
Hence, even when $x\in H$ the mild solution belongs to $\overline{H}$ but not to $H$. Moreover, still using Hypothesis \ref{ip:PC}-(i), we see that the mild solution makes sense when the initial datum belongs to $\overline{H}$.
On the other hand, thanks to Hypothesis \ref{ip:PC}-(iii), even when the initial datum is in
$\overline{H}$, the process $PX(s)$ belongs to $H$ for $s>t$ and can be written as\vspace{-0.2truecm}
\begin{equation}
PX(s)  =\overline{Pe^{(s-t)A}}x+\int_t^s\overline{Pe^{(s-r)A}}C u(r) dr +\int_t^s P e^{(s-r)A}Q^{1/2}dW(r),
\text{ \ \ \ }s\in(t,T]. \\
  \label{eq-mild-commonP}
\end{equation}\vspace{-0.2truecm}

The objective is to minimize, over all controls in $\calu$,
the following finite horizon cost ($t\in [0,T]$, $x \in H$):\vspace{-0.2truecm}
\begin{equation}\label{costoastratto-common}
J(t,x;u)=\E \left(\int_t^T \left[\ell_0(s,X^{t,x}(s))+\ell_1(u(s))\right]\,ds + \phi(X^{t,x}(T))\right)
\end{equation}
over all controls $u(\cdot) \in \calu$, under the following assumptions.
\begin{hypothesis}\label{ip-costo}
We assume that:
\begin{itemize}
  \item[(i)]
the final cost $\phi$ belongs to $C_b^P(H)$ (see Definition \ref{df:spaziphi1}).\vspace{-0.2truecm}
\item[(ii)]
for some $\eta \in [0,1)$ the {running} cost $\ell_0$ belongs to $\mathcal{S}^P_{\eta,prog}((0,T]\times \overline{H})\cap C_b((0,T]\times \overline{H})$ or is such that the map $(t,x)\mapsto t^\eta \ell_0(t,x)$ belongs to $C_b^P([0,T]\times H)$.\vspace{-0.2truecm}
  \item[(iii)]
the set $U\subset K$ is closed and bounded and the {running} cost $\ell_1:U\rightarrow \R$ is measurable and bounded from below.
\end{itemize}
\end{hypothesis}
%
%
Note that under the above assumptions, when $t\in [0,T)$, the functional $J$ is well defined also for initial states $\bar x\in \overline{H}$. We define the value function related to this control problem, as usual, as\vspace{-0.2truecm}
\begin{equation}\label{valuefunction-gen}
V(t,\bar x):= \inf_{u \in \calu}J(t,\bar x;u)\
\qquad t\in[0,T), \;\bar x \in \overline{H}
\end{equation}
while, when $t=T$, we set $V(T,x)=\phi(x)$,
for $x \in H$.

\begin{remark}\label{rm:prima-di sezHJB}
Note that in the present paper the {running} cost can depend on the
state, differently from \cite{FGFM-III}. Indeed the fixed point argument in the proof of Theorem \ref{esistenzaHJB-progr} works thanks to the extension of the``partial smoothing'' properties of the Ornstein-Uhlenbeck semigroup  proved in \cite[Section 5]{FGFM-III} to spaces of functions depending on the trajectories of the underlying semigroup $e^{tA}$, see Section \ref{sec:partsmooth-abstr-setting}, more precisely Subsections \ref{sub:OUnew} and \ref{sec:convpartsmooth-abstr-setting}.
For more on this see Remark \ref{rm:dopoproofmain}.
\end{remark}



 \subsection{The HJB equation and its solution}
 \label{sec-solHJB}
To write the HJB equation associated to our control problem we first introduce the Hamiltonian. For $p\in K$, $u \in U$, we define the current value Hamiltonian $H_{CV}:K\times U \to \R$ as
$H_{CV}(p\,;u):=\<p,u\>_{K}+\ell_1(u)$
and the (minimum value) Hamiltonian by\vspace{-0.2truecm}
\begin{equation}\label{psi1-gen}
H_{min}(p)=\inf_{u\in U}H_{CV}(p\,;u),
 \end{equation}
The HJB equation associated to the stochastic optimal control problem that we have just presented is, formally,\vspace{-0.2truecm}
\begin{equation}\label{HJBformale-common}
  \left\{\begin{array}{l}\dis
-\frac{\partial v(t,x)}{\partial t}=\call [v(t,\cdot)](x) +\ell_0(t,x)+
H_{min} (\nabla^C v(t,x)),\qquad t\in [0,T),\,
x\in H,\\[2mm]
\dis v(T,x)=\phi(x),\quad x\in H.
\end{array}\right.
\end{equation}
{Here the operator $\call$
is the generator of the
semigroup $(R_{t})_{0\leq t\leq T}$ (see (\ref{ornstein-sem-gen}))
associated to the process $Z^{t,x}$ given in \eqref{ornstein-gen}-\eqref{ornstein-mild-gen}.
Such generator is, formally, defined by taking the drift and the diffusion coefficient in \eqref{ornstein-gen}: \vspace{-0.2truecm}
\begin{equation}\label{eq:ell-gen}
 \call[f](x)=\frac{1}{2} Tr \; GG^*\nabla^2 f(x)
+ \< Ax,\nabla f(x)\>
\end{equation}}
{see \cite{DP1}[Paragraph 9.3] for more details.}
{
\begin{remark}
\label{rm:generatornew}
Note that, in the form \eqref{eq:ell-gen}, the generator $\call$ makes sense only when $x \in D(A)$, $f$ is twice differentiable in such $x$ and $GG^*\nabla^2 f(x)$ is a nuclear operator.
In some cases it is useful to rewrite it in the form
\begin{equation}\label{eq:ell-genA*}
 \call[f](x)=\frac{1}{2} Tr \; GG^*\nabla^2 f(x)
+ \< x,A^*\nabla f(x)\>.
\end{equation}
In this way the above expression of $\call$ is defined for all $x \in H$ at the price of asking more regularity on $f$, namely that $A^*\nabla f(x)$ be well defined for all $x \in H$.
\\
To give precise sense to $\call$ as the generator of the semigroup $R_t$ one does the following: first one defines an operator $\call_0$ using \eqref{eq:ell-genA*} in a suitable domain of regular functions; then one shows that such operator $\call_0$ is closable and, taking its closure in a proper ``weak way'' (recall that the semigroup $(R_t)$ is not strongly continuous) one defines $\call$. This is explained e.g.
in \cite{CerraiSem} and in \cite[Section 4.3.1 and Appendix B]{FabbriGozziSwiech}.
\end{remark}
}
{Using, formally, the method of variation of constants, and the extensions of our operators to the larger space $\overline{H}$, we arrive at the following definition of mild solution for the HJB equation.}

\begin{definition}\label{defsolmildHJB}
A function $v$
is a mild solution of the HJB equation (\ref{HJBformale-common}) if the following are satisfied:
\begin{enumerate}
\item for some $\alpha \in [0,1)$, $v(T-\cdot, \cdot)\in C^{0,1,C}_{{\alpha}}\left([0,T]\times \overline H\right)$;
\item the integral equation\vspace{-0.2truecm}
\begin{equation}\
v(t,\bar x) =R_{T-t}[\phi](\bar x)+\int_t^T R_{s-t}\left[
H_{min}(\nabla^C v(s,\cdot))\right](\bar x)\;ds
+\int_t^T R_{s-t}\left[
\ell_0(s,\cdot)\right](\bar x)\; ds,\qquad t\in [0,T),\
\bar x\in \overline{H},
\label{solmildHJB-common}
\end{equation}
is satisfied on $(0,T]\times H$ together with the final condition $V(T,x)=\phi (x)$ for $x \in H$.
\end{enumerate}
\end{definition}
\begin{remark}\label{rm:crescitapoli-HJB}
Since functions in  $C^{0,1,C}_{\alpha}\left([0,T]\times \overline
H\right)
$ are bounded (see Definition \ref{df4:Gspaces}),
the above Definition \ref{defsolmildHJB} requires, among other properties, that a mild solution is continuous and bounded up to $T$.
This constrains the assumptions on the data, e.g. it implies that the final datum $\phi$ must be continuous and bounded.
We may change this requirement in the above definition asking only polynomial growth in $x$. Our main results will remain true with straightforward modifications using the spaces with polynomial growth as in \cite[Chapter 4]{FabbriGozziSwiech}.
\end{remark}
\begin{remark}
\label{rm:dopoproofmain}
Notice that, even if we assume that $\ell_0(t,x)\in C_b^P([0,T]\times H)$, it turns out that the last convolution term in \eqref{solmildHJB-common}, i.e.\vspace{-0.2truecm}
$$
\displaystyle\int_0^t R_{t-s}\ell_0(s,\cdot)(x)\; ds=\int_0^t \int_{H}^{}\bar{\ell_0}(s,P(z+e^{(t-s)A}x))
\mathcal{N}(0,Q_t)(dz)\; ds
$$
belongs to a suitable subspace of  $\calc^P_A((0,T]; H)$, i.e. such term is, for all $\bar x\in \overline{H}$, a function of the path $y^P_{\bar x}$. So we have to look for the solution of the HJB equations in a subspace of $\calc^P_A((0,T]; H)$, this subspace turns out to be $\cals^{1,P}_{\gamma,prog}([0,T]\times \overline H)$ (or $\cals^{1,P}_{\gamma}([0,T]\times \overline H)$); for this reason we have studied the smoothing properties of the semigroup $R_t$ on these spaces in Subsections \ref{sub:OUnew} and \ref{sec:convpartsmooth-abstr-setting}.
\end{remark}

We now prove existence and uniqueness of a mild solution of the HJB equation \eqref{HJBformale-common} by a fixed point argument in the space $\cals^{1,P}_{\gamma,prog}([0,T]\times \overline H)$ (or $\cals^{1,P}_{\gamma}([0,T]\times \overline H)$), see Definitions \ref{df:spaziSeta-progr} and \ref{df:spaziSeta}.

\begin{theorem}
\label{esistenzaHJB-progr}
Let Hypotheses \ref{ip-sde-common}, \ref{ip:PC}, \ref{ip:NC}, \ref{hp:L2}, {\ref{hp:smoothingextension-conv}} and \ref{ip-costo} hold true.
Then the HJB equation (\ref{HJBformale-common})
admits a mild solution $v$ according to Definition \ref{defsolmildHJB}; $v$ is unique among the functions $w$ such that $w(T-\cdot,\cdot)\in\cals^{1,P}_{\gamma, prog}([0,T]\times \overline H)$ and it satisfies, for suitable $C_T>0$, the estimate\vspace{-0.2truecm}
\begin{equation}\label{eq:stimavmainteo}
\Vert v(T-\cdot,\cdot)\Vert_{C^{0,1,C}_{{\gamma}}}\le C_T\left(\Vert\bar\phi \Vert_\infty
+\Vert\bar\ell_0 \Vert_\infty \right).
\end{equation}
Moreover, let also Hypothesis \ref{hp:smoothingextension} hold true. Then the mild solution $v$ is unique among the functions $w$ such that $w(T-\cdot,\cdot)\in\cals^{1,P}_{\gamma}([0,T]\times \overline H)$.
Finally, in such a case, the statement remains true if, differently from Hypothesis \ref{ip-costo}, the final cost belongs to $\mathcal{S}^P_\infty(\overline{H})$ and the {running} cost $\ell_0$ belongs, for some $\eta \in [0,1)$ to
$\mathcal{S}^P_{\eta}((0,T]\times \overline{H})\cap C_b((0,T]\times \overline{H})$. \end{theorem}
\dim
We prove the first part of the Theorem.
As a first step we consider the case $T_1\geq T$, where $T_1$ is given in Hypothesis \ref{hp:smoothingextension-conv}.
We prove existence and uniqueness of a solution in $\cals^{1,P}_{\gamma,prog}([0,T]\times \overline H)$, by using a fixed point argument in it. To this aim, first we rewrite (\ref{solmildHJB-common}) in a forward way. Namely if $v$ satisfies \myref{solmildHJB-common} then, setting $w(t,x):=v(T-t,x)$ for any
$(t,x)\in[0,T]\times H$, we get that $w$ satisfies $w(0,x):=\phi(x)$, $x \in H$ and\vspace{-0.2truecm}
\begin{equation}\label{eq:wT1}
  w(t,\bar x) =R_{t}[\phi](\bar x)+\int_0^t R_{t-s}[
H_{min}(
\nabla^C w(s,\cdot))+\ell_0(s,\cdot)
](\bar x)\; ds,\qquad t\in (0,T],\
\bar x\in\overline H,\vspace{-0.2truecm}
\end{equation}
which is the mild form of the forward HJB equation\vspace{-0.2truecm}
\begin{equation}\label{HJBformaleforward-common}
  \left\{\begin{array}{l}\dis
\frac{\partial w(t,\bar x)}{\partial t}=\call [w(t,\cdot)](\bar x) +\ell_0(t,\bar x)+
H_{min} (\nabla^{C}w(t,\bar x)),\qquad t\in [0,T],\;
\bar x\in \overline H,\\[2mm]
\dis w(0,x)=\phi(x), \qquad x \in H.
\end{array}\right.
\end{equation}
Define the map $\Gamma$ on $\cals^{1,P}_{\gamma,{prog}}([0,T]\times \overline H)$  by setting, for $g\in \cals^{1,P}_{\gamma,{prog}}([0,T]\times \overline H)$, $\Gamma(g)(0,x):=\phi(x)$, $x \in H$, and\vspace{-0.2truecm}
\begin{equation}\label{Gamma}
\Gamma(g)(t,\bar x):=R_{t}[\phi](\bar x)+\int_0^t R_{t-s}[
H_{min}(\nabla^C g(s,\cdot))](\bar x){\,ds}+\int_0^t R_{t-s}[\ell_0(s,\cdot)](\bar x)\; ds,
\qquad t\in (0,T], \;\bar x \in \overline{H}.
\end{equation}
By Hypothesis \ref{ip-costo}-(i) and the last part of Proposition \ref{prop:partsmoothnew} the first term belongs to $\cals^{1,P}_{\gamma,prog}([0,T]\times \overline H)$.
Moreover, Hypothesis \ref{ip-costo}-(ii) and Proposition \ref{prop:partsmoothnew-conv}
imply that the third term belongs to $\cals^{1,P}_{\gamma,prog}([0,T]\times \overline H)$, too.
Furthermore, Hypothesis \ref{ip-costo}-(iii) implies that the Hamiltonian $H_{min}$ is Lipschitz continuous; since
$\nabla^C g$ belongs to
$\cals^{P}_{\gamma,prog}((0,T]\times \overline H;K)$ this implies that also the second term belongs to $\cals^{1,P}_{\gamma,prog}([0,T]\times \overline H)$.
Hence $\Gamma$
is well defined in $\cals^{1,P}_{\gamma,prog}([0,T]\times \overline H)$ with values in $\cals^{1,P}_{\gamma,prog}([0,T]\times \overline H)$.
\newline
As already stated in Lemma \ref{lemma:SP}, $\cals^{1,P}_{\gamma,prog}([0,T]\times \overline H)$ is a closed subspace of $C^{0,1,C}_{{\gamma}}([0,T]\times \overline H)$, and so if $\Gamma$ is a contraction, by the Contraction Mapping Principle there exists a unique (in $\cals^{1,P}_{\gamma,prog}([0,T]\times \overline H)$) mild solution of (\ref{HJBformale-common}).
Hence, we take
$g_1,g_2 \in \cals^{1,P}_{\gamma,prog}([0,T]\times \overline H)$ and we evaluate\vspace{-0.2truecm}
$$
\Vert \Gamma (g_1)-\Gamma (g_2)\Vert_{\cals^{1,P}_{\gamma}([0,T]\times \overline{H})}
=\Vert \Gamma(g_1)-\Gamma (g_2)\Vert_{C^{0,1,C}_{\gamma}([0,T]\times \overline H)}.$$
For every $(t,\bar x)\in (0,T]\times \overline H$ we get, calling $L$ the Lipschitz constant of $H_{min}$,
\begin{align*}
  \vert \Gamma (g_1)(t,\bar x)- \Gamma(g_2)(t,\bar x) \vert& =\left\vert \int_0^t R_{t-s}\left[H_{min}\left(\nabla^C g_1(s,\cdot)\right)
 -H_{min}\left(\nabla^C g_2(s,\cdot)\right)\right](\bar x)ds\right\vert\\
 &\leq \int_0^t s^{-\gamma} L \sup_{\bar y \in \overline{H}}\vert s^{{\gamma}}\nabla^C (g_1-g_2)(s,\bar y)\vert ds
  \leq \frac{1}{1-\gamma} Lt^{{1-\gamma}}\Vert g_1-g_2 \Vert_{C^{0,1,C}_{\gamma}([0,T]\times \overline{H})}.
\end{align*}
Similarly\vspace{-0.2truecm}
\begin{align*}
t^{{\gamma}}\vert \nabla^C\Gamma (g_1)(t,\bar x) &- \nabla^C\Gamma(g_2)(t,\bar x) \vert =
t^{{\gamma}}\left\vert \nabla^C\int_0^t  R_{t-s}\left[H_{min}
\left(\nabla^C g_1(s,\cdot)\right)-H_{min}\left(\nabla^C
g_2(s,\cdot)\right)\right](x)ds\right\vert\\
& \leq t^{{\gamma}} L \Vert g_1-g_2 \Vert_{C^{0,1,C}_{\gamma}([0,T]\times \overline{H})} \int_0^t (t-s)^{-\gamma}
s^{-\gamma} ds
\leq L\beta(\gamma,\gamma) t^{1-\gamma}\Vert g_1-g_2 \Vert_{C^{0,1,C}_{\gamma}([0,T]\times \overline{H})}.
\end{align*}
Hence, if $L[\beta(\gamma,\gamma)+(1-\gamma)^{-1}] T^{1-\gamma}\le 1/2$, we get\vspace{-0.2truecm}
\begin{equation}\label{stima-contr}
 \left\Vert \Gamma (g_1)-\Gamma(g_2)\right
 \Vert_{C^{0,1,C}_{\gamma}([0,T]\times \overline{H})}\leq \frac12
\left\Vert g_1-g_2\right
\Vert_{C^{0,1,C}_{\gamma}([0,T]\times \overline{H})}.
\end{equation}
Hence the map $\Gamma$ is a contraction in $\cals^{1,P}_{\gamma,prog}([0,T]\times \overline{H})$ and, if we denote by $w$ its unique fixed point, then $v:=w(T-\cdot,\cdot)$ turns out to be the wanted unique mild solution of the HJB equation (\ref{HJBformale-common}), according to Definition \ref{defsolmildHJB}.

Since the value of $T$ is independent of the initial data, the case of generic $T>0$ follows by dividing the interval $[0,T]$ into a finite number of subintervals of length $\delta$ sufficiently small, or equivalently, as done in \cite{Mas}, by taking an equivalent norm with an adequate exponential weight, such as\vspace{-0.2truecm}
\[
 \left\Vert f\right\Vert _{\mu,C^{0,1,C}_{\gamma}}
=\sup_{(t,\bar x)\in[0,T]\times \overline{H}}
\vert e^{\mu t}f(t,x)\vert+
\sup_{(t,x)\in (0,T]\times \overline{H}}  e^{\mu t}t^{\alpha}
\left\Vert \nabla^C f\left(  t,x\right)  \right\Vert _{\call(K,H)}.
\]
This concludes the proof for the case $T_1\geq T$.

\smallskip

Let now $T_1<T$ and assume, without loss of generality, that\vspace{-0.2truecm}
\begin{equation}\label{eq:condT1}
L[\beta(\gamma,\gamma)+(1-\gamma)^{-1}]  T_1^{1-\gamma}\le 1/2.
\end{equation}
Hence, we can write \eqref{HJBformale-common} as an equation for $w$ as done in \eqref{HJBformaleforward-common}
and, by the Contraction Mapping Theorem we find a unique solution
$w_1\in \cals^{1,P}_{\gamma,prog}([0,T_1]\times \overline{H})$.
Now we repeat the contraction argument in the time interval $[T_1/2, 3T_1/2]$.
To do this we first observe that, if $w$ solves
\eqref{eq:wT1}, then, for $t\in (0,T]$, $\bar x \in \overline{H}$\vspace{-0.2truecm}
\begin{align*}
w(t,\bar x) &=R_{t-T_1/2}[R_{T_1/2}[\phi]]\!+\!
\int_0^{T_1/2} \!\!\!\! R_{T_1-s}[H_{min}(\nabla^C w(s,\cdot))+\ell_0(s,\cdot)
]{(\bar x)}\; ds\!+\!\int_{T_1/2}^{t} \!\!\!\! R_{t-s}[H_{min}(\nabla^C w(s,\cdot))+\ell_0(s,\cdot)](\bar x)\; ds\\
&=R_{t-T_1/2}[w(T_1/2, \cdot)]
(\bar x)+\int_{T_1/2}^{t} R_{t-s}[H_{min}
(\nabla^C w(s,\cdot))+\ell_0(s,\cdot)]
(\bar x)\; ds
\end{align*}
This means that we aim to solve, in the space $\cals^{1,P}_{prog} ([T_1/2,3T_1/2]\times \overline{H})$
the mild equation\vspace{-0.2truecm}
\begin{equation}
\label{eq:iteratacontrazione}
w(t,\bar x) =R_{t-T_1/2}[w_1(T_1/2, \cdot)](\bar x)+\int_{T_1/2}^{t} R_{t-s}[H_{min}(\nabla^C w(s,\cdot))+\ell_0(s,\cdot)](\bar x)\; ds,
\qquad t\in [T_1/2, 3T_1/2],\;
\bar x\in\overline H,\vspace{-0.2truecm}
\end{equation}
where now the initial datum is given by
$w_1(T_1/2, \cdot)$.
We first prove that, calling $\Gamma_1(g)$ the right hand side of \eqref{eq:iteratacontrazione} with $g$ at the place of $w$, $\Gamma_1$ maps
$\cals^{1,P}_{prog} ([T_1/2,3T_1/2]\times \overline{H})$ into itself.
\newline We notice that
$w_1(T_1/2, \cdot) \notin C^P_b(H)$, indeed it does not depend on the state $x$ only through $Px$, but it depends on the whole trajectory
$y^P_{x}(\cdot \wedge (T_1/2))$. So we can not apply the smoothing given by Hypothesis \ref{ip:NC}, but we can take into account the fact that $w_1$ solves in mild sense equation \eqref{HJBformaleforward-common}, so $w_1\in \cals^{1,P}_{\gamma,prog} ([0,T_1]\times \overline{H})$ and in particular $w_1(T_1/2, \cdot)$ is bounded, continuous, differentiable and with a bounded derivative.
We use these facts to show that the map
$(t,\bar x) \mapsto R_{t-T_1/2}[w_1(T_1/2,\cdot)](\bar x)$ belongs to $\cals^{1,P}_{prog} ([T_1/2,3T_1/2]\times \overline{H})$.

We notice that since $w_1\in\cals^{1,P}_{\gamma,prog}
([0,T_1]\times \overline{H})$ there exists
$\hat f:[0,T_1]\times C^P_A((0,+\infty);H)\to\R$  such that $w_1(T_1/2, \bar x)=\hat f(T_1/2,y^P_{\bar x}(\cdot \wedge (T_1/2) )$, for all $\bar x \in \overline{H}$ so\vspace{-0.2truecm}
\begin{align*}
R_{t-T_1/2}[w_1(T_1/2, \cdot)](\bar x)=
\int_H \hat f\left(T_1/2,
\Upsilon^P_{T_1/2}
(z+\overline{e^{(t-T_1/2)A}}\bar x)\right)
\caln(0,Q_{t-T_1/2})(dz)\vspace{-0.2truecm}
\end{align*}
and this implies that the map
$(t,\bar x) \mapsto
R_{t-T_1/2}[w_1(T_1/2,\cdot)](\bar x)$
belongs to $\cals^{P}_{prog} ([T_1/2,3T_1/2]\times \overline{H})$.
Moreover, again since $w_1\in\cals^{1,P}_{\gamma,prog} ([0,T_1]\times \overline H)$ there exists $\hat{\hat{ f}}:[0,T_1]\times C^P_A((0,+\infty);H)\to\R$  such that $\nabla^Cw(T_1/2, x)=\hat{\hat{ f}}(T_1/2,y^P_{x}(\cdot \wedge (T_1/2))$, so\vspace{-0.2truecm}
\begin{align*}
\nabla^C R_{t-T_1/2}[w_1(T_1/2, \cdot)](\bar x)= R_{t-T_1/2}[\nabla^C w_1(T_1/2, \cdot)](\bar x)=
\int_H \hat{\hat {f} }\left(T_1/2,
\Upsilon^P_{T_1/2}(z+\overline{e^{(t-T_1/2)A}} x)\right)
\caln(0,Q_{t-T_1/2})(dz)\vspace{-0.2truecm}
\end{align*}and this implies that the map
$(t,\bar x) \mapsto
R_{t-T_1/2}[w_1(T_1/2,\cdot)](\bar x)$
belongs to $\cals^{P}_{prog} ([T_1/2,3T_1/2]\times \overline{H})$.

To see that, when $g\in
\cals^{1,P}_{prog} ([T_1/2,3T_1/2]\times \overline{H})$,
the convolution term $\int_{T_1/2}^t R_{t-s}[H_{min}(\nabla^C g(s,\cdot))+\ell_0(s,\cdot)](x)\; ds$
belongs to $\cals^{1,P}_{prog} ([T_1/2,3T_1/2]\times \overline{H})$
we simply apply Proposition
\ref{prop:partsmoothnew-conv} taking, in its statement and namely in formula \eqref{eq:convgdef}, $f$ equal to zero between $0$ and $T_1/2$,
and $f$ equal to $H_{min}(\nabla^C g(s,\cdot))+\ell_0(s,\cdot)$ otherwise.

Hence we conclude that $\Gamma_1$ is well defined in $\cals^{1,P}_{0,prog}([T_1/2, 3T_1/2]\times \overline H)$ with values in $\cals^{1,P}_{prog}([T_1/2, 3T_1/2]\times \overline H)$.

Now we prove the contraction property.
We take
$g_1,g_2 \in \cals^{1,P}_{prog}(
[T_1/2, 3T_1/2]\times \overline H)$. and we evaluate\vspace{-0.2truecm}
$$
\Vert \Gamma_1 (g_1)-\Gamma_1 (g_2)\Vert_{\cals^{1,P}_{prog}([T_1/2, 3T_1/2]\times \overline{H})}
=\Vert \Gamma_1(g_1)-\Gamma_1 (g_2)\Vert_{C^{0,1,C}_{0}([T_1/2, 3T_1/2]\times \overline H)}.$$
For every $(t,\bar x)\in [T_1/2,3T_1/2]\times \overline H$ we get, calling $L$ the Lipschitz constant of $H_{min}$,\vspace{-0.2truecm}
\begin{align*}
  \vert \Gamma_1 (g_1)(t,\bar x)- \Gamma_1(g_2)(t,\bar x)\vert &=
  \left\vert \int_{T_1/2}^t R_{t-s}\left[H_{min}\left(\nabla^C g_1(s,\cdot)\right)-H_{min}\left(\nabla^C g_2(s,\cdot)\right)\right](\bar x)ds\right\vert
 \\
  &\le L(t-T_1/2)\Vert g_1-g_2\Vert_{C^{0,1,C}_{0}([T_1/2, 3T_1/2]\times \overline H)}
  \le LT_1{\Vert g_1-g_2}\Vert_{C^{0,1,C}_{0}([T_1/2, 3T_1/2]\times \overline H)}
\end{align*}
and\vspace{-0.6truecm}
\begin{align*}
&\vert \nabla^C\Gamma_1 (g_1)(t,\bar x) - \nabla^C\Gamma_1(g_2)(t,\bar x) \vert
=
\left\vert \nabla^C\int_{T_1/2}^t  R_{t-s}\left[H_{min}
\left(\nabla^C g_1(s,\cdot)\right)-H_{min}\left(\nabla^C
g_2(s,\cdot)\right)\right](x)ds\right\vert\\
&+L\left(\int_{T_1/2}^t (t-s)^{-\gamma} ds\right)
\Vert g_1-g_2 \Vert_{C^{0,1,C}_{0}([[T_1/2, 3T_1/2]]\times \overline{H})}
\leq
L(1-\gamma)^{-1}T_1^{1-\gamma}
\Vert g_1-g_2 \Vert_{C^{0,1,C}_{\gamma}([T_1/2, 3T_1/2]\times \overline{H})}.
\vspace{-0.3truecm}
\end{align*}
Now we observe that condition \eqref{eq:condT1} implies that $L[T_1+(1-\gamma)^{-1}T_1^{1-\gamma}]\le 1/2$.
This implies that $\Gamma_1$ is a contraction
on $\cals^{1,P}_{prog}([T_1/2,3T_1/2]\times \overline{H})$ so we get existence and uniqueness of a solution $w_2$ in such space.
This solution must coincide with $w_1$ on the interval $[T_1/2,T_1]$ by repeating the same contraction argument in such time interval.
Hence, calling $w_3$ the concatenation of $w_1$ and $w_2$ we get $w_3$ is the unique solution of
\eqref{eq:wT1} in $\cals^{1,P}_{\gamma,prog}([0,3T_1/2]\times \overline{H})$.
%
The conclusion of the proof follows simply iterating the procedure to cover all the interval $[0,T]$.

Next we briefly sketch the proof of the second part of the Theorem, we assume here {Hypothesis \ref{hp:smoothingextension}}.
Define the map $\Gamma$ as in \eqref{Gamma}, considering now $g\in \cals^{1,P}_{\gamma}([0,T]\times \overline H)$.
By Hypothesis \ref{ip-costo}-(i) and Proposition \ref{prop:partsmoothnew} we have that $R_{t}[\phi]\in \cals^{1,P}_{\gamma}([0,T]\times \overline H)$.
Moreover, Hypothesis \ref{ip-costo}-(ii) and Proposition \ref{prop:partsmoothnew-conv}
imply that the term $\int_0^t R_{t-s}[\ell_0(s,\cdot)]\,ds$ belongs to $\cals^{1,P}_{\gamma}([0,T]\times \overline H)$, too.
Furthermore, as argued before the Hamiltonian $H_{min}$ is Lipschitz continuous and so with at most linear growth.
These facts, together with the fact that $\nabla^C g$ belongs to $\cals^{P}_{\gamma}((0,T]\times \overline H;K)$ imply that $\Gamma$ is well defined in $\cals^{1,P}_{\gamma}([0,T]\times \overline H)$ with values in $\cals^{1,P}_{\gamma}([0,T]\times \overline H)$.
\newline As already stated in Lemma \ref{lemma:SP}, $\cals^{1,P}_{\gamma}([0,T]\times \overline H)$ is a closed subspace of $C^{0,1,C}_{{\gamma}}([0,T]\times \overline H)$, and so as in the first part of the Theorem it suffices to prove that $\Gamma$ is a contraction.
This is done taking $g_1,g_2 \in \cals^{1,P}_{\gamma}([0,T]\times \overline H)$ and repeating in the same way the steps to arrive at \eqref{stima-contr}.
As above, the case of generic $T>0$ follows by dividing the interval $[0,T]$ into a finite number of subintervals of sufficiently small length, or equivalently, as done in \cite{Mas}, by taking an equivalent norm with an adequate exponential weight.
\qed



{
\begin{remark}\label{rem:esistenzaHJB-LocLipEtAl}
Looking at the previous literature on mild solutions
it seems possible to generalize Theorem \ref{esistenzaHJB-progr} in various directions (with increasing order of difficulty):
\vspace{-0.3truecm}
\begin{itemize}
\item First of all, as already said in Remark \ref{rm:crescitapoli-HJB}, in the case when the data $\phi$ and $\ell$ have polynomial growth in the state variable $x$. In this case it seems quite straightforward to apply the arguments used in \cite{Ce95} and
      \cite[Chapter 4, Sections 4.3-4.4-4.5]{FabbriGozziSwiech}.
\vspace{-0.3truecm}
  \item Second, when the control space $U$ is not bounded. In this case the Hamiltonian $H_{min}$ is not globally Lipschitz. However, in such case it seems possible to apply the techniques used in \cite{G2}, \cite{FMloclip}, \cite[Section 4.7]{FabbriGozziSwiech} to prove suitable apriori estimates which allow to get existence and uniqueness of the HJB equation.
\vspace{-0.3truecm}
  \item Third when the operator $G$ is state dependent.
  As noted in \cite{PZ} (see also
  \cite[Section 4.3.3]{FabbriGozziSwiech}) under suitable assumptions on $G$, the corresponding transition semigroup still enjoys smoothing properties. It is then reasonable to prove that some partial smoothing holds in this case and use it to treat HJB equations also in this case.
\vspace{-0.3truecm}
\end{itemize}
Even if the above generalizations seem doable, their careful proof would require another nontrivial paper and are left for future research.
\end{remark}
}

\vspace{-0.2truecm}


\section{Examples}
\label{sec:verifica_ipotesi}

In this Section we provide the details on the examples mentioned in the introduction and prove that the crucial assumptions \ref{ip:PC}, \ref{ip:NC}, \ref{hp:L2} and \ref{hp:smoothingextension-conv} are verified in such motivating problems: stochastic heat equations with boundary control and stochastic equations with delay in the control; for stochastic heat equations with boundary control we also prove that Hypothesis \ref{hp:smoothingextension} is verified.
We underline the fact that we are not able to prove in an abstract general setting that Hypothesis \ref{ip:NC} implies Hypothesis \ref{hp:smoothingextension-conv}, so in each model considered we give references where Hypothesis \ref{ip:NC} is verified, and we verify Hypothesis \ref{hp:smoothingextension-conv}. Only in Section \ref{SSE:HEATEQUATION}  also Hypothesis \ref{hp:smoothingextension} is verified, and, as \ref{hp:smoothingextension-conv}, it must be proved directly because it is not implied by \ref{ip:NC}.

\vspace{-0.3truecm}

\subsection{Heat equations with boundary control}
\label{SSE:HEATEQUATION}

We consider a stochastic heat equation {in an open connected set with smooth boundary $\calo\subseteq \R^d$ ($d=1,2,3$) with Dirichlet boundary conditions and with boundary control.}  We reformulate it as a stochastic evolution equation in the space $H:=L^2(\calo)$, following the path outlined in \cite[Appendix C]{FabbriGozziSwiech}, see also \cite[section 3.1]{FGFM-III}.
The state equation is
\vspace{-0.2truecm}
\begin{equation}\label{eqDiri-abstr-contr}
d X(s)= {A}X(s)\,ds +
( -{A})D u(s)\,dt+{G}dW(s),
\qquad X(t)=x.
\vspace{-0.2truecm}
\end{equation}
where $W(\cdot)$ is a cylindrical noise in $H$ and the initial condition $x(\cdot)\in H$, the operator ${A}$ in $H$ is defined by
\vspace{-0.2truecm}
$$
\cald({A})=H^2(\calo)\cap H^1_0(\calo)
 \qquad
{A}y = \Delta y
 {\rm \;\; for\; \;} y\in \cald({A}),
\vspace{-0.2truecm}
 $$
 where $H^2(\calo)$ and $H^1_0(\calo)$ are the usual Sobolev spaces.

The operator ${A}$ is self-adjoint and diagonal with strictly negative eigenvalues $\{-\lambda_n\}_{n\in \N}$, with $\lambda_n\sim n^{2/d}$ as $n \to +\infty$, we denote by $\{e_n\}_{n\in \N}$ the correspnding orthonormal basis of eigenvectors of ${A}$
Notice that we also have $
\overline{e^{t{A}}}=e^{t\overline{{A}}}$ where $\overline{{A}}:\overline H \to  H
$ is the canonical extension of ${A}$ to $\overline H$.
We define
\vspace{-0.2truecm}
\begin{equation}\label{notazioneB}
{B}:=( -{A})D,
\vspace{-0.2truecm}
\end{equation}
{where $D:L^2(\partial \calo)\rightarrow H$ is the linear trace operator and it is defined by setting
 $D a=f$ where $f$ is the unique solution of the Dirichlet problem
$\Delta f(\xi)=0$, ($\xi\in \calo$),
and $f(\xi)=a(\xi)$, ($\xi\in \partial\calo$).
The operator ${B}$, defined in $K=L^2(\partial\calo)$,
does not take values in $H=L^2(\calo)$.}
Indeed for all $\varepsilon>0$,
$D\in \call\left(K, \cald((-{A})^{1/4-\varepsilon})\right)$
(see again \cite[Appendix C]{FabbriGozziSwiech}).
So
\vspace{-0.2truecm}
$$
{B}=(-{A})^{3/4+\varepsilon}(-{A})^{1/4-\varepsilon}D:
K \rightarrow \cald((-{A})^{-3/4-\varepsilon}),
\vspace{-0.2truecm}
$$
where for $\gamma>0$ we denote by $\cald((-{A})^{-\gamma})$
the completion of $H$ with respect to the norm
$|\cdot|_{-\gamma}=|(-{A})^{-\gamma}\cdot|_H$. Hence
for all $\eps \in (0,1/4)$, we have
$B_0\in \call\left(K,\cald((-{A})^{-3/4-\varepsilon})\right)$ and we
may say that ${A}$ is unbounded on $H$, in the sense that
its image is not contained in $H$ but in
$\cald((-{A}^{-3/4-\varepsilon})=
H^{-3/2-2\varepsilon}(\calo)$.
We fix a small $\eps\in (0,1/4)$ and we set $\overline{H}:=\cald((-{A})^{-3/4-\varepsilon})$, see also \cite[section 3.1]{FGFM-III}. Here we can extend immediately $e^{t{A}}$ to
$\overline{e^{t{A}}}:\overline H \to  H$
thanks to its eigenfunction expansion.
The unique mild solution of \eqref{eqDiri-abstr-contr} (which exists thanks e.g. to
\cite[Theorem 1.141]{FabbriGozziSwiech}) is denoted by $X(\cdot;t,x,u)$ and is
\vspace{-0.2truecm}
$$
X(s;t,x,u)=
e^{(s-t){A}}x+\int_t^s e^{(s-r){A}}{A} u(r) dr +\int_t^se^{(s-r){A}}{G}dW(r)
,\text{ \ \ \ }s\in[t,T].
\vspace{-0.2truecm}
$$
Let us consider the uncontrolled version of equation \eqref{eqDiri-abstr-contr}, i.e. let $u\equiv 0$ in \eqref{eqDiri-abstr-contr}, and let us denote by $X(s;t,x)$ the uncontrolled process; the transition semigroup related is defined on functions $\phi \in B_b(H,\R)$, as
\vspace{-0.2truecm}
\begin{equation}\label{transition-heat}
R_{s}[\phi](x)=\E \phi(X(s;0,x)), \qquad 0\leq s\leq T.
\vspace{-0.2truecm}
\end{equation}
Now we consider the optimal control problem related to the stochastic heat equation with boundary control in its abstract reformulation \eqref{eqDiri-abstr-contr}, also in order to introduce the class of function on which we study the partial smoothing, defined by means of the operator $P$. For any given $t\in [0,T]$ and $x \in H$, the objective is to minimize, over all control strategies in $\calu$, the following finite horizon cost:
\vspace{-0.2truecm}
\begin{equation}\label{costoastratto}
J(t,x;u)=\E \left[\int_t^T \left[\ell_0(s,X(s;t,x,u)))+\ell_1(u(s))\right]\,ds + \phi(X(T;t,x,u))\right].
\vspace{-0.2truecm}
\end{equation}
Notice that towards \cite{FGFM-III} we are able to consider also running cost depending on the state due to the new results contained in Section \ref{sec:convpartsmooth-abstr-setting}.
We assume the following:
\begin{hypothesis}\label{hp:BCcost}
\begin{itemize}
\item[]
  \item[(i)]$\ell_0:[0,T]\times H\rightarrow \R$ is measurable and $\forall t\in[0,T], $ $\ell_0(t,\cdot)$ is such that for a suitable finite linearly independent set $\{\v_1,\dots,\v_N\}\subseteq \cald((-{A})^{\eta})$ (with $\eta>1/4+\eps $) and a suitable $\bar\ell_0 \in C_b([0,T]\times\R^n)$ we have
$\ell_0(t,x)=\bar\ell_0\left(t,\<x,\v_1\>_H,\dots,
\<x,\v_N\>_H\right)$.
  \item[(ii)] $\ell_1:U\to \R$ is measurable and bounded from below
  \item[(iii)] $\phi:H\to \R$ is such that, for the finite set $\{\v_1,\dots,\v_N\}\subseteq \cald((-{A})^{\eta})$ introduced at point (i) and a suitable $\bar\phi \in C_b(\R^n)$ we have
$\phi(x)=\bar\phi\left(\<x,\v_1\>_H,\dots,\<x,\v_N\>_H
\right)$.
\end{itemize}
\end{hypothesis}

Notice that by Hypothesis \ref{hp:BCcost} we are considering special functions defined on $D((-{A}^{\eta}),\,\eta>\frac{1}{4}$, defined by means of  the projection on the span of $\<\v_1,...,\v_N\>$, denoted by $P$ and defined by
\vspace{-0.3truecm}
\begin{equation}\label{P-n-gen}
P:H\rightarrow H,\quad Px=\sum_{i=1}^N\<x,\v_i\>\v_i, \, \forall x \in H.
\vspace{-0.2truecm}
\end{equation}
Note that $P=P^*$, that
$\overline{Pe^{t{A}}}
=P\overline{e^{t{A}}}=Pe^{t\overline{{A}}}$ for $t>0$, and that
$P_\eta:=(-{A})^{\eta}P$ is a continuous operator on $H$.
Moreover we set
\vspace{-0.1truecm}
\begin{equation}\label{barQ_t-heat}
\overline Q_t:= PQ_tP=P (-A)^{-1-\beta}(I-e^{2tA})P.
\vspace{-0.2truecm}
\end{equation}
%
%

We now verify that our assumptions are true in this case.
First of all, Hypotheses \ref{ip-sde-common} \ref{ip:PC}, \ref{ip:NC} are proved in \cite[Remark 5.4, Remark 5.8, Lemma A.2]{FGFM-III}.
We verify Hypothesis \ref{hp:L2} for $N=1$ as the case $N>1$ is completely similar.
\vspace{-0.2truecm}
$$
\left\vert\overline{Pe^{t{A}}}\bar x\right\vert
=\vert\<e^{t\overline{{A}}}\bar x,\v\>\v\vert
\le \vert \v \vert_H \left\vert\<A^{-\eta}e^{t\overline{{A}}}\bar x, A^{\eta}\v\>\right\vert
\le \vert \v \vert_H \vert A^{\eta} \v \vert_H
\vert A^{-\eta}e^{t\overline{{A}}}\bar x \vert_H.
\vspace{-0.2truecm}
$$
Now, for $\alpha=3/4+\eps -\eta$,
\vspace{-0.2truecm}
$$
\vert A^{-\eta}e^{t\overline{{A}}}\bar x \vert_H
=
\vert A^{\alpha}e^{t\overline{{A}}}A^{-\eta-\alpha}\bar x \vert_H
\le Ct^{-\alpha}\vert A^{-\eta-\alpha}\bar x \vert_H=Ct^{-\alpha}\vert \bar x \vert_{\overline{H}}
\vspace{-0.2truecm}
$$
and the claim follows since $\eta>1/4+\eps$.
Finally we verify Hypothesis \ref{hp:smoothingextension}, which also imply the weaker Hypothesis  \ref{hp:smoothingextension-conv}
(see Remark \ref{lemma:ip:contr-sollevata-equiv-1}). Notice that if  $\<\v_1,...,\v_N\>$ are eigenvectors of ${A}$, then $P$ and $e^{t{A}}$ commute and so \ref{ip:NC} implies  \ref{hp:smoothingextension}, see remark \ref{rm:NCwhenPcommute}, but for general $\<\v_1,...,\v_N\>$ $P$ and $e^{t{A}}$ don't commute and so we have to prove Hypothesis \ref{hp:smoothingextension}.

\begin{proposition}\label{prop-ipsollevata-heat}
 Let us consider the uncontrolled version of the stochastic heat equation \eqref{eqDiri-abstr-contr}, let $R_{s}$, $0\leq s \leq T$ the transition semigroup related defined in \eqref{transition-heat},
and let ${B}$ the boundary control operator defined in \ref{notazioneB} and $P$ as in
\eqref{P-n-gen}. Then Hypotheses \ref{hp:smoothingextension}
is satisfied with $C=B_0$ and $\gamma=1-\delta$ for some $\delta \in (0,1/4)$.
 \end{proposition}
 \dim In this proof we consider the case of the projection on the space generated by only one element $\v\in D((-A)^\eta)$, namely $P:H\rightarrow H,\; Px=\<x,\v\>\v$, for all $x \in H,\, P=P^*$, the extension to a map as in \eqref{P-n-gen} being straightforward.
\newline In order to prove Hypothesis \ref{hp:smoothingextension}, point (i), we will prove \eqref{eq:NCdualhigh}. Indeed for any $t>0$ and for $z\in L_\rho^2([0,+\infty),H)$,
\begin{align}\label{conti1heat}
&\<\Upsilon^P_\infty   Q_{t} (\Upsilon^P_\infty )^*z,z\>_{L_\rho^2([0,+\infty),H)}
=\int_0^{t}\vert\int_0^{+\infty}e^{-\rho\tau}Q^{1/2}e^{(\tau+r)\overline{{A}}^*}P^*z(\tau) \,d\tau\vert^2\, dr\\ \nonumber
&\quad=\int_0^{t}\vert\int_0^{+\infty}e^{-\rho\tau}Q^{1/2}e^{(\tau+r)\overline{{A}}^*}\<z(\tau),\v\>\v \,d\tau\vert^2\, dr=\int_0^{t}\vert Q^{1/2}e^{r\overline{{A}}^*}\<\int_0^{+\infty}e^{-\rho\tau}e^{\tau\overline{{A}}^*}z(\tau)d\tau\, ,\v\>\v \vert^2\, dr\\ \nonumber
&\quad =\vert Q_{t}^{1/2}P \int_0^{+\infty}e^{-\rho\tau}e^{\tau\overline{{A}}^*}z(\tau)d\tau \vert^2\geq C t^{1-\delta} \vert {B}^* e^{(t-s)\overline{{A}}^*}P\int_0^{+\infty}e^{-\rho\tau}e^{\tau\overline{{A}}^*}z(\tau)d\tau\vert^2
\end{align}
where in the last passage we have applied \eqref{eq:NCdualbasic}, that $\eta >1/4+\eps$ and
\cite[Lemma A.2]{FGFM-III}. {Namely in \cite[Lemma A.2]{FGFM-III} it is proved that with $\eta> \frac14+\varepsilon$ and for $\delta\in(0,\frac14)$,
$ \operatorname{Im}\overline{Pe^{tA}}(-{A}D)\subset \operatorname{Im}\bar Q_t^{1/2}$,
and
$\Vert \bar Q_t^{-1/2}\overline{Pe^{tA}}(-{A}D)\Vert \leq C t^{-1+\delta}
$,
where $\bar Q_t:=PQ_tP.$
  }
In order to conclude it suffices to notice that, recalling also Proposition \ref{pr:propcalS}-(iv),\vspace{-0.2truecm}
$$\vert (B_0^*e^{(t-s)\overline{{A}^*}})(\Upsilon^P_s)^* z\vert_K= \vert {B}^* e^{(t)\overline{{A}}^*}P
\int_0^{+\infty}e^{\tau\overline{{A}}^*}z(\tau)
d\tau\vert^2.
\vspace{-0.6truecm}
$$
\qed
\vspace{-0.5truecm}

\vspace{-0.2truecm}

{
\begin{remark}
\label{rm:ApproxBoundary}
Hypothesis \ref{hp:BCcost}, (i) and (iii) undoubtedly restricts the set of {running} costs to which our result can be applied.
In this respect we make a couple of observations. On one hand, the dependence on a finite number of components with respect to vectors in $\cald((-A)^{\eta})$ arise from the coupled need of spatial regularity of the semigroup $R_t$ and integrability of $\hat\Lambda^{P,C}_s$. The coupling does not allow to choose e.g. $P=(-A)^{-\eta}$ as one could think at a first glance.
On the other hand, due to the results of \cite[Section B.6]{FabbriGozziSwiech}, in particular Lemma B.78, the set of such costs is dense, under the $\calk$-convergence introduced
in Definition \ref{k-conv}, into the set of continuous and bounded functions. This means that the {running} cost $\ell_0$ can be chosen close to any reasonable data in applied examples. Of course this does not guarantee that the corresponding solutions are close but in any case represents a first step towards to the understanding of such problems.
\end{remark}
}

\vspace{-0.6truecm}
\subsection{SDEs with delay in the control}\label{SSE:DELAYEQUATION}


We consider the following controlled stochastic
differential equation in $\R^n$ with delay in the control:
\vspace{-0.2truecm}
\begin{equation}
\left\{
\begin{array}
[c]{l}%
dy(s)  =a_0 y(s) {ds}+b_0 u(s) ds +\displaystyle\int_{-d}^0 u(s+\xi)b_1(d\xi)
+\sigma dW(s)
,\text{ \ \ \ }s\in [t,T] \\
y(t)  =y_0,\qquad
u(t+\xi)=u_0(\xi), \quad \xi \in [-d,0),
\end{array}
\right.  \label{eq-contr-rit}
\vspace{-0.2truecm}
\end{equation}
where we assume the following.

\begin{hypothesis}\label{hp:delaystate}
\begin{itemize}
\item[]
  \item[(i)] $W$ is a standard Brownian motion in $\R^k$, and $(\calf_t)_{t\geq 0}$ is the
augmented filtration generated by $W$;\vspace{-0.2truecm}
  \item[(ii)] the control strategy $u$ belongs to $\calu$ where
$\calu:=\left\lbrace u(\cdot):(\Omega\times [0,T]\to U):\;
\hbox{predictable} \right\rbrace$,
with $U$ a closed and bounded subset of $\R^m$;\vspace{-0.2truecm}
\item[(iii)] $a_0\in \call(\R^n;\R^n)$, $b_0 \in \call(\R^m;\R^n)$, $\sigma\in \call(\R^k;\R^n)$, $d>0$;\vspace{-0.2truecm}
  \item[(iv)] $b_1$
  is an $m\times n$ matrix of signed Radon measures on $[-d,0]$ (i.e. it is an element of the dual space of $C([-d,0],\call(\R^m;\R^n))$).
\end{itemize}
\end{hypothesis}
Given any initial datum $(y_0,u_0)\in {\R^n\times L^2([-d,0],\R^n)}$ and any admissible control $u\in \calu$ equation (\ref{eq-contr-rit}) admits a unique strong (in the probabilistic sense) solution which is continuous and predictable
(see e.g. \cite{IkedaWatanabe} Ch. 4, Sect. 2-3).


Now, using the approach of \cite{VK} (see \cite{GM} for the stochastic case, see also \cite{FGFM-III} where the case of unbounded control operator is considered), we reformulate equation (\ref{eq-contr-rit}) as an abstract stochastic differential equation in the Hilbert space $H=\R^n\times L^2([-d,0],\R^n)$.
The operator $A_1 : \cald(A_1) \subset H \rightarrow H$ is defined as follows: for $x=(x_0,x_1)\in H$,
\vspace{-0.2truecm}
\begin{equation}\label{A1}A_1x=( a_0 x_0 +x_1(0), -x_1'), \quad \cald(A_1)=\left\lbrace x\in H :x_1\in W^{1,2}([-d,0],\R^n), x_1(-d)=0 \right\rbrace.
\vspace{-0.2truecm}
\end{equation}
We denote by $A_1^*$ the adjoint operator of $A_1$:
\vspace{-0.2truecm}
\begin{equation}
 \label{Astar}
A_1^{*}x=( a_0 x_0, x_1'), \quad \cald(A_1^{*})=\left\lbrace x\in H:x_1\in W^{1,2}([-d,0],\R^n), x_1(0)=x_0 \right\rbrace .
\vspace{-0.2truecm}
\end{equation}
We denote by $e^{tA_1}$ the $C_0$-semigroup generated by $A_1$. For $x\in H$ we have
\vspace{-0.2truecm}
\begin{equation}
e^{tA_1} \left(\begin{array}{l}x_0 \\x_1\end{array}\right)=
\left(
\begin{array}
[c]{ll}%
e^{ta_0 }x_0+\int_{-d}^{0}1_{[-t,0]}{(s)} e^{(t+s)a_0 } x_1(s)ds \\[3mm]
x_1(\cdot-t)1_{[-d+t,0]}(\cdot).
\end{array}
\right)  \label{semigroup}
\vspace{-0.2truecm}
\end{equation}
Similarly, denoting by $e^{tA_1^*}=(e^{tA_1})^*$ the $C_0$-semigroup generated by $A_1^*$,
we have for
$z=\left(z_0,z_1\right)\in H$
\vspace{-0.2truecm}
\begin{equation}
e^{tA_1^*} \left(\begin{array}{l}z_0 \\z_1\end{array}\right)=
\left(
\begin{array}[c]{ll}
e^{t a_0^* }z_0 \\[3mm]
e^{(\cdot+t) a_0^* }z_0 1_{[-t,0]}(\cdot) +z_1(\cdot+t)1_{[-d,-t)}(\cdot).
\end{array}
\right)  \label{semigroupadjoint}
\vspace{-0.2truecm}
\end{equation}
The infinite dimensional noise operator is defined as
\vspace{-0.2truecm}
\begin{equation}
 \label{G}
G:\R^{k}\rightarrow H,\qquad Gy=(\sigma y, 0), \; y\in\R^k.
\vspace{-0.2truecm}
\end{equation}
The control operator $B_1$ is defined as
(here the control space is $K:=\R^m$ and we denote by $C'([-d,0],\R^n)$
the dual space of $C([-d,0],\R^n)$)
\vspace{-0.2truecm}
\begin{equation}
 \label{Bnotbdd}
\begin{array}{c}
B_1:\R^{m}\rightarrow \R^n \times C'([-d,0],\R^n),
   \\[2mm]
(B_1u)_0=b_0 u,
\quad \<f,(B_1u)_1\>_{C,C'} = \int_{-d}^0 f(\xi)b_1(d\xi )u , \quad u\in\R^m, \quad
f \in C([-d,0],\R^n).
\end{array}
\vspace{-0.2truecm}
\end{equation}
We will choose $\overline{H}:=\R^n \times C'([-d,0],\R^n)$ as $B_1$ takes values there.
The adjoint $B_1^*$ is
\vspace{-0.2truecm}
\begin{equation}
 \label{B*notbdd}
\begin{array}{c}
B_1^*:\R^n \times C{''}([-d,0],\R^n) \rightarrow \R^{m},
   \\[2mm]
B_1^*(x_0,x_1)=
b^*_0 x_0+\dis\int_{-d}^0 b_1^*(d\xi)x_1(\xi)d\xi,
\; (x_0,x_1)\in \R^n \times C([-d,0],\R^n),
\end{array}
\vspace{-0.2truecm}
\end{equation}
where we have denoted by $C{''}([-d,0],\R^n)$ the dual space of $C'([-d,0],\R^n)$, which contains $C([-d,0],\R^n)$. If $b_1$ is absolutely continuous with respect to the Lebesgue measure on $[-d,0]$ then $B_1$ is a bounded operator
from $\R^n$ to $H$, and we can easily write $e^{tA_1}B_1$: see \cite{FGFM-I}.
\\
If $b_1$ is as in Hypothesis \ref{hp:delaystate}-(iv), then $B_1$ is unbounded. Still it is possible to write $e^{tA_1}B_1$ by extending the semigroup, by extrapolation, to $\R^n\times C'([-d,0];\R^n)$. In particular the first component is written as

\vspace{-0.2truecm}
\begin{equation}\label{eq:etAB}
\left(e^{tA_1}B_1\right)_0:\R^m \to \R^n,\qquad    \left(e^{tA_1}B_1\right)_0 u=
    e^{ta_0}b_0u+ \int_{-d}^0 1_{[-t,0]}e^{(t+r)a_0}b_1(dr)u,\qquad u \in \R^m.
\vspace{-0.2truecm}
\end{equation}
Let us now define the predictable process
$Y=(Y_0,Y_1):\Omega \times [0,T]\to H$ as
\vspace{-0.2truecm}
$$
Y_0(s)=y(s), \qquad Y_1(s)(\xi)=\int_{-d}^\xi b_1(d\zeta)u(\zeta+s-\xi),
\vspace{-0.2truecm}
$$
where $y$ is the solution of \eqref{eq-contr-rit}
and $u\in \calu$ is the control process.
By \cite[Proposition 2]{GM}, the process $Y$
is the unique mild solution of the abstract evolution equation in $H$
\vspace{-0.2truecm}
\begin{equation}
\left\{
\begin{array}
[c]{l}
dY(s)  =A_1Y(s) ds+B_1u(s) ds+GdW(s)
,\text{ \ \ \ }t\in[ 0,T] \\
Y(0)  =x=(x_0,x_1),
\end{array}
\right.   \label{eq-astr}%
\vspace{-0.2truecm}
\end{equation}
where $x_1(\xi)=\int_{-d}^\xi b_1(d\zeta)u_0(\zeta-\xi)$, for $\xi\in [-d,0)$.
Note that we have $x_1\in L^2([-d,0];\R^n)$
The mild (or integral) form of (\ref{eq-astr}) is
\vspace{-0.2truecm}
\begin{equation}
Y(s)  =e^{(s-t)A_1}x+\int_t^se^{(s-r)A_1}B_1 u(r) dr +\int_t^se^{(s-r)A_1}GdW(r)
,\text{ \ \ \ }s\in[t,T]. \\
  \label{eq-astr-mild}%
\vspace{-0.2truecm}
\end{equation}


We now introduce a (finite horizon) cost that we have to minimize over all admissible control strategies in $\calu$: we are able to handle also the case of a running cost depending on the state, so generalizing the cost considered in \cite{FGFM-III}:
\vspace{-0.2truecm}
\begin{equation}\label{eq:costoconcretodelay}
\bar J(t,y_0,u_0;u(\cdot))=\E \left[\int_t^T \left[\bar\ell_0(s,y(s;t,x,u))+\ell_1(u(s))\right]\,ds + \bar\phi(y(T;t,x,u))\right]
\vspace{-0.2truecm}
\end{equation}
under the following assumption.
\begin{hypothesis}\label{hp:delaycost}
\begin{itemize}
 \item[(i)] $\bar\ell_0:[0,T]\times \R^n\rightarrow \R$, is continuous and bounded.
\vspace{-0.2truecm}
  \item[(ii)] $\ell_1:U\to \R$ is measurable  and bounded from below.
\vspace{-0.2truecm}
  \item[(iii)] $\bar\phi:\R^n\to \R$ is continuous and bounded.
\end{itemize}
\end{hypothesis}
Such cost functional, using the infinite dimensional reformulation given above, can be rewritten as
\vspace{-0.2truecm}
\begin{equation}\label{eq:costoastrattodelay}
J(t,x;u(\cdot))=\E \left[\int_t^T \left[\ell_0(s,Y(s;t,x,u))+\ell_1(u(s))\right]\,ds +
\phi(Y(T;t,x,u))\right]
\vspace{-0.2truecm}
\end{equation}
where $\ell_0:[0,T]\times H\rightarrow \R$ and $\phi:H\to \R$ are defined respectively as $\ell_0(t,(x_0,x_1))=\bar\ell_0(t,x_0)$ and $\phi(x_0,x_1)=\bar\phi(x_0)$ for all $x=(x_0,x_1)\in H$.

Now we set $A=A_1$, $C=B_1$ as from \eqref{Bnotbdd}, and $P(x_0,x_1)=(x_0,0)$, and we check that our Hypotheses are verified.

Hypotheses \ref{ip-sde-common} is clearly satisfied.

Concerning Hypothesis \ref{ip:PC}, Point (i) follows from our choice of $H$ and $\overline{H}$ above and from the extension of the semigroup generated by $A_1$ by extrapolation (see e.g. \cite[Ch. II.5]{EngelNagelBook}). Point (ii) is immediate from the definition of $B_1$ and the choice of $\overline{H}$.
For point (iii) we note that
$\operatorname{Im \, }P=\R^n \times \{0\}$
and that, by \myref{semigroup}, we have, for $x=(x_0,x_1) \in H$,
\vspace{-0.2truecm}
$$
Pe^{tA}x=\left(
e^{ta_0 }x_0+\int_{-d}^{0}1_{[-t,0]} e^{(t+s)a_0 } x_1(s)ds,
0 \right)
\vspace{-0.2truecm}
$$
Hence, also in this case, we can extend immediately $Pe^{tA}$ to
$\overline{Pe^{tA}}:\overline H \to H$
by setting, for $x=(x_0,x_1) \in \R^n \times C'([-d,0];\R^n)$,
$\overline{Pe^{tA}}x=\left(
e^{ta_0 }x_0+\int_{-d}^{0}1_{[-t,0]} e^{(t+\xi)a_0 } x_1(d\xi),
0 \right)$.
From this expression we also immediately deduce that Hypothesis \ref{hp:L2} holds.

The validity of Hypothesis \ref{ip:NC} in this case has been proved in
\cite[Lemma A.3]{FGFM-III}.
under the following additional assumption (which is always verified e.g. if $\sigma\sigma^T$ is invertible).

\begin{hypothesis}
\label{hp:new}
We have
\vspace{-0.2truecm}
\begin{equation}\label{ip-delay}
 \operatorname{Im}\left(e^{ta_0}b_0 +\int_{-d}^0 1_{[-t,0]}e^{(t+r)a_0}b_1(dr)
\right)
\subseteq\operatorname{Im}\sigma,
\quad \forall t>0.
\vspace{-0.2truecm}
\end{equation}
\end{hypothesis}

Now we go to the (more difficult) proof that Hypothesis \ref{hp:smoothingextension-conv} holds.
First of all we need an additional assumption on the measure $b_1$, namely besides Hypothesis \ref{hp:delaystate}, point (iv), we have to assume the following.
\begin{hypothesis}\label{ip-b1}
There exists $0<\varepsilon \leq d $ such that $b_1((-\varepsilon, 0])=0.$
\end{hypothesis}

As examples of measure satisfying hypothesis \ref{ip-b1} we mention the Dirac measure $\delta_{-d}$ and the linear combination of Dirac measures $\delta_{-d_1},..., \delta_{-d_n}$ with $0<d_1<...<d_n\leq d$.
\newline Notice that also assuming  Hypothesis \ref{ip-b1} we are not able to
prove Hypothesis \ref{hp:smoothingextension}, see Remark \ref{rm-delay-final} below.

\begin{remark}
Notice that Hypotheses \ref{hp:delaystate}-(iv) and \ref{ip-b1} on $b_1$ cover the very common case of pointwise delay, which is technically complicated to deal with: indeed it gives rise to an unbounded control operator $B$. This case has been treated also in \cite[section 3.2]{FGFM-III}, the difference is that in the present paper we are able to consider also running cost depending on the state, see \eqref{eq:costoconcretodelay}.
The case of pointwise delay is not the only one we can treat, in particular we recall that, in \cite{FGFM-I} and \cite{FGFM-II} the case of $b_1$ absolutely continuous with respect to the Lebesgue measure has been treated assuming
$
b_1(d\xi)=\bar b_1(\xi)d(\xi),\;
\bar b_1\in L^2([-d,0],\call(\R^m;\R^n)).
$
We are able to consider also such a case when  $\bar b_1(\xi)=0$ for $-\varepsilon\leq \xi\leq 0.$
\end{remark}

Before to go to the proof we need some preparations to compute the operators
$\Upsilon^P_s$ and the related ones.
First notice that, as it is also true for stochastic heat equations with boundary control treated in section \ref{SSE:HEATEQUATION} we have
$\Imm  \overline{Pe^{tA}}\subseteq \Imm P$.
Moreover for the present case $P$ can be immediately extended to $\overline{P}:\overline{H}\to H$ so $\overline{Pe^{tA}}=\overline{P}\,\,\overline{e^{tA}}$.
In particular we get that, for any $z\in H$,
\vspace{-0.2truecm}
$$
\< \overline {Pe^{tA}}y, z\>_H=\<e^{ta_0}y_0,z_0\>_{\R^n}+\int_{-d}^{0}1_{[-t,0]}(\xi) e^{(t+\xi)a_0 }z_0y_1(d\xi)=\<y,( \overline {Pe^{tA}})^* z\>_{\overline H, \overline{H}^\prime}.
\vspace{-0.2truecm}
$$
from which we deduce that
\vspace{-0.2truecm}
$$
(\overline { Pe^{t A}})^*:H\rightarrow \overline{H}, \; ( Pe^{t A})^*z=\left(\begin{array}{l}
e^{t a_0^*}z_0\\
1_{[-t,0]}(\cdot) e^{(t+\cdot)a_0^* }z_0
\end{array}\right)
\vspace{-0.2truecm}
$$
where the second component has to be considered as an element of $C^{\prime\prime}([-d,0], \mathbb R^n)$.
%
%
This implies that, for any $y\in \overline H$ and for any $z\in  L^2(0,s;H)$, (abbreviated with $L^2_s$)
\vspace{-0.2truecm}
\begin{align*}
\<\Upsilon^P_s y, z\>_{L^2_s}&=\int_0^s\<y,(\overline {Pe^{rA}})^* z(r)\>_{\overline H, \overline{H}^\prime}\,dr
=
\<y,\int_0^{s}( \overline {Pe^{r A}})^* z(r)\, dr\>_{\overline H, \overline{H}^\prime}\\
&=\<y_0, \int_{0}^{s}e^{r a_0^*}z_0(r)\, dr)\>_{\R^n}+\<y_1, \int_{0}^s1_{[-r,0]}(\cdot) e^{(r+\cdot)a_0^* }z_0(r)\,dr\>_{C^{\prime}([-d,0], \mathbb R^n),C^{\prime\prime}([-d,0], \mathbb R^n}.
\vspace{-0.2truecm}
\end{align*}
Notice that, by the dominated convergence theorem, the function
\vspace{-0.2truecm}
$$
\xi\mapsto  \int_{0}^{s}1_{[-r,0]}(\xi) e^{(r+\xi)a_0^* }z_0(r)(dr)=\int_{-\xi}^{s} e^{(r+\xi)a_0^* }z_0(r)dr
\vspace{-0.2truecm}
$$
is continuous and so the duality between $C^{\prime}([-d,0], \mathbb R^n)$ and $C^{\prime\prime}([-d,0], \mathbb R^n)$ reduces to the duality between $C^{\prime}([-d,0], \mathbb R^n)$ and $C([-d,0], \mathbb R^n)$.
Hence for any $z\in L^2_s$
\vspace{-0.2truecm}
\begin{equation}\label{eq:aggiuntonew}
(\Upsilon^P_s)^* z
=
\left(\int_{0}^{s}e^{r a_0^*}z_0(r)\, dr,
 \int_{0}^s1_{[-r,0]}(\cdot) e^{(r+\cdot)a_0^* }z_0(r)\,dr
\right).
\vspace{-0.2truecm}
\end{equation}

\begin{proposition}\label{prop-ipsollevata-delay}
 Let us consider
the abstract reformulation \eqref{eq-astr} of the stochastic equation with delay in the control \eqref{eq-contr-rit}, let Hypotheses \ref{hp:delaystate}, \ref{hp:new}, \ref{ip-b1} hold true and let $\sigma$ be invertible {so implying that Hypothesis \ref{hp:delaycost} is satisfied.}.
Then for any $t<\varepsilon$
(where $ \varepsilon$  is given in Hypothesis \ref{ip-b1}), Hypothesis \ref{hp:smoothingextension-conv} is satisfied with $\gamma = 1/2$.
 \end{proposition}
\dim
For $z\in L^2(0,s;H)$, using \eqref{eq:aggiuntonew}, \eqref{semigroupadjoint}, and \eqref{B*notbdd}, we get
\vspace{-0.2truecm}
\begin{align}\label{conti2}
&\left\vert  B_1^* \overline{e^{(t-s)A^*}(}\Upsilon^P_s)^* z\right\vert
=
\left\vert B_1^* \overline{e^{(t-s)A^*}}
\left(\int_{0}^{s}e^{\rho a_0^*}z_0(\rho)\, d\rho,
 \int_{0}^s1_{[-\rho,0]}(\cdot) e^{(\rho+\cdot)a_0^* }z_0(\rho)\,d\rho\right)
\right\vert\\[1mm]
\nonumber
&=\left\vert\int_0^s  B_1^* \left(
e^{(t-s+\rho)a_0^*}z_0(\rho)d\rho,
e^{(t-s+\rho+\cdot)
a_0^*}1_{[-(t-s+\rho),0]}(\cdot)z_0(\rho)
 d\rho\right)
\right\vert\\[1mm]
\nonumber
&=\left\vert\int_0^s b_0^* e^{(t-s+\rho)a_0^*} z_0(\rho)\, d\rho
 +\int_0^s \int_{-d}^0 b_1^*(d\xi) e^{(t-s+\rho+\xi)a_0^*}  z_0(\rho)1_{[-(t-s+\rho),0]}(\xi)\, d\rho\right\vert
\\[1mm] \nonumber
&\leq \left\vert  b_0^* e^{(t-s)a_0^*}\int_0^s e^{\rho a_0^*} z_0(\rho)\, d\rho\right\vert
 +\left\vert \int_{-d}^0 b_1^*(d\xi) \int_0^s e^{(t-s+\rho+\xi)a_0^*}  z_0(\rho)1_{[-(t-s+\rho),0]}(\xi)\, d\rho\right\vert
 =I+II
\vspace{-0.2truecm}
\end{align}
Concerning the first term we recall that, writing the operator $Q_t$ in a matrix form, we get, for all $t\ge 0$, $(x_0,x_1)\in H$
\vspace{-0.2truecm}
\begin{equation}\label{Q-Q^0}
Q_t\left(\begin{array}{l}
x_0\\
x_1
\end{array}\right)=\left(\begin{array}{ll}
Q^0_t&0\\
0&0
\end{array}\right)\left(\begin{array}{l}
x_0\\
x_1
\end{array}\right)= \left(\begin{array}{l}
Q^0_tx_0\\
\quad 0
\end{array}\right).
\vspace{-0.2truecm}
\end{equation}
As argued in \cite[Appendix A.2]{FGFM-III}
thanks to \eqref{ip-delay}, we have, for some $C>0$
\vspace{-0.2truecm}
$$
\operatorname{Im}e^{ta_0}b_0\subset \operatorname{Im}(Q^0_t)^{1/2} \text{ with }\Vert (Q^0_t)^{-1/2} e^{ta_0}b_0 \Vert\leq C t^{-1/2},\qquad \forall t>0,
\vspace{-0.2truecm}
$$
i.e., for all $k\in \R^n$,
$\vert b_0^*e^{t a_0^*}k\vert \leq C  t^{-1/2} \vert (Q^0_{t-s})^{1/2}k \vert$.
Taking into account \eqref{Q-Q^0} and letting $k=\displaystyle\int_0^s e^{\rho a_0^*} z_0(\rho)\, d\rho$
we get
\vspace{-0.2truecm}
$$
I=\left\vert b_0^*e^{(t-s)a_0^*}\int_0^s e^{\rho a_0^*} z_0(\rho)\, d\rho\right\vert
\leq C  (t-s)^{-1/2} \left\vert (Q^0_{t-s})^{1/2}\int_0^s e^{\rho a_0^*} z_0(\rho)\, d\rho \right\vert =C  (t-s)^{-1/2} \left\vert \left(\Upsilon^P_sQ_{t-s} (\Upsilon^P_s)^*\right)^{1/2}z\right\vert
\vspace{-0.2truecm}
$$
where in the last passage we have used the fact that
\vspace{-0.2truecm}
\begin{align*}
\<\Upsilon^P_sQ_{t-s}( \Upsilon^P_s)^*z,z\>_{L^2_\rho}=
\<Q_{t-s}( \Upsilon^P_s)^*z,\Upsilon^P_s)^*z\>_{\overline H}
=\<Q^0_{t-s}\int_0^s e^{\rho a_0^*} z_0(\rho)\, d\rho,\int_0^s e^{\rho a_0^*} z_0(\rho)\, d\rho\>_{\R^n}.
\vspace{-0.2truecm}
\end{align*}
It remains to estimate $II$:
\vspace{-0.2truecm}
\begin{align}\label{conti3}
II&=\left\vert \int_{-d}^0 b_1^*(d\xi) e^{(t-s+\xi)a_0^*}\int_0^s e^{\rho a_0^*}  z_0(\rho)1_{[-(t-s+\rho),0]}(\xi)\, d\rho\right\vert
\\[1mm] \nonumber
&=\left\vert \int_{-d}^0 b_1^*(d\xi) e^{(t-s+\xi)a_0^*}\int_{0\vee[-(t-s)-\xi]}^{s\vee [-(t-s)-\xi]} e^{\rho a_0^*}  z_0(\rho)1_{[-(t-s+\rho),0]}(\xi)\, d\rho\right\vert,\nonumber
\vspace{-0.2truecm}
\end{align}
and so it is clear that if $-(t-s)-\xi>s$, that is if $t<-\xi$, the integral is null, $II=0$. So if, according to Hypothesis \ref{ip-b1} , $b_1((-\varepsilon,0])=0$ and $t<\varepsilon$, we get $II=0$ and we can conclude that Hypothesis \ref{hp:smoothingextension-conv} holds true and this concludes the proof.
\qed

\vspace{-0.2truecm}

\begin{remark}
\label{rm-delay-final}
We notice that Hypothesis \ref{ip-b1} is related to the choice of working in the space $L^2(0,s;H)$, indeed in the proof of Proposition \ref{prop-ipsollevata-delay} the point is that the integral term $\displaystyle \int_0^s e^{\rho a_0^*}z_0(\rho)  \,d\rho$ cannot be estimated in terms of the norm of $z_0$ in $L^2(0,s;\R^n)$, and the integral $\displaystyle \int_0^s e^{\rho a_0^*}z_0(\rho)  \,d\rho$ appears in the computations of the adjoint operator of $\Upsilon^P_s$ considered as an operator taking values in $L^2(0,s;H)$. On the other side the choice of working in the Hilbert space $L^2(0,s;H)$ seems to be unavoidable in order to prove Propositions \ref{prop:partsmoothnew} and \ref{prop:partsmoothnew-conv} which allow us to deal with state dependent costs.
\newline We stress the fact that even if we assume Hypothesis \ref{ip-b1}, we are not able to verify Hypothesis \ref{hp:smoothingextension}. Indeed, to do this, we should consider, for $z\in L^2_\rho(0,+\infty;H)$
\vspace{-0.4truecm}
\begin{align*}
 &\vert B_1^* \overline{e^{(t)A^*}}(\Upsilon^P_\infty)^* z\vert=\vert B_1^* \overline{e^{(t)A^*}}\int_0^{+\infty} e^{-\rho\tau}\left(\overline{Pe^{\tau A}}\right)^* z(\tau)\, d\tau\vert\\ \nonumber
&=\vert\int_0^{+\infty} B_1^* e^{-\rho\tau}\overline{e^{(t\tau)A^*}}P^* z(\tau)\, d\tau\vert=\vert\int_0^{+\infty}  e^{-\rho\tau}B_1^* \left(\begin{array}{l}
e^{(t+\tau-)a_0^*}z_0(\rho)\\
e^{(t+\tau+\cdot )a_0^*}1_{[-(t+\tau),0]}(\cdot)z_0(\tau)
\end{array}\right)\,
 d\tau\vert\\ \nonumber
 &\leq\vert  b_0^* e^{ta_0^*}\int_0^{+\infty}  e^{-\rho\tau}e^{\tau a_0^*} z_0(\tau)\, d\tau\vert
 +\vert \int_{-d}^0 b_1^*(d\xi) \int_0^{+\infty}e^{-\rho\tau}e^{(t+\tau+\xi)a_0^*}  z_0(\tau)1_{[-(t+\tau),0]}(\xi)\, d\tau\vert
 =I+II
\vspace{-0.4truecm}
\end{align*}
Here, even letting Hypothesis \ref{ip-b1} hold true, the second term $II$ is not null if $\xi+t+\tau>0$, that is when the indicator function $1_{[-(t+\tau),0]}(\xi)=1$, unlike what we can achieve in \eqref{conti3} there is no restriction on $\tau$ that makes the integral identically $0$; consequently there is no way to bound $II$.
\end{remark}


{\begin{remark}
\label{rm:ApproxDelay}
\begin{itemize}
  \item The reader may find somehow restrictive the Hypothesis \ref{ip-b1} above. However, when dealing with applied models this is rarely so. On one hand all the (very common) cases where $b_1$ is a combination of Dirac deltas fits into our setting (recall that if there is a Dirac delta at $0$ this goes into the ``present'' part of the state). On the other side any absolutely continuous $b_1$ can be approximated, in the $L^2$ norm, by a sequence $b_{1,n}$ satisfying
      Hypothesis \ref{ip-b1}.
      Of course, as pointed out also in Remark       \ref{rm:ApproxBoundary} this does not guarantee that the corresponding solutions are close but in any case represents a first step towards to the understanding of such problems.
  \item An example of state equation where our setting applies is the one of the paper \cite{CarmonaEtAl18} (there the authors consider a game between banks, here we look only at the case of one agent for simplicity). In such paper the state equation is of the type
      $dy(s)= u(s) - u(s-d) + \sigma dW(s)$ for $s\in [t,T]$,
      which falls within our assumptions.
      Moreover in such paper the cost is quadratic: this allows the authors to find explicit solution of the HJB equation. Our technique can be used to find optimal feedback control in when the running cost and the final cost depend on the state in a more general way. Note that, using the first part of Remark \ref{rm:crescitapoli-HJB} also cases where such costs have polynomial growth could be treated.
\vspace{-0.3truecm}
\end{itemize}
\end{remark}
}


{
\section{Verification Theorem and Optimal Feedbacks}\label{sec-verifica}
}

\vspace{-0.2truecm}

{
The aim of this section is to provide a verification theorem (Theorem \ref{teorema controllo}) and the existence of optimal feedback controls (Theorem \ref{teo su controllo feedback}) for our problem.
These results will also imply that the mild solution $v$ of the HJB equation (\ref{HJBformale-common}), built in Theorem \ref{esistenzaHJB-progr}, is equal to the
value function $V$ of our optimal control problem.
The proof of such results proceeds along the lines of
\cite[Sections 4.5 and 4.8]{FabbriGozziSwiech} and \cite[Sections 4 and 5]{FGFM-II}: however such procedure  must be carefully adapted to this different context and the needed adaptation depends, in particular, on structural properties of the operators $A$ and $C$ (in particular that, for all $x \in \cald(A^*)$ we have $x \in \cald(C^*)$ and $\vert C^*x\vert \leq c \vert A^*x\vert $, see Remark \ref{rm:defpistrictstrong} for details).
Such property holds in the example of boundary control (Subsection \ref{SSE:HEATEQUATION}) while it does not in the example of delay in the control (Subsection \ref{SSE:DELAYEQUATION}). This means, in particular, that the proof in the first case is more similar to the one of \cite[Sections 4.5 and 4.8]{FabbriGozziSwiech}
while in the second case it requires a careful modification of the arguments used in
\cite[Sections 4 and 5]{FGFM-II}.
Given the above we decided to expose the details of the second (more involved) proof while we explain the ideas of the first one in the subsequent Remark \ref{rm:defpistrictstrong}\footnote{{We also explain, in Remark \ref{rm:delaysecondderivative} how this second proof could possibly be adapted to cover also the case of boundary control.}}. For simplicity we consider here the case when the image of $P$ is finite dimensional, which in any case holds in our examples.
This section is divided in three parts: the first about the approximation of the solution of HJB, the second about the verification theorem, the third about the feedback controls.}


{
\subsection{$\calk$-convergence and approximations of solutions of the HJB equation}
\label{sec-Kstrong}
}


{We first introduce the notion of $\calk$-convergence, following \cite[Definition 4.1]{FGFM-II}, and the references therein.}
{
\begin{definition}\label{k-conv} Let $\overline H$ be a real and separable Hilbert space.
A sequence $(f_n)_{n\geq 0}\in C_b(\overline H)$
is said to be $\calk$-convergent to a function $f\in C_b(\overline H)$ (and we shall write
$f_n\overset{\calk}{\rightarrow} f$ or $f=\calk-\lim_{n\rightarrow\infty}f_n$) if for any compact set
$\calk\subset\overline H$
\vspace{-0.3truecm}
$$
\sup_{n\in\N}\Vert f_n\Vert_{C_b(\overline H)}<+\infty \quad {\rm and } \quad
\lim_{n\rightarrow\infty}\sup_{x\in\calk}\vert f(x)-f_n(x)\vert =0.
\vspace{-0.3truecm}
$$
\noindent
Similarly, given $I\subseteq \R$, a sequence
$(f_n)_{n\geq 0}\in C_b(I \times \overline H))$ is said to be $\calk$-convergent to a function $f\in C_b(I\times \overline H))$ (and we shall write again
$f_n\overset{\calk}{\rightarrow} f$ or $f=\calk-\lim_{n\rightarrow\infty}f_n$) if for any compact subsets $I_0\subseteq I$ and
$\calk\subset \overline H$
\vspace{-0.3truecm}$$
\sup_{n\in\N}\Vert f_n\Vert_{C_b(I\times\overline H)}<+\infty \quad {\rm and } \quad
\lim_{n\rightarrow\infty}\sup_{(t,x)\in I_0\times\calk}\vert f(t,x)-f_n(t,x)\vert=0.
\vspace{-0.3truecm}
$$
\end{definition}
}
{Now we consider the following forward
Kolmogorov equation (with unknown $w$), which is obtained substituting, in
\eqref{HJBformaleforward-common}, the Hamiltonian $H_{min}$ with a generic known term $\calf$ and extending the operators to $\overline{H}$:
\vspace{-0.3truecm}
\begin{equation}\label{eq:Kforward}
  \left\{\begin{array}{l}\dis
\partial_t w(t,\bar x)=\overline{\call}[ w(t,\cdot)](\bar x)
+\calf(t,\bar x) +\ell_0(t,\bar x)
\qquad t\in (0,T],\,
\bar x\in \overline H,
\\
\dis w(0,x)=\phi(x), \qquad x \in H.
\end{array}\right.
\vspace{-0.3truecm}
\end{equation}
}
{Here $\phi\in C^P_b(H)$ and the functions $\calf, \ell_0$ both satisfy Hypothesis \ref{ip-costo}-(ii).
Moreover $\overline{\call}$ is defined as follows.
Let first $\overline{A}:\cald(\overline{A})\subseteq \overline{H}\to \overline{H}$ be the generator
of the extended semigroup $\overline{e^{tA}}$
as from Hypothesis \ref{ip:PC}-(i). Clearly $\overline{A}$
is an extension of $A$ and
$\cald({A})\subseteq\cald(\overline{A})$.
Then let, for $f \in C_b(\overline{H})$ sufficiently regular, and $\bar{x}\in \cald(\overline{A})$,
\vspace{-0.3truecm}
\begin{equation}\label{eq:overlineL}
\overline{\call}[f](\bar x):=
\frac{1}{2} Tr \; GG^*\nabla^2 f(\bar x)
+ \< \overline{A}\bar x,\nabla f(\bar x)\>
\vspace{-0.3truecm}
\end{equation}
Clearly $\overline{\call}[f]$ is found taking the extension of $\call [f]$ in \eqref{eq:ell-gen} (see also Remark \ref{rm:generatornew}) to $\overline{H}$.
\\
We now provide a definition of regular solution which is suitable for our purposes.
}

\vspace{-0.3truecm}

{
\begin{definition}\label{df:strictandpistrong}
Let $K$ be a Hilbert space and take an operator $C\in \call(K,\overline{H})$ as in Hypothesis \ref{ip:PC}-(ii).  By a $C$-\textit{strict solution} of the Kolmogorov equation (\ref{eq:Kforward}) we mean
a function $w$ such that
\vspace{-0.3truecm}
\begin{equation}\label{propr-strict-sol}
 \left\lbrace
\begin{array}{l}
w \in C_b([0,T]\times \overline H)\quad  {\rm and}\quad w(0, x)=\phi( x) \quad \forall  x \in {H};
\\
w(t,\cdot) \in UC^2_b(\overline{H}), \; \forall t\in [0,T]; \qquad
w(\cdot, \bar x)\in C^{1}([0,T]),\; \forall \bar x\in \overline{H};
\\
\nabla^C w\in C_b([0,T]\times \overline H, K),
\qquad
GG^*\nabla^2 w \in C_b([0,T]\times \overline H, \call_1(\overline{H})),
\\
w\text{ satisfies (\ref{eq:Kforward}) }
\forall t\in (0,T], \;\bar x\in D(\overline{A})
\end{array}
\right.
\vspace{-0.3truecm}
\end{equation}
\end{definition}
Concerning the last request in \eqref{propr-strict-sol}, recall that here $\call$ is given by \eqref{eq:ell-gen} and this expression make sense only for $x\in \cald(\overline{A})$.}


{
\begin{remark}
\label{rm:classicalvsstrict}
We can see the above definition as a modification
of the definition of \emph{strict solution} given in \cite[Definition 9.24]{DP1}. The core business here is the fact that here we give sense to the equation only for $\bar x\in \cald(\overline{A})$, while in the definition of \emph{classical solution} given in \cite[Definitions B.81 and B.82]{FabbriGozziSwiech} the equation makes sense for all $\bar x \in \overline{H}$ requiring more regularity on the gradient of the solution (namely that it belongs to $\cald(\overline{A}\;^*$).
\end{remark}
}


{We now prove a useful approximation result.}


{
\begin{lemma}\label{lm:approximation}
Let Hypotheses \ref{ip-sde-common}, \ref{ip:PC}, \ref{ip:NC}, \ref{hp:L2}, \ref{hp:smoothingextension-conv} and \ref{ip-costo} hold true. Let also the image of $P$ be finite dimensional.
Let $v$ be the mild solution of the HJB equation  (\ref{HJBformale-common}) according to Definition \ref{defsolmildHJB} and set $w(t,\bar x)=v(T-t,\bar x)$ for $(t,\bar x)\in (0,T]\times \overline H$.
Let moreover $\bar \phi$ be the map associated to the final datum $\phi$ along Definition \ref{df:spaziphi1}; similarly let $\hat\ell_0$ be the map associated to the running cost $\ell_0$ along Definition \ref{df:SP} (or $\bar\ell_0$ the map associated to $\ell_0$ along Definition \ref{df:spaziphi1} when the second of the Hypothesis \ref{ip-costo}-(ii) is verified).
\\
Then there exist three sequences of functions $(\bar\phi_n)$,
$(\hat \ell_{0,n})$ (alternatively
$\bar \ell_{0,n}$
if the second of the Hypothesis \ref{ip-costo}-(ii) is verified) and $(\calf_n)$ such that, for all $n\in \N$,
\vspace{-0.3truecm}    \begin{equation*}\label{eq:defapproxphif}
        \bar\phi_n \in C_c^\infty(\overline H),
        \qquad
        \hat \ell_{0,n} \in C_c^\infty((0,T]\times \calc^P_A((0,+\infty);\overline H) ), \qquad
        \calf_n  \in C_c^\infty((0,T]\times \calc^P_A((0,+\infty);\overline H) )
\vspace{-0.3truecm}
\end{equation*}
(alternatively $\bar \ell_{0,n}\in  C_c^\infty ((0,T]\times\overline H)$)
and
\vspace{-0.3truecm}
\begin{equation}\label{eq:convapproxphif}
        \bar\phi_n \rightarrow \bar\phi ,
        \qquad         \hat \ell_{0,n} \rightarrow \hat \ell_0 , (\text{alternatively  }\bar \ell_{0,n}\rightarrow \bar\ell_0)
        \qquad\calf_n \to H_{min} (\nabla^C w)
\vspace{-0.3truecm}
    \end{equation}
in the sense of $\calk$-convergence.
Moreover, defining: $\phi_n(x)=\bar\phi_n(Px)$ for $x \in H$:, $\ell_{0,n}(s,\bar x)=\hat\ell_{0,n}(s,y^P_{\bar x}(\cdot \wedge s))$ for $\bar x \in \overline{H}$ (alternatively $\ell_{0,n}(s,x)=\bar\ell_{0,n}(s,Px)$ for $x \in H$)
and
\vspace{-0.3truecm}
\begin{equation*}\label{v_n-mild}
w_n(t,\bar x):= R_t [\phi_n](\bar x) + \int_0^t R_{t-s}[\calf_n (s,\cdot)+\ell_{0,n}(s,\cdot)] (\bar x)ds
\qquad t \in(0,T],\; \bar x \in \overline{H},
\vspace{-0.3truecm}
\end{equation*}
the following hold:
\vspace{-0.3truecm}
\begin{itemize}
  \item $w_n \in UC^{1,2}_b ([0,T]\times \overline{H})  \cap \cals^{1,P}_{\gamma,prog}([0,T]\times \overline{H})$;
\vspace{-0.3truecm}
  \item $w_n$ is a strict solution of
  \myref{eq:Kforward} with $\phi_n$ in place of $\phi$,
$\calf_n$ in place of $\calf $ and $ \ell_{0,n}$ in place of $ \ell_{0}$, provided we can prove that $GG^*\nabla^2 w \in C_b([0,T]\times \overline{H}, \call_1(\overline{H}))$.
\vspace{-0.3truecm}
  \item we have, in the sense of $\calk$-convergence,
$w_n \rightarrow w$ and $t^{{\gamma}}\nabla^C w_n \to t^{{\gamma}}\nabla^C w$.
\end{itemize}
\end{lemma}
}

{
\dim
We divide the proof in three steps.}

{{\bf Step 1: choosing the three approximating sequences.}
Since the image of $P$ is finite dimensional we choose $\bar\phi_n$ and $\bar\ell_{0,n}$ to be
the standard approximations by convolution of $\bar \phi$ and $\bar \ell_0$\footnote{This can still be done if the image of $P$ is infinite dimensional choosing such sequences as from \cite[Lemma B.78]{FabbriGozziSwiech}.}.
The definition of $\hat\ell_{0,n}$ and of $\calf_n$ is similar and we present the procedure for the definition of $\calf_n$.  We first notice that, since $w \in \cals^{1,P}_{\gamma,prog}([0,T]\times \overline{H})$, then the function
$(0,T]\times \overline{H}\to \R$,
$(t,\bar x) \mapsto  \calf(t,\bar x):= H_{min} (\nabla^C w(t,\bar x))$
has the property that there exist $\hat f:(0,T]\times \calc^P_A ((0,+\infty);\R) \to \R$, such that the map $(t,\bar x)\mapsto t^\gamma\hat f(t,y^P_{\bar x}(\cdot\wedge t))$ is bounded and
$\calf (t,\bar x)=\hat f(t,y^P_{\bar x}(\cdot\wedge t))$.
We notice that for each fixed $t\geq 0$ the function $\hat f$ is defined in the Banach space $\calc^P_A ((0,t]$, which admits a Schauder basis, hence the approximation by convolution procedure introduced in \cite{PZ} can be extended also to Banach spaces, see \cite{MasBanach} where the same argument is used in the proof of Theorem 4.3, page 404. So let $\hat f_n$ be the approximation by convolution of $\hat f$ and define
$\calf_n (t,\bar x)=\hat f_n(t,y^P_{\bar x}(\cdot\wedge t))$.}

{
{\bf Step 2: proving that $w_n\in UC^{1,2}_b([0,T]\times H) \cap \cals^{1,P}_{\gamma,prog}([0,T]\times \bar H)$ and
that it is a strict solution.}
The fact that $w_n\in\cals^{1,P}_{\gamma,prog}([0,T]\times \overline H)$ follows immediately from \myref{v_n-mild} and Theorem \ref{esistenzaHJB-progr}.
\newline Differentiability of first and second order with respect to the variable $x$ follows by applying the dominated convergence theorem,
after explicitly differentiating
$R_{t}\left[\phi_n\right]$
and the convolution term containing $\calf_n$,
and $\ell_{0,n}$, see \cite{FGFM-II}, Lemma 4.3, \textit{Step 2}. Indeed the approximating coefficients $\phi_n,\,\hat{f_n} ,\,\ell_{0,n}$ turn out to be infinitely many times differentiable, and so we can differentiate  $ R_{t}\left[\phi_n\right]  $ and the convolution terms containing $ \hat{f_n} $ and $\ell_{0,n}$ differentiating $ \phi_n$, $ \hat{f_n} $ and $\ell_{0,n}$ respectively.
The proof that $w_n$ is differentiable with respect to time
is completely analogous to what is done in
\cite[Section 4.4 and B.7]{FabbriGozziSwiech},
(see also \cite[Theorems 9.23 and 9.25]{DP1}).
By Theorem 5.3 in \cite{CeGo}, see also Theorem 7.5.1 in \cite{DP3} for Kolmogorov equations,
we finally conclude that $w_n$
is a strict solution to equation \ref{eq:Kforward}.}

{{\bf Step 3: proof of the convergences.}
First we can  prove that the sequences $(w_n)$ and $(t^{\gamma}\nabla^C w_n)$
are bounded uniformly with respect to $n$ with computations similar to the ones in
\cite{FGFM-II}[Lemma 4.3].
\newline Now we prove the convergences, we restrict to the case when $\ell_0$ and $\ell_{0,n}$ are defined in terms of $\hat\ell_0$ and $\hat\ell_{0,n}$ respectively, being the case when they are defined in terms of $\bar\ell_0$ and $\bar\ell_{0,n}$ completely similar. We have
\vspace{-0.2truecm}
\begin{align*}
 \vert & w_n(t,\bar x)-w(t,\bar x)\vert 
 \leq \vert R_{t}[\phi-\phi_n](\bar x)\vert+\vert\int_0^t R_{t-s}[\calf(s,\cdot)-\calf_n(s,\cdot))](\bar x)\; ds\vert +\vert\int_0^t R_{t-s}[\ell_0(s,\cdot)-\ell_{0,n}(s,\cdot)](\bar x)\; ds\vert\\
& \leq \int_{\overline H}\vert \bar\phi\left(z_1+\overline{Pe^{tA}}\bar x\right)-\bar\phi_n\left(z_1+\overline{Pe^{tA}}\bar x\right)\caln(0,P Q_tP^*)(dz_1)\\
&\quad +\int_0^t \int_{H} \vert\hat f\left(\Upsilon^P_s(z+\overline{e^{(t-s)A}}\bar x \right)-\hat f_n\left(\Upsilon^P_s(z+\overline{e^{(t-s)A}}\bar x \right\vert)\caln(0,Q_{t-s})(dz)\\
 &\quad +\int_0^t \int_{H} \vert\hat\ell_0 \left(\Upsilon^P_s(z+\overline{e^{(t-s)A}}\bar x \right)-\hat \ell_{0,n}\left(\Upsilon^P_s(z+\overline{e^{(t-s)A}}\bar x \right\vert)\caln(0,Q_{t-s})(dz)\\
 & = \int_{\overline H}\vert \bar\phi\left(z_1+\overline{Pe^{tA}}\bar x\right)-\bar\phi_n\left(z_1+\overline{Pe^{tA}}\bar x\right)\caln(0,P Q_tP^*)(dz_1)\\
&\quad +\int_0^t \int_{L^2_s} \vert\hat f\left(z_1+\Upsilon^P_s \overline{e^{(t-s)A}}\bar x \right)-\hat f_n\left(z_1+\Upsilon^P_s\overline{e^{(t-s)A}}\bar x \right)\vert\caln(0,\Upsilon^P_sQ_{t-s}(\Upsilon^P_s)^*)(dz_1)\\
 &\quad +\int_0^t \int_{H} \vert\hat\ell_0 \left(z_1+(\Upsilon^P_s\overline{e^{(t-s)A}}\bar x \right)-\hat \ell_{0,n}\left(z_1+\Upsilon^P_s\overline{e^{(t-s)A}}\bar x \right\vert)\caln(0,\Upsilon^P_sQ_{t-s}(\Upsilon^P_s)^*)(dz_1)\\
\vspace{-0.4truecm}
\end{align*}
Since, for every compact set $\calk \subset \overline H$ the set
$\{\overline{Pe^{tA}}\bar x, \; t \in [0,T],\, x \in \calk\}$ is compact in $H$, hence in $\overline{H}$, and the set $\{\Upsilon^P_s \overline{e^{(t-s)A}}\bar x, \;0\leq s\leq t,\, t \in [0,T],\, x \in \calk\}$ is compact in $\calc^P_A((0,t];H)$ since it is the image of a compact set through the continuous ( see Lemma
\ref{lm:UpsilonPT} ) map $\Upsilon^P_s$.
Then by the Dominated Convergence Theorem we get that for any compact set $\calk \subset \overline H$
\vspace{-0.2truecm}
\begin{equation*}
\label{eq:convvnnew}
\sup_{(t,x)\in[0,T]\times\calk} \vert w_n(t,x)-w(t,x)\vert\rightarrow 0.
\vspace{-0.2truecm}
\end{equation*}
We can prove a similar result for the convergence of $\nabla^C w_n(t,\bar x)$ to $\nabla^C w(t,\bar x)$ on the lines of the proof of Lemma 4.3, \textit{Step 3}, in \cite{FGFM-II}, and applying
the last part of Proposition \ref{prop:partsmoothnew} for the term involving $\bar \phi,\, \bar\phi_n$, and Proposition \ref{prop:partsmoothnew-conv} for the other two terms. We finally get
\[
\sup_{\bar x\in(0,T]\times\calk} \vert t^{\gamma}\nabla^C w_n(t,\bar x)-t^{\gamma}\nabla^C w(t,\bar x)\vert\rightarrow 0, \qquad  \hbox{for any compact set $\calk\subset\overline H$.}
\]
\qed}
\medskip

\vspace{-0.5truecm}

{\begin{remark}\label{rm:defpistrictstrong}
The approximation result proved just above is needed to prove the fundamental identity,
(see next Proposition \ref{prop rel fond}) which is the key step to get the verification theorem and the existence of optimal feedback controls.
On this we think it is useful to observe the following.
\\
Along the lines of \cite[Remark 4.4]{FGFM-II},
the previous approximating procedure does not give rise in general to
functions $w_n$ with $\nabla w_n\in D(A^*)$, so we only find approximating strict solutions
and not classical solutions
(as from Remark \ref{rm:classicalvsstrict}).
\\
The reason for our choice follows from the fact that, for our purposes, we need that the approximants of the data $\phi$, $\ell_0$, $\calf$ remain in the space where we have ``partial'' smoothing. If we would use the approximation through classical solutions (as in \cite[Section B.7.3]{FabbriGozziSwiech}), the fact that $\nabla w_n\in D(A^*)$ is true when there exists $c_0>0$ such that
\vspace{-0.2truecm}
\begin{equation}\label{eq:DCDAstar}
\hbox{for every $\bar x \in \cald(\overline{A}\,^*)$ we have $\vert C^*\bar x\vert \leq c_0 \vert \overline{A}\,^*x\vert$.}
\vspace{-0.2truecm}
\end{equation}
 as pointed out in \cite[Example 4.137 and Theorem 4.148]{FabbriGozziSwiech}.
Relation \eqref{eq:DCDAstar} is indeed true in our example with boundary control (Subsection \ref{SSE:HEATEQUATION});
however this is not true in our example with delay in the control (Subsection \ref{SSE:DELAYEQUATION}).
The consequence is that, in the case of boundary control
we could prove the results adapting, with some modifications the approximation through classical solution adopted in \cite[Section B.7.3]{FabbriGozziSwiech}.
On the other hand, to cover the example of delay in the control we need to change the approach working with strict solution on the line of what is done in \cite{FGFM-II}.
Note finally that, since the state process $X$ does not belong to $D(\overline{A})$ in general, this fact will force us to introduce suitable regularizations $X_k$ of $X$ (see the proof of Proposition \ref{prop rel fond}).
\end{remark}
}

\vspace{-0.5truecm}

{\begin{remark}\label{rm:delaynotgood}
We now explain why the example with delay in the control
of Subsection \ref{SSE:DELAYEQUATION} does not satisfy the property \eqref{eq:DCDAstar}. Recall that $A=A_1$ $C=B_1$ in this case. We take any element
$(x_0,x_1)\in \{\R^n\times C^1([-d,0],\R^n)\;|\; x_0=x_1(0)\}\subseteq \cald(A_1^*)\subseteq \cald(\overline{A}_1^{\;*})$.
For such elements we have
\vspace{-0.3truecm}
\begin{equation*}
\vert B_1^*(x_0,x_1) \vert_{\R^m} =
\left\vert b^*_0 x_0+\dis\int_{-d}^0 b_1^*(d\xi)x_1(\xi)d\xi\right\vert_{\R^m},
\vspace{-0.3truecm}
\end{equation*}
see \eqref{B*notbdd}, and
$\left\vert \overline{A}_1^{\;*}(x_0,x_1)\right\vert^2_{\overline{H}}
=
\left\vert A_1^*(x_0,x_1) \right\vert_{\overline{H}}^2
=
\vert a_0x_0\vert_{\R^n}^2+ \sup_{\theta \in [-d,0]}\vert x_1^\prime(\theta)\vert^2_{\R^n}$.
see \eqref{Astar}. If we take, for example, $m=n=1$,
$a_0=0$, $b_0=1$, $b_1=\delta_{-d}$ and the sequence
$x^k=(k,k+\xi)$, we have
$|B_1^* x^k|=2k-d$ and
$\vert A_1^* x^k\vert_{\overline{H}}=1$. Hence
property \eqref{eq:DCDAstar} does not hold in this case.
\end{remark}
}

\vspace{-0.5truecm}

{\begin{remark}\label{rm:delaysecondderivative}
We note that the machinery of Lemma \ref{lm:approximation} works if we can prove that $GG^*\nabla^2 w \in C_b([0,T]\times \overline H, \call_1(H))$. This is immediate for the case of delay in the control presented in Subsection \ref{SSE:DELAYEQUATION} since $G$ has finite dimensional image.
\\
On the other hand in the boundary control case of Subsection \ref{SSE:HEATEQUATION} this is not obvious.
A way to prove it in such case could be to take $\bar \phi_n,\,\hat\calf_n$ and $\hat\ell_{0,n} $ (or $\bar\ell_{0,n}$)
using the cylindrical functions, introduced in \cite{FabbriGozziSwiech}[Lemma B.78]
with respect to a suitable orthonormal basis contained in $\cald(A^*)$. This may not be obvious as one has to extend the arguments of \cite{FabbriGozziSwiech}[Lemma B.78]
to the Banach space setting in a more general way than what is done in \cite{MasBanach}.
\end{remark}
}

\vspace{-0.5truecm}

{\subsection{The Fundamental Identity and the Verification Theorem}}

\vspace{-0.2truecm}

{Now we go back to the control problem of Section \ref{sec-setting}.
We rewrite here for the reader convenience the state equation (\ref{eq-common-contr}), with generic initial datum in $\overline{H}$:
\vspace{-0.3truecm}
\begin{equation*}
  \begin{array}{l}
  \dis
d X(s)= AX(s)\,ds+Cu(s)\,ds +GdW(s), \qquad s\in [t,T],
\qquad X(t)=\bar x\in \overline{H}.
\end{array}
\vspace{-0.3truecm}
\end{equation*}
and the objective functional in \eqref{costoastratto-common}
\vspace{-0.2truecm}
\begin{equation*}
J(t,\bar x;u)=\E \left(\int_t^T \left[\ell_0(s,X^{t,x}(s))+\ell_1(u(s))\right]\,ds + \phi(X^{t,x}(T))\right),\qquad t\in [0,T), \; \bar x \in \overline{H}.
\vspace{-0.2truecm}
\end{equation*}
We underline that throughout this subsection and the following one, in order to avoid further technical difficulties, we keep the probability space $(\Omega,\calf,\P)$
fixed. Nothing would change if we work in the weak formulation, where the probability space can change (see e.g. \cite{FabbriGozziSwiech}[Chapter 2] for more on strong and weak formulations
in infinite dimension).
\\
We now prove the fundamental identity.}
{
\begin{proposition}\label{prop rel fond}
Let Hypotheses \ref{ip-sde-common}, \ref{ip:PC}, \ref{ip:NC}, \ref{hp:L2}, \ref{hp:smoothingextension-conv} and \ref{ip-costo} hold true. Let also the image of $P$ be finite dimensional.
Let $v$ be the mild solution of the HJB equation (\ref{HJBformale-common})
according to Definition \ref{defsolmildHJB}.
Then for every $t\in[0,T) $ and $\bar x\in \overline H$, and for every admissible control $u$, we have the fundamental identity
\vspace{-0.2truecm}
\begin{equation}\label{relfond}
 v(t,\bar x)
=J(t,\bar x;u)+\E\int_t^T \left[H_{min}(\nabla^C v(s,X(s)))
- H_{CV}(\nabla^C v(s,X(s));u(s))\right]\,ds.
\vspace{-0.2truecm}
 \end{equation}
\end{proposition}
}
{
\dim The proof can be performed following the lines of the proof of Proposition 5.1 in \cite{FGFM-II}.
During this proof, to avoid heavy notation we will write
$A$ for $\overline{A}$.
\\
Take any admissible state-control couple $(X(\cdot),u(\cdot))$,
and let, for $t\in(0,T]$, $\bar x \in \overline{H}$, $v_n(t,\bar x):=w_n(T-t,\bar x)$ where $(w_n)_n$ is the approximating sequence
of strict solutions defined
in Lemma \ref{lm:approximation}.  In order to apply the Ito formula to
$v_n(t, X(t))$, we need to regularize $X(t)$ since it does not live in $D(A)$; to this aim we define for $k \in \N$, sufficiently large,
\vspace{-0.2truecm}
\begin{equation*}\label{eq:Xkdef}
X_k(s)=k(k-A)^{-1}X(s).
\vspace{-0.2truecm}
\end{equation*}
The process $X_k$ is in $D(A)$, it converges to $X$ ($\P$-a.s. and $s\in [t,T]$ a.e.)
and it is a strong solution\footnote{{Here we mean strong in the probabilistic sense
and also in the sense of \cite{DP1}, Section 5.6.}}
of the Cauchy problem
\vspace{-0.2truecm}
\begin{equation*}
\left\{
\begin{array}
[c]{l}
dX_k(s)  =AX_k(s) ds+C_ku(s) ds+G_kdW_s
,\text{ \ \ \ }s\in [t,T] \\
X_k(t)=\bar x_k,
\end{array}
\right.
\vspace{-0.2truecm}
\end{equation*}
where $C_k=k(k-A)^{-1}C$, $G_k=k(k-A)^{-1}G$ and $\bar x_k=k(k-A)^{-1}\bar x$.
Now observe that the operator $C_k$ is continuous in $\overline H$, hence we can apply Dynkin's formula (see e.g.
\cite[Section 1.7]{FabbriGozziSwiech} or \cite[Section 4.5]{DP1})
to $v_{n}(s, X_k(s))$ in the interval $[t,T]$.}
{Using the Kolmogorov equation (\ref{eq:Kforward}), whose strict solution is $w_n$,
we then write:
\vspace{-0.2truecm}
\begin{equation}\label{quasirelfondv^nk}
\E\phi_n(X_k(T)) - v_n(t,\bar x_k)=\E\int_t^T \left[\calf_n (s,X_k(s))+ \ell_{0,n}(s,X_k(s))+\<C_k u(s),\nabla v_n(s,X_k(s))\>_{\overline{H}} \right]ds
\vspace{-0.2truecm}
\end{equation}
Now we pass to the limit for $k\rightarrow \infty$, on the lines of what is done in the proof of the fundamental identity in \cite{FGFM-II}.
Applying the Dominated Convergence Theorem to all terms but the last we get
\vspace{-0.2truecm}
\begin{align}\label{quasirelfondv^n}
&\E\phi_n(X(T)) - v_n(t,\bar x)\\ \nonumber&=\E\int_t^T [\calf_n (s,X(s))+\ell_{0,n}(s,X(s))]ds
+\lim_{k \to + \infty}\E\int_t^T\<C_k u(s),\nabla v_n(s,X_k(s))\>_{\overline{H}}]ds.
\vspace{-0.2truecm}
\end{align}
}
{Concerning the last term we argue as in the proof of fundamental identity in \cite{FGFM-II} to get
\vspace{-0.2truecm}
\begin{equation}
  \label{eq:limnablaC}
\lim_{k \to + \infty}\E\int_t^T\<C_k u(s),\nabla v_n(s,X_k(s))\>_{\overline{H}} ds
=
\E\int_t^T\<u(s),\nabla^C v_n(s,X(s))\>_{K}ds
\vspace{-0.2truecm}
\end{equation}}
{Now we let $n\rightarrow\infty$. By Lemma \ref{lm:approximation}, we know that
$$v_n(t,\bar x)\rightarrow v(t,\bar x) \quad\text{  and  }\quad(T-t)^{\gamma} \nabla^Bv_n(t,\bar x)\rightarrow (T-t)^{\gamma}\nabla^C v(t,x)
$$
in the sense of $\calk$-convergence. This, together with the properties of the approximating sequences found in Lemma \ref{lm:approximation} allows to apply Dominated Convergence into \eqref{quasirelfondv^nk} and \eqref{eq:limnablaC} getting}
{
\vspace{-0.2truecm}
$$
E\phi(X(T)) - v(t,\bar x)=
\E\int_t^T [H_{min}(\nabla^C v(s,X(s)))+ \ell_0(s,X(s))]ds
+\E\int_t^T \<u(s),\nabla^Cv(s,X(s))\>_{K}\,ds
\vspace{-0.2truecm}
$$
}
{Now, adding and subtracting $\E\dis\int_t^T  \ell_1(u(s))ds$ and rearranging the terms, we get
\vspace{-0.2truecm}
\begin{align*}
 v(t,\bar x)&
=\E\phi(X(T))+\E\int_t^T [\ell_0(s,X(s))+ \ell_1(u(s))]ds\\
&+\E\int_t^T \left[H_{min}(\nabla^C v(s,X(s)))-
H_{CV}(\nabla^C v(s,X(s));u(s))\right]
\,ds
\vspace{-0.2truecm}
\end{align*}
which immediately gives the claim.
\qed}




{We can now state our Verification Theorem i.e.
a sufficient condition of optimality given in term of the
mild solution $v$ of the HJB equation (\ref{HJBformale-common}). We omit the proof as it follows immediately form the previous Proposition \ref{prop rel fond}.
\begin{theorem}
\label{teorema controllo}
Let Hypotheses \ref{ip-sde-common}, \ref{ip:PC}, \ref{ip:NC}, \ref{hp:L2}, \ref{hp:smoothingextension-conv} and \ref{ip-costo} hold true. Let also the image of $P$ be finite dimensional.
Let $v$ be the mild solution of the HJB equation (\ref{HJBformale-common})
according to Definition \ref{defsolmildHJB}.
Then the following holds.
\vspace{-0.2truecm}
\begin{itemize}
 \item For all $(t,\bar x)\in [0,T)\times \overline H$ we have
$v(t,\bar x) \le V(t,\bar x)$, where $V$ is the value function
defined in \eqref{valuefunction-gen}.
\vspace{-0.2truecm}
\item
Let $t\in [0,T)$ and $\bar x\in \overline H$ be fixed.
If, for an admissible control $u^*$, we
have, calling $X^*$ the corresponding state,
\vspace{-0.2truecm}
$$
u^*(s)\in \arg\min_{u\in U}H_{CV}(\nabla^C v(s,X^*(s);u)
\vspace{-0.2truecm}
$$
for a.e. $s\in [t,T]$, $\P$-a.s., then the pair $(u^*,X^*)$ is
optimal for the control problem starting from $\bar x$ at time $t$
and $v(t,\bar x)=V(t,\bar x)=J(t,\bar x;u^*)$.
\end{itemize}
\end{theorem}}

\subsection{{Optimal feedback controls and $v=V$}}
\label{sec:contr-feedback}
{We now prove the existence of optimal feedback controls.}
{Under the Hypotheses
of Theorem \ref{teorema controllo} we define, for $(s,\bar x)\in [0,T)\times \overline H$,
the {\em feedback map}
\vspace{-0.2truecm}
\begin{equation}\label{defdiPsi}
\Psi(s,\bar x):=\arg \min_{u\in U} H_{CV}(\nabla^Cv(s,\bar x);u),
\vspace{-0.2truecm}
\end{equation}
where, as usual, $v$ is the solution of the HJB equation \myref{HJBformale-common}.
Given any $(t,\bar x)\in [0,T)\times \overline H$,
the so-called Closed Loop Equation
is written, formally, as
\vspace{-0.2truecm}
\begin{equation}\label{cleinclusion}
\left\{
\begin{array}
[c]{l}
dX(s) \in AX(s) ds+C\Psi\left(s,X(s)\right) ds+GdW_s
,\text{ \ \ \ }s\in [t,T) \\
Y(t)=\bar x.
\end{array}
\right.
\vspace{-0.2truecm}
\end{equation}
First of all we have the following straightforward corollary whose proof is immediate from Theorem \ref{teorema controllo}.}
{
\begin{corollary}
\label{cr:optimalfeedback}
Let the assumptions of Theorem \ref{teorema controllo} hold true.
Let $v$ be the mild solution of \myref{HJBformale-common}.
Fix $(t,\bar x)\in [0,T)\times \overline H$ and assume that, on $[t,T)\times \overline{H}$, the map $\Psi$
defined in (\ref{defdiPsi}) admits a measurable selection $\psi$
such that the Closed Loop Equation
\vspace{-0.2truecm}
\begin{equation}
\label{eq:CLEselection}
\left \{
\begin{array}{l}
d X(s) = AX(s)d s+C\psi\left(s,X(s)\right) ds+GdW_s
,\text{ \ \ \ }s\in [t,T) \\
X(t)=x.
\end{array}
\right.
\vspace{-0.2truecm}
\end{equation}
has a mild solution $X_\psi(\cdot;t,\bar x)$  (in the sense of \cite[p.187]{DP1}).
Define, for $s \in [t,T)$, $u_\psi (s)=\psi(s,X_\psi(s;t,\bar x))$.
Then the couple
$(u_\psi(\cdot),X_\psi(\cdot;t,\bar x))$ is optimal at
$(t,\bar x)$ and $v(t,\bar x)=V(t,\bar x)$.
If, finally, $\Psi(t,\bar x)$ is always a singleton and the mild solution
of \myref{eq:CLEselection} is unique,
then the optimal control is unique.
\end{corollary}
}
{
We now give sufficient conditions to verify the assumptions of Corollary
\ref{cr:optimalfeedback}.
First of all define
\vspace{-0.2truecm}
 \begin{equation}\label{defdigammagrandebis}
\Gamma(p):=\left\{ u\in U: \<p,u\>+\ell_1(u)= H_{min }(p)\right\}.
\vspace{-0.2truecm}
\end{equation}
Observe that, under mild additional conditions on $U$ and $\ell_1$
(for example taking $U$ compact or $\ell_1$ of superlinear growth),
the set $\Gamma$ is nonempty for all $p \in K$.
If this is the case then, by \cite{AubFr}, Theorems 8.2.10 and
8.2.11, $\Gamma$ admits a measurable selection, i.e. there exists
a measurable function $\gamma:K \rightarrow U$ with
$\gamma(z)\in \Gamma(z)$ for every $z\in K$.
Since we assumed that $H_{min }$ is Lipschitz continuous, then $\Gamma$, and so
$\gamma$, must be uniformly bounded.
In some cases
studied in the literature this is enough to find an optimal feedback
but not in our case due to the lack of the so-called ``structure condition''
(see on this \cite[Remark 5.6]{FGFM-II}).
Hence to prove existence of a mild solution
of the closed loop
equation \myref{eq:CLEselection}, as requested in Corollary \ref{cr:optimalfeedback},
we need more regularity of the feedback
term $\psi(s,\bar x)=\gamma(\nabla^C v(s,\bar x))$.
\\
Consequently, beyond the smooth assumptions on the coefficients required in Theorem \ref{teo su controllo feedback} below, which give the regularity of $\nabla^C v(t,x)$, we also need the following assumption about the map $\Gamma$.}

{
\begin{hypothesis}\label{hp:lipsel}
The set-valued map $\Gamma$ defined in (\ref{defdigammagrandebis})
is always non empty; moreover it
admits a Lipschitz continuous selection $\gamma$.
\end{hypothesis}}

%

\vspace{-0.2truecm}

{Taking the selection $\gamma$ from Hypothesis \ref{hp:lipsel}
we consider the closed loop equation
\vspace{-0.2truecm}
\begin{equation}\label{cle}
\left\{
\begin{array}
[c]{l}
dX(s)  =AX(s) ds+C\gamma\left(\nabla^C v(s,X(s))\right) ds+GdW_s
,\text{ \ \ \ }s\in[t,T] \\
X(t)=\bar x
\end{array}
\right.
\end{equation}
}
{
and we have the following result.}

\vspace{-0.2truecm}

{
\begin{theorem}\label{teo su controllo feedback}
Let Hypotheses \ref{ip-sde-common}, \ref{ip:PC}, \ref{ip:NC}, \ref{hp:L2}, \ref{hp:smoothingextension-conv}, \ref{ip-costo} and \ref{hp:lipsel} hold true. Let also the image of $P$ be finite dimensional.
Let $v$ be the mild solution of the HJB equation (\ref{HJBformale-common}).
Fix any $(t,\bar x)\in [0,T)\times\overline H$.
Assume also that, $\forall\,s\in[0,T]$,
$\bar \phi : H\rightarrow \R$ and $\hat\ell_0(s, \cdot):\calc_p^A((0,s];\overline{H})\rightarrow \R$  are Lipschitz continuous.
Then the closed loop equation (\ref{cle}) admits a unique mild solution $X_\gamma(\cdot;t,x)$ (in the sense of \cite[p.187]{DP1}); setting, for $s\in [t,T]$,
$u_\gamma(s):=\gamma\left(\nabla^{C} v(s,X_\gamma(s;t,x) )\right)$,
we obtain an optimal control at $(t,x)$ which is unique if $\Gamma$ is always a singleton. Moreover $v(t,x)=V(t,x)$.
\end{theorem}
}
{
\noindent {\bf Proof.}
In order to prove the existence and uniqueness of the mild solution of \eqref{cle} we first need to prove that, since $\bar\phi$ and $\hat \ell_0$ are Lipschitz continuous in $x$, uniformly in $t$, then also $\nabla^Cv$ is.
To prove this one has to modify to proof of Theorem \ref{esistenzaHJB-progr} making the fixed point in the space of functions $f$ such that $f(T-\cdot)\in \cals^{1,P}_{\gamma,prog}([0,T]\times \overline{H})$
and $(T-t)^\gamma\nabla^C f(t,\cdot)$ is Lipschitz continuous, uniformly in $t\in[0,T]$.
Such modifications can be done through straightforward modifications and we omit it.}
{
Given the above we can apply a fixed point theorem to the following integral form of \eqref{cle}:
\vspace{-0.2truecm}
\begin{equation*}\label{clemild}
X(s)= e^{(s-t)A}x+
\int_t^s e^{(s-r)A}G\,dW_r
+\int_{t}^s e^{(s-r)A}C\gamma(\nabla^{C}v(r ,\overline{X}_r))dr.
\vspace{-0.2truecm}
\end{equation*} and we get that it admits a unique solution, exactly as in \cite[Theorem 5.7]{FGFM-II}.
\qed
}

\vspace{-0.4truecm}

{
\begin{remark}\label{rm:verifica-con}
We notice that the identification $v=V$ can be done,
using an approximation procedure, also
in cases when we do not know if optimal feedback controls exist, so also when we do not assume further regularity (Lipschitz continuity in Theorem \ref{teo su controllo feedback}) on the coefficients. The proof goes along the lines of what is done in \cite{FGFM-II}[Theorem 5.9].
\end{remark}
}

\vspace{-0.4truecm}

{
\begin{remark}\label{rm:paper3}
We underline that the verification Theorem \ref{teorema controllo} and the existence of optimal controls proved in Theorem \ref{teo su controllo feedback} are valid also for the optimal control problem treated in \cite{FGFM-III}, where the running cost does not depend on the state.
\end{remark}
}

\appendix

\section{Notation and C-Derivatives}\label{section-prel}

\subsection{Basic notation and spaces}\label{subsection-notation}
For the reader's convenience we collect here the basic notation used throughout the paper.
\\
Let $H$ be a Hilbert space.
The norm of an element $x$ in $H$ will be denoted by
$\left|  x\right|_{H}$ or simply $\left|x\right|$,
if no confusion is possible, and by
$\left\langle \cdot,\cdot\right\rangle _H$,
or simply by $\left\langle \cdot,\cdot\right\rangle$ we denote the inner product in $H$.
We denote by $H^{\ast}$ the dual space of $H$.
If $K$ is another Hilbert space, $\call(H,K)$ denotes the
space of bounded linear operators from $H$ to $K$ endowed with the usual operator norm.
All Hilbert spaces are assumed to be real and separable.

Let $E$ be a Banach space. As for the Hilbert space case, the norm of an element $x$ in $E$ will be
denoted by $\left|x\right|_{E}$ or simply $\left|x\right|$,
if no confusion is possible.
We denote by $E^{\ast}$ the dual space of $E$,
and by $\left\langle \cdot,\cdot\right\rangle_{E^*,E}$
the duality between $E$ and $E^*$.

If $F$ is another Banach space, $\call(E,F)$ denotes the
space of bounded linear operators from $E$ to $F$ endowed with the usual operator norm.
All Banach spaces are assumed to be real and separable.

In what follows we will often meet inverses of operators which are not
one-to-one. Let $Q\in \call\left(H,K\right)$.
Then $H_{0}=\ker Q$ is a closed subspace of $H$. Let
$H_{1}:=[\ker Q]^\perp$ be the orthogonal complement of $H_{0}$ in $H$: $H_{1}$ is closed, too.
Denote by $Q_{1}$ the restriction of $Q$ to $H_{1}$: $Q_{1}$ is
one-to-one and $\operatorname{Im}Q_{1}=\operatorname{Im}Q$.
For $k\in \operatorname{Im}Q$, we define $Q^{-1}$ by setting
\(
Q^{-1}\left(k\right)  :=Q_{1}^{-1}\left(k\right)
\).
The operator $Q^{-1}:\operatorname{Im}Q\rightarrow H$ is called the pseudoinverse of $Q$. $Q^{-1}$ is linear and closed but in general not continuous.
Note that if $k\in\operatorname{Im}Q$, then
$Q_{1}^{-1}\left(  k\right)\in[\ker Q]^\perp$.
is the unique element of
\(
\left\{ h  :Q\left(  h\right)  =k\right\}
\)
with minimal norm (see e.g. \cite{Z}, p.209),

Next we introduce some spaces of functions.
Let $H$ and $Z$ be real separable Banach spaces.
By $B_b(H,Z)$ (respectively $C_b(H,Z)$, $UC_b(H,Z)$) we denote the space of all functions
$f:H\rightarrow Z$ which are Borel measurable and bounded (respectively continuous
and bounded, uniformly continuous and bounded).

Given an interval $I\subseteq \R$ we denote by
$C(I\times H,Z)$ (respectively $C_b(I\times H,Z)$)
the space of all functions $f:I \times H\rightarrow Z$
which are continuous (respectively continuous and bounded).
$C^{0,1}(I\times H,Z)$ is the space of functions
$ f\in C(I\times H)$ such that, for all $t\in I$,
$f(t,\cdot)$ is continuously Fr\'echet differentiable {with $\nabla f $ continuous}.
By $UC_{b}^{1,2}(I\times H,Z)$
we denote the linear space of the mappings $f:I\times H \to Z$
which are uniformly continuous and bounded
together with their first time derivative $f_t$ and their first and second space
derivatives $\nabla f,\nabla^2f$.

Finally, if the arrival space $Z=\R$ we do not write it in all the above spaces.




\subsection{$C$-derivatives}\label{subsection-C-directionalderivatives}

In this subsection we recall the definition of $C$-directional derivatives, see \cite{FGFM-III} where $C$ is a suitable linear operator possibly unbounded. We remark that these generalized derivatives has been introduced to cover a wider class than the one considered in \cite{FedericoGozziJDE},\cite[Section 4.2.1]{FabbriGozziSwiech}, (which was already an extension of the definition of directional derivatives given in \cite[Section 2]{Mas}, \cite{FTGgrad} in the case when $C$ is a bounded operator, see also \cite{FGFM-I,FGFM-II}).

We recall the definition of $C$-derivative,
see \cite[Definition 4.1]{FGFM-III}: the operator $C$ is unbounded in the sense that it does not take its values in the state space $H$ but in a larger Banach space $\overline H$ such that $H\subset \overline H$ with continuous embedding. We send the reader to
\cite[Remarks 4.2-4.3-4.4]{FGFM-III} for more explanations and comparison with other definitions.

\begin{definition}
\label{df4:Gderunbounded} Let $H, \,Z,\, K$ and $\overline H$ be
Banach spaces such that $H\subset \overline H$ with continuous embedding.
Let $C:K\rightarrow \overline H$ be a linear and bounded operator.
\begin{itemize}
\item[(i)]
Let $k\in K$, let $h=Ck$ and let $f:\overline H\rightarrow Z$.
We say that $f$ admits $C$-directional derivative
at a point $x\in \overline{H}$ in the direction $k\in K$
(and we denote it by $\nabla^{C}f(x;k)$) if the limit,
in the norm topology of $Z$,
\vspace{-0.3truecm}
\begin{equation}\label{Cderivatabis}
 \nabla^{C}f(x;k):=\lim_{s\rightarrow 0}
\frac{f(x+s Ck)-f(x)}{s},
\vspace{-0.3truecm}
\end{equation}
exists.
\item[(ii)]
Let $f:\overline H\rightarrow Z$.
We say that $f$ is $C$-G\^ateaux differentiable
at a point $x\in \overline{H}$ if $f$ admits the $C$-directional derivative in every
direction $k\in K$ and there exists a {\bf bounded} linear operator,
the $C$-G\^ateaux derivative $\nabla^C f(x)\in L(K,Z)$, such that $\nabla^{C}f(x;k)  =\nabla^{C}f(x)k$
for all $k \in K$. We say that $f$ is $C$-G\^ateaux
differentiable on $H$ (respectively $\overline{H}$) if it is $C$-G\^ateaux differentiable at every point $x\in
H$ (respectively $x\in\overline{H}$).
\item[(iii)]
Let $f:\overline H\rightarrow Z$.
We say that $f$ is $C$-Fr\'echet differentiable
at a point $x\in \overline{H}$ if it is $C$-G\^ateaux differentiable and if the limit in (\ref{Cderivatabis}) is uniform for $k$ in the unit ball of $K$. In this case
we call $\nabla^C f(x)$ the $C$-Fr\'echet derivative (or simply the $C$-derivative) of $f$ at $x$. We say that $f$ is $C$-Fr\'echet differentiable on $H$ (respectively $\overline{H}$) if it is $C$-Fr\'echet differentiable at every point $x\in H$ (respectively $x\in\overline{H}$).
\end{itemize}
\end{definition}

We now give the definition of suitable spaces of $C$-differentiable functions, which have been introduced in \cite[Definition 4.5]{FGFM-III}. In the following for a function $f:\overline H\rightarrow Z$ we denote by $f\!\!\mid_{H}$ the restriction of $f$ to $H$. Similarly, when $I$ is an interval in $\R$, for $f:I\times \overline H\rightarrow Z$ we denote by $f\!\!\mid_{I\times H}$ the restriction of $f$ to $I\times H$.

\begin{definition}
\label{df4:Gspaces}
Let $I$ be an interval in $\R$, let $H$, $\overline{H}$, $K$ and $Z$ be suitable real Banach spaces. Moreover let $H\subset \overline{H}$ with continuous inclusion,
and let $C\in \call(K,\overline{H})$.
\begin{itemize}
\item
We call $C^{1,C}_{b}(\overline{H},Z)$ the space of all continuous and bounded functions $f:\overline H\to Z$ which admit continuous and bounded $C$-Fr\'echet derivative.
Moreover we call $C^{0,1,C}_b(I\times \overline{H},Z)$ the space of
continuous and bounded functions $f:I\times \overline H\to Z$ such that, for every $t\in I$, $f(t,\cdot)\in C^{1,C}_b(\overline{H},Z)$ and
$\nabla^C f\in C_b\left(I\times \overline{H},L(K,Z)\right)$
When $Z=\R$ we write  $C^{1,C}_{b}(\overline{H})$ instead of $C^{1,C}_{b}(\overline{H},Z)$, and it turns out that if $f\in C^{1,C}_{b}(\overline{H})$, then  $\nabla^C f\in C_b\left(I\times \overline{H},K)\right)$.
\item
For any $\alpha\in{[0,1]}$ and $T>T_0\geq0$ (this time $I$ is equal to $[T_0,T]$) we denote by
$C^{0,1,C}_{\alpha}([T_0,T]\times \overline{H},Z)$
the space of functions
$f\in C_b([T_0,T]\times H,Z)\cap
C^{0,1,C}_b((T_0,T]\times \overline{H},Z)$\footnote{Note that here
$f(t,\cdot)$ is well defined only in $H$ when $t=T_0$, while for $t>T_0$ it is defined over $\overline{H}$. The reason is that the Ornstein-Uhlenbeck semigroup in our examples and in our setting (and consequently the solution of the HJB equation) satisfy the same property.}
such that
the map $(t,x)\mapsto t^{\alpha} \nabla^C f(t,x)$
belongs to $C_b((T_0,T]\times \overline{H},\call(K,Z))$.
When $Z=\R$ we write  $C^{0,1,C}_{\alpha}([T_0,T]\times \overline{H})$ instead of $C^{0,1,C}_{\alpha}([T_0,T]\times \overline{H},Z)$, and it turns out that, if
$f\in C^{0,1,C}_{\alpha}([T_0,T]\times \overline{H})$, then the map $(t,x)\mapsto t^{\alpha} \nabla^C f(t,x)$ belongs to
$C_b\left((T_0,T]\times \overline{H},K)\right)$.
The space $C^{0,1,C}_{\alpha}([T_0,T]\times \overline{H},Z)$
is a Banach space when endowed with the norm
\vspace{-0.3truecm}
\[
 \left\Vert f\right\Vert _{C^{0,1,C}_{\alpha}([T_0,T]\times \overline{H})  }=\sup_{(t,x)\in (T_0,T]\times \overline{H}}
\vert f(t,x)\vert+
\sup_{(t,x)\in (T_0,T]\times \overline{H}}  t^{\alpha }\left\Vert \nabla^C f(t,x)\right\Vert_{\call(K,Z)}.
\vspace{-0.2truecm}
\]
When clear from the context we will write simply
$\left\Vert f\right\Vert _{C^{0,1,C}_{\alpha}}$.


\end{itemize}
\end{definition}


%


{\noindent\textbf{Acknowledgements}: The authors would like to thank the anonymous referees for valuable remarks which led to an improved version of the paper.
}

\end{document}